\documentclass[reprint,twocolumn,10pt,amsmath,amssymb,aps,pre,superscriptaddress,floatfix]{revtex4-2}
\usepackage[section]{placeins}
\usepackage[T1]{fontenc}
\usepackage[utf8]{inputenc}
\usepackage[english]{babel}
\usepackage[margin=1.5cm]{geometry}
\usepackage{fancyhdr}
\usepackage{soul}
\usepackage{graphicx}
\usepackage{caption}
\usepackage{subcaption}
\usepackage{amsmath, amssymb, amsfonts, amsthm, stmaryrd}
\usepackage{enumitem}
\usepackage{tikz}
\usepackage{indentfirst}
\usetikzlibrary{positioning, matrix, arrows.meta}
\usepackage{placeins}
\usepackage{csquotes}
\usepackage[toc,page,title]{appendix}
\usepackage{csquotes}
\usepackage[most]{tcolorbox}
\usepackage{xcolor}
\usepackage{enumitem}
\usepackage{setspace}

\usetikzlibrary{positioning, matrix, arrows.meta}
\usepackage{tabularx}   % tabularx environment + X column
\usepackage{booktabs}   % \toprule \midrule \bottomrule
\usepackage[toc,page,title]{appendix}
\usepackage[table]{xcolor} % \rowcolor, \rowcolors, colored cells/rows

\usepackage{array}

\newtheorem{theorem}{Theorem}

\newtheorem{definition}[theorem]{Definition}
\newtheorem{corollary}[theorem]{Corollary}

\newtheorem{proposition}[theorem]{Proposition}
\definecolor{blueheader}{HTML}{DCEEFB}
\definecolor{beigeback}{HTML}{FDF6EC}

\newtcbtheorem{algo}{Algorithm}{
  colback=beigeback,
  colframe=blueheader,
  coltitle=black,
  fonttitle=\bfseries,
  boxrule=0.8pt,
  arc=2mm,
  left=2mm,
  right=2mm,
  top=1mm,
  bottom=1mm
}{alg}
\begin{document}
\setlength{\parindent}{1em}

\title{Mean First Passage Times for Loop-Erased Random Walks and Erased-Loop Size Distributions on Graphs}

\author{Komlan Fiagbe}
\email{komlan.fiagbe@unamur.be}
\affiliation{Department of Mathematics \& Namur Institute for Complex Systems - naXys, University of Namur, Rue Joseph Grafé 2, 5000 Namur, Belgium}

\author{Jean-François de Kemmeter}
\affiliation{Department of Mathematics, Royal Military Academy, Brussels, Belgium.}

\author{Timoteo Carletti}
% \email{timoteo.carletti@unamur.be}
\affiliation{Department of Mathematics \& Namur Institute for Complex Systems - naXys, University of Namur, Rue Joseph Grafé 2, 5000 Namur, Belgium}
 \begin{abstract}
Loop-erased random walks provide a natural way to extract self-avoiding
trajectories from ordinary random walks by chronologically removing loops. They
are deeply connected with uniform spanning trees: the loop-erased path between
two vertices has the same distribution law as the unique path connecting these vertices in a
random spanning tree. In this work, we use this connection as a computational
framework for finite graphs. First, using the classical transfer-current theorem, we obtain an explicit formula for the
mean first passage time of loop-erased random walks. The resulting expression
returns the expected length of the loop-erased path as a weighted sum over simple
paths, where the weights are exact determinantal probabilities obtained from the
spanning-tree kernel. Second, we express the
distribution of erased loop sizes in the same framework. This distribution measures the scale of the
redundant excursions removed during loop erasure and provides information that is
complementary to effective path length. The approach is applied to arbitrary
finite graphs and then specialized to complete graphs, cycle graphs, and
motif-based networks, where explicit formulas and asymptotic behaviors are
obtained. Numerical experiments show strong agreement between the theoretical
predictions and Monte Carlo simulations. The results provide a unified method for
quantifying both the effective paths kept by loop erasure and the loops removed
during random exploration of networks.
\end{abstract}   

\maketitle
\section{Introduction}
Random walks are a central tool in the analysis of networks. They provide mathematical models for diffusion~\cite{masuda2017random}, community detection~\cite{lambiotte2014random}, search~\cite{dipatti2015optimal} and ranking nodes~\cite{page1999pagerank}. Networks exploration can be obtained via metrics such as hitting times~\cite{wang2008first,sun2022hitting,zhang2011mean,lee2014estimating}, commute times, cover times, and effective resistances~\cite{DoyleSnell2000_RandomWalksElectricNetworks,Grimmett2012_ProbabilityOnGraphsNotes,LyonsPeres2016ProbabilityTreesNetworks}. These quantities describe how a walker moves through a graph and how the topology of the graph affects accessibility between vertices. In many situations, however, the full random-walk trajectory contains many redundant steps. A walker may return several times to already visited vertices before reaching a target. If the goal is to extract the effective route followed by the walker, it is natural to erase these redundant excursions. This is the principle of the {\em Loop-Erased Random Walk} (LERW): starting from an ordinary random-walk path, loops are removed in chronological order, leaving a self-avoiding path~\cite{Lawler1999LERW,Barlow2014LoopErasedWalksUST}. LERW were introduced and developed as fundamental objects in probability theory. Lawler's work established LERW as a canonical method for producing self-avoiding paths from random walks~\cite{Lawler1999LERW}. The importance of LERW increased further through its connection with uniform spanning trees. In particular, the loop-erased path between two vertices has the same distribution as the unique path connecting these vertices in a uniform spanning tree. This connection is the basis of Wilson's algorithm, which generates uniform spanning trees efficiently by using loop-erased random walks~\cite{wilson1996generating}. The relationship between LERW, uniform spanning trees, and spanning forests has since become a major theme in probability on graphs and networks~\cite{Pemantle2004Choosing,BenjaminiLyonsPeresSchramm2001_UniformSpanningForests,LyonsPeres2016ProbabilityTreesNetworks}.
A large part of the theoretical literature on LERW has focused on infinite graphs, lattices, and scaling limits. In the planar case, Lawler, Schramm, and Werner proved the conformal invariance of loop-erased random walks and uniform spanning trees, connecting these objects to Schramm--Loewner evolution~\cite{LawlerSchrammWerner2004ConformalInvarianceLERWUST}. Other works have studied LERW and spanning trees on lattices and more general graph structures, emphasizing their probabilistic and geometric properties ~\cite{Barlow2014LoopErasedWalksUST,LyonsPeres2016ProbabilityTreesNetworks}. These studies provide a deep theoretical foundation, but many questions remain when one moves from regular lattices to arbitrary finite graphs and 
% $applied$ 
empirical network models.
The spanning-tree point of view is particularly useful because it turns questions about random trajectories into questions about paths in random trees. Edge-inclusion probabilities in random spanning trees can be described by determinantal formulas involving transfer-current kernels~\cite{Lyons2003DeterminantalProbabilityMeasures,KasselWu2015_TransferCurrentPatternFields}. This gives exact probability weights for events involving sets of edges, rather than relying only on empirical sampling. In the present work, we exploit this determinantal structure to compute metrics associated with loop-erased random walks on finite graphs. In particular, the probability of a simple path appearing as the loop-erased trajectory can be expressed through the determinant of the corresponding transfer-current submatrix.

The aim of this work is to use this framework to study loop-erased random walks on finite graphs from a quantitative and applicative perspective. We focus on two complementary observables. The first one is the {\em Mean First Passage Time} (MFPT) of the loop-erased random walk. Unlike the usual first passage time, which counts all steps taken by the original random walk before reaching a target, the loop-erased version measures the length of the self-avoiding path that remains after all loops have been removed. It therefore gives an effective exploration length between two vertices. The second observable is the {\em distribution of erased loop sizes}. During the loop-erasure procedure, each return to an already visited vertex creates a loop that is removed from the trajectory. The size of this erased loop measures the scale of the redundant excursion. This distribution gives information that is different from the length of the final loop-erased path: it describes what the random walk explored and then discarded.
Previous works have studied erased-loop statistics mainly in lattice and low-dimensional settings. Dhar and Dhar investigated the distribution of erased-loop sizes for loop-erased random walks, and Agrawal and Dhar extended this study to two- and three-dimensional systems~\cite{DharDhar1997ErasedLoops,AgrawalDhar2001ErasedLoops}. These works show that erased loops contain important information about the geometry of the underlying space and about the redundant parts of random-walk exploration. In contrast, the present study considers erased-loop statistics on arbitrary finite graphs and relates them to spanning-tree probabilities. This provides a finite-graph framework that is suitable for network applications and allows both theoretical computation and numerical validation.
The complete graph case connects the study of LERW to classical results on random labeled trees. The distance between two vertices in a uniformly random labeled tree was studied by Meir and Moon, and related enumerative results are collected in Moon's work on labeled trees~\cite{MeirMoon1970Distance,Moon1970CountingLabelledTrees}. These results allow explicit formulas and asymptotic expressions for the mean length of loop-erased paths on complete graphs. Cycle graphs, on the other hand, provide a simple setting in which the role of global cyclic structure can be isolated using elementary graph-theoretic arguments~\cite{West1996IntroductionGraphTheory}. Together, these graph classes serve as reference models for understanding how density and cyclic structure affect LERW observables.

From an applied network perspective, loop-erased quantities provide a complementary description of graph exploration. Shortest-path distances ignore the random exploration process, while ordinary hitting times include all repeated visits and local oscillations. Loop-erased paths lie between these two viewpoints: they are generated by random exploration, but they retain only the self-avoiding effective route. The erased-loop distribution adds further information by quantifying the scale of the redundant cyclic excursions produced by the graph topology. This is particularly relevant for networks with short loops, motifs, and local subgraph structure. Network motifs were introduced as simple building blocks of complex networks~\cite{Milo2002R}, and later models incorporated arbitrary distributions of subgraphs or short loops into random graph constructions~\cite{KarrerNewman2010RandomGraphsSubgraphs,Newman2019SpectraShortLoops}. Such structures can strongly influence the loops encountered and erased by a random walk.
For computational purposes, the present approach combines exact spanning-tree probabilities with graph-search procedures. The enumeration of candidate paths and cycles is related to classical depth-first search methods and algorithms for listing paths and cycles in graphs~\cite{Tarjanrober1972,sedgewick2001algorithms,ferreira2012optimal}. Although the number of simple paths can grow rapidly with graph size, the proposed method is effective for sparse networks, small graphs, and motif-based structures where the relevant paths can be enumerated or sampled. In this sense, the method differs from a purely Monte Carlo treatment: simulations may be used to generate candidate paths or to validate the results, but the path weights themselves are obtained from the transfer-current kernel.

A possible application-oriented interpretation concerns agents exploring an initially
unknown environment. A foraging animal or an autonomous robot may reach a target
after repeated visits to the same locations. If its visited route has been
recorded, chronological loop erasure provides a self-avoiding path that can
be followed in reverse towards the starting point, avoiding the deleted
excursions. This interpretation is conceptual: loop erasure neither identifies
a globally shortest path nor establishes that biological navigators implement
this algorithm. It requires sufficient route memory and reliable retracing
of previously traversed edges. The retained-path length and erased-loop
statistics thus quantify different aspects of such exploration-derived
navigation.

The contribution of this paper is therefore twofold. Building on the known connection between LERW and spanning trees and on the transfer-current theorem, we express the mean first passage time of loop-erased random walks on finite graphs. The resulting expression returns the expected length of the loop-erased path as a weighted sum over simple paths,
where the weights are exact determinantal probabilities obtained from the spanning-tree kernel. Second, We use the established Dhar--Dhar correspondence~\cite{DharDhar1997ErasedLoops} to express erased-loop size statistics through the same kernel. Our finite-graph treatment includes a sampled-path estimator with a captured-mass error bound, explicit specializations to complete and cycle graphs, and decompositions for motif-based networks. The mean loop-erased path length describes the part of the random-walk trajectory that is kept, while the erased-loop size distribution describes the loops removed. Numerical comparisons with Monte Carlo simulations are presented throughout the paper to support the theoretical findings.

 \section{MFPT of LERW on an arbitrary graph}
 \label{sec:MFPTLERW}
Let \(G=(V,E)\) be a finite, connected, simple, undirected graph, and let
\(u,v\in V\) be two distinct vertices. Consider a random walk
\((X_t)_{t\geq 0}\) on \(G\), started at \(u\), and stopped upon first
hitting \(v\). Its hitting time is defined by
\[
    \tau_v := \inf\{t\geq 0 : X_t = v\}.
\]
The stopped random walk trajectory is then the, possibly self-intersecting,
path
\[
    \omega = (X_0,X_1,\dots,X_{\tau_v}).
\]
Applying chronological loop erasure to \(\omega\), that is, erasing loops in
the order in which they are formed, yields a self-avoiding path from \(u\) to
\(v\), denoted by
\(
    \mathrm{LE}(\omega).
\)
This path is called the loop-erased random walk, or LERW, from \(u\) to \(v\)
~\cite{Lawler1999LERW}.
In the unweighted case, the passage time of the LERW is simply its length,
namely the number of edges in the loop-erased path
\[
    T_{u,v}^{\mathrm{LERW}}
    := \bigl|\mathrm{LE}(\omega)\bigr|_E ,
\]
where $\bigl|\mathrm{A}\bigr|_E$ denotes the number of elements of $E$ belonging to $A\subset E$.

The same definition is kept even if \(G\) is endowed with edge weights \(c:E\to(0,\infty)\).
The mean first passage time for the loop-erased random walk is therefore
\begin{equation*}
\mathrm{MFPT}^{\mathrm{LERW}}(u,v)
:=
\mathbb{E}\left[T_{u, v}^{\mathrm{LERW}}\right]\, ,
\end{equation*}
and by introducing the probability, $\mathbb{P}\!\left(\mathrm{LE}(\omega)=\gamma\right)$, that the loop erased walk of $\omega$ is given by $\gamma\in \Gamma(u,v)$, the set of all simple paths from $u$ to $v$, we can obtain
\begin{equation}
\label{eq:deMFTP}
    \mathrm{MFPT}^{\mathrm{LERW}}(u,v)
=
\sum_{\gamma\in\Gamma(u,v)} |\gamma|_E \mathbb{P}\left(\mathrm{LE}(\omega)=\gamma\right)\, .
\end{equation}
One can introduce the matrix $Y$ indexed by (oriented) edges of the graph, called the ``transfer-current'' or ``spanning-tree kernel'', with the following determinantal property ~\cite{Lyons2003DeterminantalProbabilityMeasures,LyonsPeres2016ProbabilityTreesNetworks,KasselWu2015_TransferCurrentPatternFields}: for any finite set of (distinct) edges $S\subseteq E$,
\begin{align*}
\mathbb{P}(S \subseteq T) \;=\; \det\!\big(Y_{S,S}\big),
\end{align*}
where $T$ is a random spanning tree drawn from the spanning-tree measure
(uniform in the unweighted case, and weighted by $\prod_{e\in T} c_e$ in the conductance-weighted case),
and $Y_{S,S}$ denotes the submatrix of $Y$ restricted to the edges in $S$. 
A standard and convenient construction of such a kernel is obtained as follows:
choose an arbitrary orientation of edges and let $B$ be the incidence matrix
\begin{equation}
    \label{eq:Binc}
B(i,\ell)=\begin{cases}
    -1 & \text{if $\ell=[i,*]$}\\
     1 & \text{if $\ell=[*,i]$}\\
     0 & \text{otherwise}
\end{cases}
\end{equation}
where $\ell=[i,j]$ denotes the oriented edge with tail $i$ and head $j$. Let $C=\mathrm{diag}(c_e)$, i.e., the diagonal matrix with the edge weights on the diagonal, and let 
% \textcolor{red}{[If not wrong $B$ has dimension $\vert V \vert \times \vert E \vert$, with $\vert V \vert$ the number of nodes and $\vert E \vert$ the number of edges. In this case, shouldn't be the weighted Laplace matrix be given by $L = B C B^T$ (so that the dimensions agree) ?]}
\[
L := BCB^\top 
\]
be the weighted Laplace matrix. By using the Moore-Penrose pseudoinverse $L^{\dagger}$, we can define
\begin{align}\label{SEC1EQ1}
 Y :=\; C^{1/2} B^\top L^{\dagger} B C^{1/2}\, .  
\end{align}
% \textcolor{red}{[Same remark as before regarding the dimensions of the matrix.]}
Consider the unique simple path between $u$ and $v$ inside a spanning tree $T$, say $\mathrm{Path}_{T}(u,v)$. In a tree, the event ``$\mathrm{Path}_{T}(u,v)=\gamma$'' is exactly the event
``all edges of $\gamma$ are present'':
\[
\{\mathrm{Path}_{T}(u,v)=\gamma\}
\;=\;
\{E(\gamma)\subseteq T\},
\]
because $\gamma$ already connects $u$ to $v$, and a tree cannot contain an alternative $u\text{--}v$ route without creating a cycle. The crucial link is that the LERW path from $u$ to $v$ has the same distribution law than the $u\text{--}v$ path in $T$
under the spanning-tree measure produced by loop-erasing random-walk segments ~\cite{LawlerSchrammWerner2004ConformalInvarianceLERWUST,Barlow2014LoopErasedWalksUST,Pemantle2004Choosing}.
Therefore, for each simple path $\gamma$, it holds true
\begin{align}\label{SEC1EQ2}
 \mathbb{P}\!\left(\mathrm{LE}(\omega)=\gamma\right)
&=
\mathbb{P}\!\left(\mathrm{Path}_{T}(u,v)=\gamma\right)\\
&=
\mathbb{P}\!\left(E(\gamma)\subseteq T\right)\notag\\
&=
\det\!\big(Y_{E(\gamma),E(\gamma)}\big).
\end{align}
Let $E(\gamma)$ be the set of edges of the simple path $\gamma$ and define the weight
\begin{align}\label{SEC1EQ3}
 p_\gamma :=\det\big(Y_{E(\gamma),E(\gamma)}\big)\, .
\end{align}
Then the MFPT of the LERW~\eqref{eq:deMFTP} can be rewritten as
\begin{equation}
\label{eq:defMFPT2}
  \mathrm{MFPT}^{\mathrm{LERW}}(u,v)
=
\sum\limits_{\gamma\in\Gamma(u,v)} p_\gamma\,|\gamma|_E\, .  
\end{equation}\begin{algo}{$\widehat{\mathrm{MFPT}}^{\mathrm{LERW}}(u,v)$}{1}
\textbf{1. Kernel precomputation.}\\
Form the weighted Laplacian \(L := BCB^\top \), compute its Moore--Penrose pseudoinverse \(L^\dagger\), and define the edge kernel
\[
Y := C^{1/2} B^\top L^\dagger B C^{1/2}.
\]

\vspace{0.2cm}

\textbf{2. Candidate LERW paths.}\\
For \(i=1,\dots,N\): simulate a random walk from \(u\) until the first hitting time of \(v\); apply loop erasure to obtain a simple path \(\gamma\).
Store each distinct \(\gamma\) in a set \(\mathcal{G}\).

\vspace{0.2cm}

\textbf{3. Determinantal weights.}\\
For each \(\gamma\in\mathcal{G}\): extract its edge set \(E(\gamma)\), form the principal submatrix \(Y_{E(\gamma),E(\gamma)}\), and set
\[
p_\gamma \;=\; \det\!\big(Y_{E(\gamma),E(\gamma)}\big).
\]

\vspace{0.2cm}

\textbf{4. MFPT estimate.}\\
Return
\[
\widehat{\mathrm{MFPT}}^{\mathrm{LERW}}(u,v)
\;=\;
\frac{\sum_{\gamma\in\mathcal{G}} p_\gamma\,|\gamma|_E}{\sum_{\gamma\in\mathcal{G}} p_\gamma}.
\]

\end{algo}

The estimator returned by Algorithm~\ref{alg:1}
is not necessarily unbiased for a finite
number of sampled trajectories, since
undiscovered paths are excluded from the
normalized weighted average.
However, the cumulative probability mass
$M_N=\sum_{\gamma\in\mathcal{G}_N}p_\gamma$
provides a direct measure of the
completeness of the discovered path set.
In Appendix~\ref{APP_MFPT_ERROR}, we establish
the deterministic error bound
\begin{align*}
\left|
\mathrm{MFPT}^{\mathrm{LERW}}(u,v)
-
\widehat{\mathrm{MFPT}}_N^{\mathrm{LERW}}(u,v)
\right|
\leq\\
\leq (1-M_N)
\left(
|V|-1-d_G(u,v)
\right)    
\end{align*}
where $d_G(u,v)$ is the length of the shortest path between $u$ and $v$ in $G$.
We also derive an interval containing
the exact MFPT and prove that the
estimator converges almost surely to
the exact value as $N\to\infty$.
Let us observe that this procedure does not estimate the MFPT by taking the empirical average of the sampled LERW lengths. Instead, sampling is used only to discover candidate simple paths \(\gamma\).
Each discovered path is then assigned its theoretical probability weight \(p_\gamma\) through the determinant of the kernel \(Y\).
During Step~2, one may also monitor the cumulative mass \(\sum_{\gamma\in\mathcal{G}} p_\gamma\), which provides an estimate of the fraction of total probability already captured by the currently discovered paths.
Furthermore, the LERW-based path generation can be replaced by an explicit enumeration of simple paths, for instance using a depth-first search (DFS) strategy ~\cite{sedgewick2001algorithms}. We know that the number of paths can grow exponentially as the size of the graph increases; however this method performs very well on large sparse networks, where it is possible to efficiently enumerate the paths between two nodes, as well as on small networks.
\begin{figure}
    \centering
    \includegraphics[width=1\linewidth]{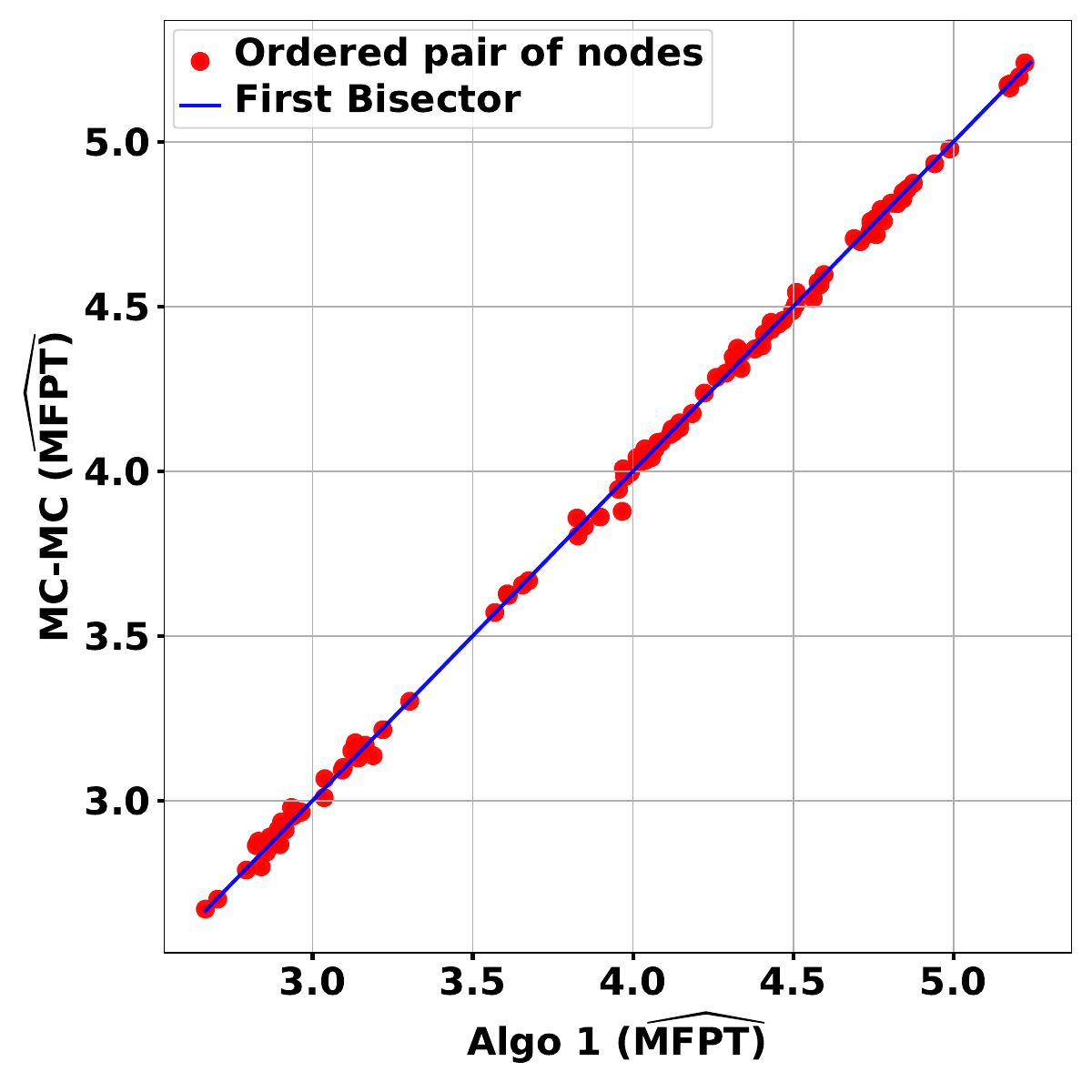}
    \caption{Mean First Passage Time of Loop Erased Random Walks between the different pairs of nodes of a Watts-Strogatz graph ~\cite{Article32} with $15$ nodes, a rewiring probability of $0.1$, and a mean degree of $4$. The number of samples in the Monte Carlo approach is $40000$.  On the x-axis, we have the MFPT obtained with the algorithm \ref{alg:1} and on the y-axis, we have the MFPT obtained with an MC-MC approach with empirical average.}
    \label{KERNEL}
\end{figure}
For the numerical comparisons against the first bisector, we report the maximum absolute deviation
\begin{equation}\label{eq:max_abs_dev}
\varepsilon_{\max}=\max_i |y_i-x_i|,
\end{equation}
where $x_i$ and $y_i$ are the coordinates of the plotted points.

\begin{figure}
    \centering
    \includegraphics[width=1\linewidth]{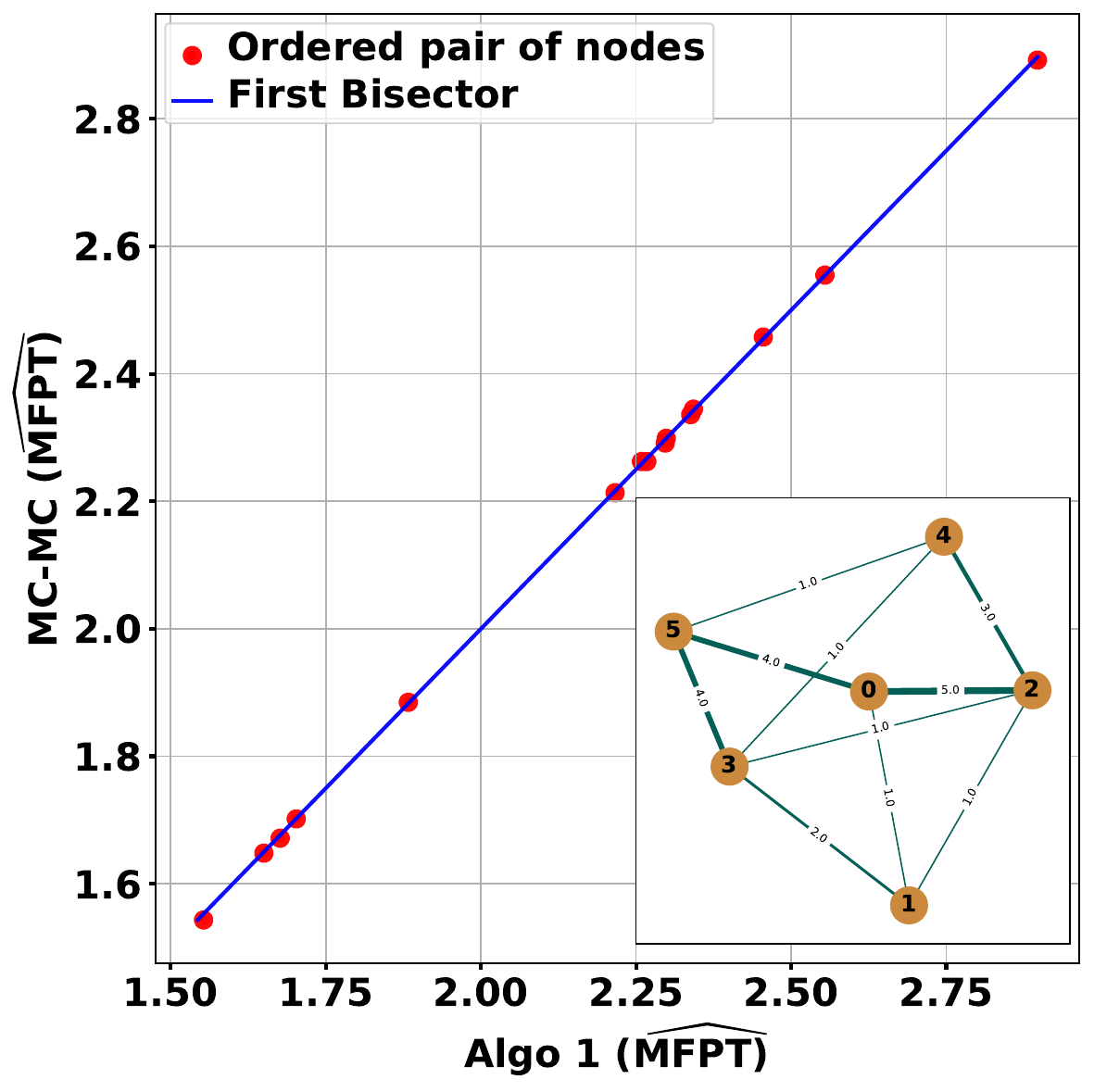}
    \caption{Main Panel shows the Mean First Passage Time of Loop Erased Random Walks between the different pairs of nodes of a weighted graph. The number of samples in the Monte Carlo approach is $40000$. On the x-axis, we have the MFPT obtained with the algorithm~\ref{alg:1} and on the y-axis, we have the MFPT obtained with an MC-MC approach with empirical average. Inset shows the example of weighted graph we used. Weights are reported on each link whose thickness increases as the weight gets larger.}
    \label{WEIGHT}
\end{figure}

We apply this method to a Watts-Strogatz graph~~\cite{Article32} with $15$ nodes, a rewiring probability of $0.1$, and a mean degree of $4$. The results obtained by using the algorithm~\ref{alg:1} are shown in Fig.~\ref{KERNEL} and compared with those obtained with an MC-MC approach relying on the computation of empirical averages. The two approaches agree closely, with $\varepsilon_{\max}=3.6\times 10^{-2}$.
Another application is on a toy example of weighted graph shown in the inset of Fig.~\ref{WEIGHT}. Let us recall that in the weighted graph, $ C^{1/2} \neq I$, the edge conductances modify the random-walk transition probabilities and hence the distribution of loop-erased paths. If passage cost is of interest, an independent traversal-cost function $t_e$ must be specified and inserted in the path sum. The results are reported in the main panel of Fig.~\ref{WEIGHT} ; the largest absolute discrepancy between the two estimates is $\varepsilon_{\max}=9.6\times10^{-2}$.
The computational cost depends on Laplacian pseudoinverse computation and, most importantly, the number of distinct simple paths examined; spanning trees need not be explicitly enumerated. We have to enumerate paths between nodes and this can require exponential time in the worst case, since the number of such paths may itself be exponential. This approach is therefore especially effective for small graphs or sparse networks. Having quantified the effective part of the random walk that remains after loop erasure, we now turn to the complementary question: what is removed during this procedure? In the next section, we study the distribution of erased loop sizes and show that it can be described using the same spanning-tree and transfer-current framework.

\section{Distribution of erased loops sizes of LERW on an arbitrary graph}
\label{sec:distsize}
Let $G=(V,E)$ be a finite connected, simple, undirected, unweighted graph.

\begin{definition}[Spanning Forest]
A \emph{spanning forest} of $G$ is a subgraph $F = (V,E_F)$ such that $E_F \subseteq E$, $F$ is acyclic.
\end{definition}

\begin{definition}[$k$-Component Spanning Forest and spanning tree]
A spanning forest $F$ of $G$ is called a $k$-component spanning forest if $F$ has exactly $k$ connected components, counting isolated vertices as single-vertex components. A spanning tree of $G$ is a $1$-component spanning forest of $G$.
We denote by
\[
\mathcal{T}(G) := \{ \text{spanning trees of } G \}
\]
the set of all spanning trees of $G$.

\end{definition}
\begin{definition}[$k$-Excess Subgraph and spanning unicycle]
A subgraph $H = (V,E_H)$ of $G$ is called a $k$-excess subgraph ~\cite{SunWilson2015SandpilesUnicycles}
if there exists a spanning tree $T = (V,E_T)$ of $G$ such that
\[
E_H = E_T \cup X,
\]
where $X \subseteq E \setminus E_T$ and $|X| = k$. A spanning unicycle is a $1$-excess subgraph of $G$. We denote by
\[
\mathcal{U}_k(G) := \{ \text{$k$-excess subgraphs of } G \}
\]
the set of all $k$-excess subgraphs of $G$.
\end{definition}

Fix a vertex $o\in V$.
Let $(X_n)_{n\ge0}$ be a random walk on $G$ started at $o$.
For $N\ge0$, denote by
\[
\omega_N=(X_0,X_1,\dots,X_N)
\]
the random walk path up to time $N$. As already observed, given a finite path $\omega=(v_0,\dots,v_N)$, its
loop-erasure, denoted $\mathrm{LE}(\omega)$,
is obtained by chronologically removing loops:
whenever $v_i=v_j$ with $i<j$, delete the segment
$(v_{i+1},\dots,v_j)$ and continue. 
The resulting path is self-avoiding. Denote by
\[
\eta_N := \mathrm{LE}(\omega_N)
\]
the loop-erased random walk (LERW) up to time $N$. 
The LERW determines a spanning tree structure in the following sense.
If we let $N$ to grow and stop the process when the walk first visits all vertices,
the loop-erased path construction induces a spanning tree
(see Wilson's algorithm) ~\cite{wilson1996generating}. Thus LERW and spanning trees are
probabilistically equivalent objects on finite graphs.
Fix $\ell\ge 2$.
We study the probability that the $(N+1)$-th step of the random walk
creates a loop of size $\ell$ that is erased by the loop-erasure.

\begin{definition}[Loop of size $\ell$]
 Let \(G=(V,E)\) be a finite, connected, simple, undirected graph, and let \((X_n)_{n\geq0}\) be a random walk on \(G\), started at a vertex \(o\in V\). For \(N\geq0\), denote by
\[
\omega_N=(X_0,X_1,\ldots,X_N)
\]
the random-walk trajectory up to time \(N\), and by
\[
\eta_N=\operatorname{LE}(\omega_N)
=(Y_0^{(N)},\ldots,Y_{m_N}^{(N)})
\]
its chronological loop erasure.
At step \(N+1\), a loop of size \(\ell\) is formed and erased if the new vertex \(X_{N+1}\) already belongs to \(\eta_N\), say
\[
X_{N+1}=Y_j^{(N)}
\]
for some \(0\leq j\leq m_N\). Its size is
\[
\ell=m_N-j+1.
\]
On a simple graph, \(\ell\geq2\).
\end{definition}
Let \(P_N(\ell)\) denote the probability that step \(N+1\) creates and erases a loop of size \(\ell\). We define the long-time erased-loop probability by
\[
P(\ell)=\lim_{N\to\infty}P_N(\ell),
\]
whenever this limit exists.
This measures the expected frequency per random-walk step
of erased loops of size $\ell$. The correspondence below
is understood with the long-time convention of
\cite{DharDhar1997ErasedLoops}.
% Let \(T\) be a uniformly selected spanning tree of \(G\), and let
% \(e=\{u,v\}\) be an edge selected uniformly from \(E(G)\). We define
% \(P_{\mathrm{st}}(\ell)\) as the probability that \(d_T(u,v)+1=\ell\), using size two
% for an edge already in \(T\). We have then the following theorem.
% Let $P_{st}(\ell)$ denote the probability that a spanning tree contains exactly $\ell-1$ edges of some simple cycle of length $\ell$. 
Let \(T\) be a uniformly random spanning tree of \(G\), and let \(e=\{u,v\}\) be an edge chosen uniformly at random from \(E(G)\). We define
\[
P_{\mathrm{st}}(\ell)
=\Pr\bigl(d_T(u,v)+1=\ell\bigr),
\]
where \(d_T(u,v)\) denotes the distance between \(u\) and \(v\) in \(T\). For \(e\notin T\), the quantity \(d_T(u,v)+1\) is the length of the unique cycle formed by adding \(e\) to \(T\). If \(e\in T\), then \(d_T(u,v)=1\), and we assign the value \(\ell=2\) by convention. We have then the following theorem.
\begin{theorem}[Dhar--Dhar correspondence ~\cite{DharDhar1997ErasedLoops}]\label{SEC2TH1}
For a finite connected simple undirected unweighted graph $G$, with
the long-time erased-loop convention of~\cite{DharDhar1997ErasedLoops},
\begin{equation}\label{SEC2EQ1}
  P(\ell)
=\frac{P_{st}(\ell)}{\ell}.  
\end{equation}
\end{theorem}
Hence $P(\ell)$ is determined entirely by
spanning tree edge inclusion probabilities. The established correspondence and its use here are discussed in Appendix~\ref{APPA}. The question now is how to determine $P_{st}(\ell)$. To achieve this goal, let us consider separately the cases $\ell > 2$ and $\ell = 2$.

\textbf{Case $\ell > 2$.} Let $C \subset G$ be a fixed simple cycle of length $\ell > 2$ with
\[
E(C) = \{e_1,\dots,e_\ell\}.
\]
Let $\mathcal T(G)$ denote the set of spanning trees of $G$, and define
$P_{st}(\ell)$ as the edge-sampled probability of a fundamental cycle of length $\ell$.
For a fixed cycle $C$, consider the event
\[
|E(T)\cap E(C)|=\ell-1,
\]
that is, $T$ contains all but one edge of $C$.
Equivalently,
\[
E(T)\cap E(C)=E(C)\setminus\{e\}
\quad\text{for some } e\in E(C).
\]
If this occurs, then adding the missing edge $e$ to $T$ creates exactly
one cycle, namely $C$. Indeed, since $T$ is a tree, adding any edge
produces a unique cycle, and here that cycle must be $C$. Hence
\[
U := T\cup\{e\}
\]
is a spanning unicycle whose unique cycle is $C$.
Conversely, let $U$ be a spanning unicycle whose unique cycle is $C$.
Removing any edge $e\in E(C)$ yields a spanning tree
\[
T = U\setminus\{e\}\in\mathcal T(G),
\]
and this tree contains exactly $\ell-1$ edges of $C$.
Moreover, each unicycle $U$ with cycle $C$ corresponds to exactly
$\ell$ pairs $(T,e)$, since there are $\ell$ possible edges in $C$
that can be removed.
Therefore,
\begin{align*}
  P_{st}(\ell)
&=
\frac{1}{|E|\,|\mathcal T(G)|}
\sum_{\substack{C:\, |C|=\ell}}
\#\Bigl\{T\in\mathcal T(G):\\&\quad
|E(T)\cap E(C)|=\ell-1\Bigr\}
.  
\end{align*}

\textbf{Case $\ell = 2$.}
Let $u\neq v$ be two distinct vertices and let
\[
C=\{(u,v),(v,u)\}.
\]
A cycle of length $2$ corresponds to traversing the same edge in
opposite directions (in the directed sense), rather than to a
simple undirected cycle.
In this case,
\[
|E(T)\cap E(C)|=1
\]
meaning that $T$ contains the undirected edge $(u,v)$.
If we ``add'' this edge again, no new physical cycle is created;
the length-two case is the backtracking convention. Hence
\[
P_{st}(2)
=
\frac{\sum_{\substack{C:\, |C|=2}}
\#\Bigl\{
T\in\mathcal T(G):
|E(T)\cap E(C)|=1
\Bigr\}}{|\mathcal T(G)|}
.
\]

For each spanning tree $T\in\mathcal T(G)$, let
\[
E\setminus E(T)
=
\{e\in E(G): e\notin E(T)\}
\]
denote the set of non-tree edges (chords).
If $e=(u,v)\in E\setminus E(T)$, then adding $e$ to $T$
creates a unique cycle, called the fundamental cycle of $e$
with respect to $T$. Its length is
\[
|C_T(e)| = d_T(u,v)+1,
\]
where $d_T(u,v)$ is the graph distance between $u$ and $v$ in $T$.
This identity extends to all edges $e=(u,v)\in E$:
if $e\in E(T)$, we may interpret it as producing a cycle of length $2$.
Hence, for every $\ell\ge 2$,
\[
\mathbf 1_{\{|C_T(e)|=\ell\}}
=
\mathbf 1_{\{d_T(u,v)=\ell-1\}}.
\]
Now consider the following two-stage experiment:
\begin{itemize}
\item Choose $T$ uniformly at random from $\mathcal T(G)$.
\item Choose $e$ uniformly at random from $E(G)$.
\end{itemize}
With the backtracking convention when $e\in T$, the event
$L_T(e)=\ell$ is precisely $\{d_T(u,v)=\ell-1\}$.
Therefore,
\[
P_{st}(\ell)
=
\frac{1}{|\mathcal T(G)|}
\sum_{T\in\mathcal T(G)}
\frac{1}{|E|}
\sum_{e=(u,v)\in E}
\mathbf 1_{\{|C_T(e)|=\ell\}}.
\]

By using $|C_T(e)|=\ell$ if and only if $d_T(u,v)=\ell-1$, we obtain
\begin{align*}
P_{st}(\ell)
&=
\frac{1}{|\mathcal T(G)|}
\sum_{T\in\mathcal T(G)}
\frac{1}{|E|}
\sum_{(u,v)\in E}
\mathbf 1_{\{d_T(u,v)=\ell-1\}} \\
&=
\frac{1}{|E|}
\sum_{(u,v)\in E}
\frac{1}{|\mathcal T(G)|}
\sum_{T\in\mathcal T(G)}
\mathbf 1_{\{d_T(u,v)=\ell-1\}}.
\end{align*}
Particularly, we have
\begin{align}
    P_{st}(2)
&=
\frac{1}{|\mathcal T(G)|}
\sum_{T\in\mathcal T(G)}
\frac{1}{|E|}
\sum_{(u,v)\in E}
\mathbf 1_{\{d_T(u,v)=1\}} \\
&=\frac{1}{|\mathcal T(G)|}\sum_{T\in\mathcal T(G)}\frac{|V|-1}{|E|}\\
&=\frac{|V|-1}{|E|}.
\end{align}
From that we get 
\begin{equation}
    P(2)=\frac{|V|-1}{2|E|}.
\end{equation}

We use now the property of the matrix $Y$ defined in~\eqref{SEC1EQ1}, the relations \eqref{SEC1EQ2} and the definition \eqref{SEC1EQ3}. Recall that for $e=(u,v)\in E$ and a spanning tree $T$,
the fundamental cycle created by adding $e$ satisfies
\[
|C_T(e)|=\ell
\quad\Longleftrightarrow\quad
d_T(u,v)=\ell-1.
\]
In a tree, the graph distance $d_T(u,v)$ equals the length of the
unique simple path between $u$ and $v$, that is,
\[
d_T(u,v)
=
|\mathrm{Path}_T(u,v)|.
\]
Therefore,
\[
\mathbf 1_{\{d_T(u,v)=\ell-1\}}
=
\sum_{\delta\in\Gamma_{\ell-1}(u,v)}
\mathbf 1_{\{\mathrm{Path}_T(u,v)=\delta\}},
\]
where $\Gamma_{\ell-1}(u,v)$ denotes the set of simple
$u$--$v$ paths of length $\ell-1$ in $G$.
Taking expectation with respect to the uniform spanning tree measure,
we obtain
\[
\frac{\sum_{T\in\mathcal T(G)}
\mathbf 1_{\{d_T(u,v)=\ell-1\}}}{|\mathcal T(G)|}
=\sum_{\delta\in\Gamma_{\ell-1}(u,v)}
\mathbb P\!\left(\mathrm{Path}_{T}(u,v)=\delta\right).
\]
By using the determinantal formula
\[
\mathbb P\!\left(\mathrm{Path}_{T}(u,v)=\delta\right)
=
\det\!\big(Y_{E(\delta),E(\delta)}\big)
=
p_\delta,
\]
we finally obtain
\begin{align}
P_{st}(\ell)
&=
\frac{1}{|E|}
\sum_{(u,v)\in E}~~
\sum_{\delta\in\Gamma_{\ell-1}(u,v)}
p_\delta.
\end{align}
We can then derive
\begin{align}
P(\ell)
&=
\frac{1}{\ell|E|}
\sum_{(u,v)\in E}~~
\sum_{\delta\in\Gamma_{\ell-1}(u,v)}
p_\delta.
\end{align}
We define $\pi(\ell)$ as the proportion of loop of size $\ell$ erased.
We hence have the following corollary to Theorem \ref{SEC2TH1}.
\begin{figure}
    \centering
    \includegraphics[width=1\linewidth]{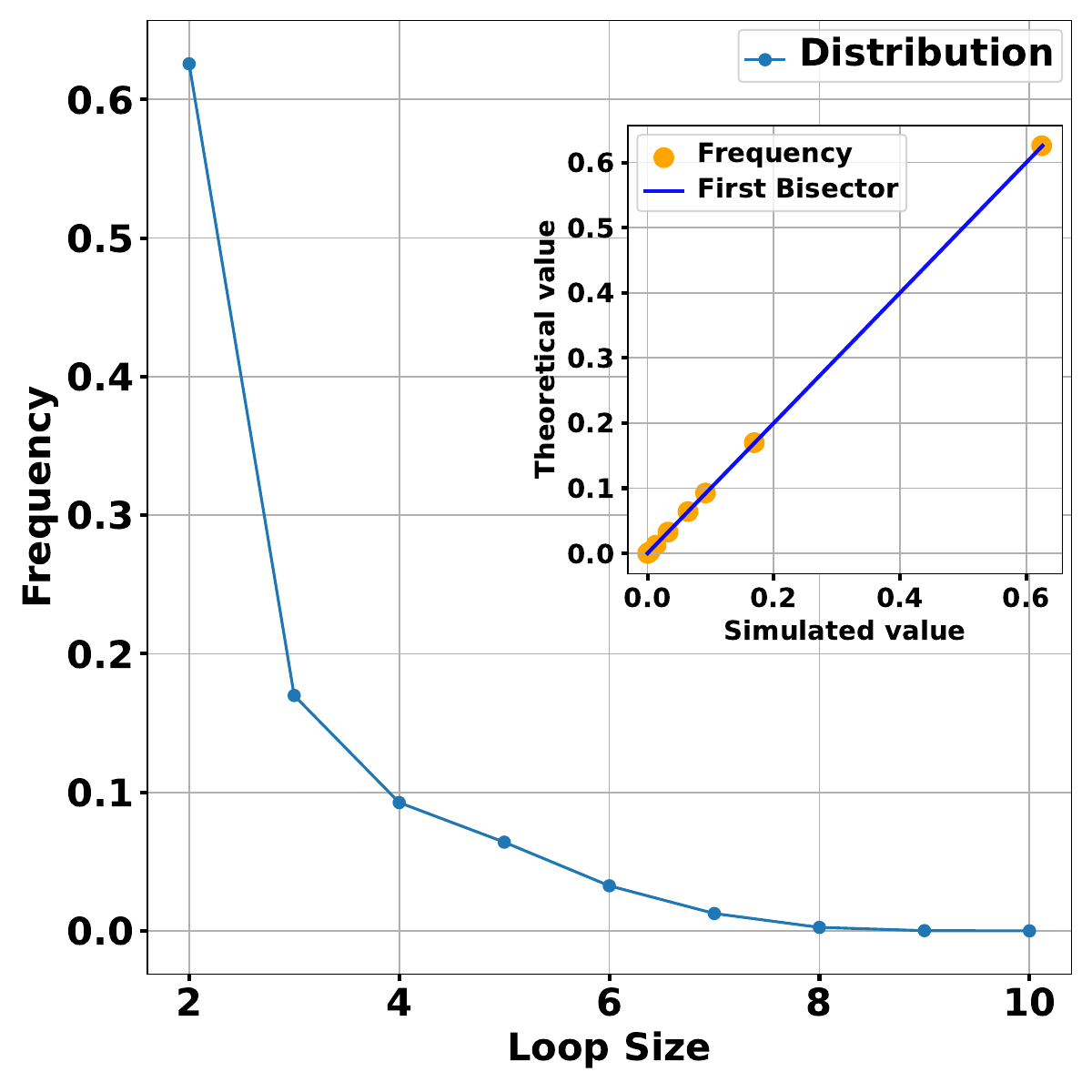}
    \caption{Main panel shows the distribution of erased loops obtained with the relation \ref{SEC2EQ2} on a Watts-Strogatz graph with $10$ nodes, a rewiring probability of $0.5$ and a mean degree of $4$. Inset shows the comparison between Frequencies obtained by the relation \ref{SEC2EQ2} and Frequencies observed with true loop erase random walk on this graph. On the y-axis, we have the frequency obtained with the formula \ref{SEC2EQ2} and on the x-axis, we have the frequency obtained by counting size of erased loops. }
    \label{LOOP}
\end{figure}

\begin{corollary}
For any unweighted, simple, undirected, connected graph $G$,
\begin{equation}\label{SEC2EQ2}
   \pi(\ell)
=\displaystyle\frac{\frac{1}{\ell}
\sum_{(u,v)\in E}~~
\sum_{\delta\in\Gamma_{\ell-1}(u,v)}
p_\delta}{\sum_{k\geq 2}\frac{1}{k}
\sum_{(u,v)\in E}~~
\sum_{\beta \in\Gamma_{k-1}(u,v)}
p_\beta}. 
\end{equation}
\end{corollary}

\begin{proof}
Since $P(\ell)$ is the stationary probability that a newly added step
creates and erases a loop of size $\ell$,
this does not directly give the proportion of a loop of size $\ell$ erased, $\pi(\ell)$, because some new steps do not create a loop. We need to normalize. So we have
\begin{align*}
\pi(\ell)&= \frac{P(\ell)}{\sum_{k\geq 2}P(k)}\\
&=\displaystyle\frac{\frac{1}{\ell |E|}
\sum_{(u,v)\in E}~~
\sum_{\delta\in\Gamma_{\ell-1}(u,v)}
p_\delta}{\sum_{k\geq 2}\frac{1}{k|E|}
\sum_{(u,v)\in E}~~
\sum_{\beta \in\Gamma_{k-1}(u,v)}
p_\beta}\\
&=\displaystyle\frac{\frac{1}{\ell}
\sum_{(u,v)\in E}~~
\sum_{\delta\in\Gamma_{\ell-1}(u,v)}
p_\delta}{\sum_{k\geq 2}\frac{1}{k}
\sum_{(u,v)\in E}~~
\sum_{\beta \in\Gamma_{k-1}(u,v)}
p_\beta}.
\end{align*}
That concludes the proof.
\end{proof}

We used Eq.~\ref{SEC2EQ2} to compute the distribution of erased loops on a Watts--Strogatz graph with $10$ nodes, a rewiring probability of $0.5$, and a mean degree of $4$. 
The resulting distribution is shown in the main panel of Fig.~\ref{LOOP}. 
To validate the theoretical prediction, we compared it with empirical frequencies obtained by counting the sizes of erased loops generated by loop-erased random walks. 
The comparison is presented in the inset of Fig.~\ref{LOOP}. 
The theoretical distribution closely matches the empirical frequencies 
obtained from simulations, showing that formula~\ref{SEC2EQ2} accurately captures 
the distribution of erased loop sizes in this network model ($\varepsilon_{\max}=1.5\times 10^{-3}$ in the inset).

\section{Case of Complete graphs }

\begin{definition}[Complete graph]
The complete graph on $n$ vertices ~\cite{West1996IntroductionGraphTheory}, denoted by $K_n$, is a simple graph with a set
\[
V(K_n)=\{v_1,\dots,v_n\}
\]
and edge set
\[
E(K_n)=\bigl\{\{v_i,v_j\}\;:\;1\le i<j\le n\bigr\}.
\]
Similarly, every pair of distinct vertices is joined by an edge.
\end{definition}

Let us now introduce Theorem~\ref{SEC3TH1} from~\cite{Moon1970CountingLabelledTrees, MeirMoon1970Distance} (for the sake of completeness a proof is given in Appendix~\ref{APPB}) useful in the scope of this work; in particular, we are interested in the following result~\ref{cor:MFTPKn} obtained as a corollary of the former theorem.

\begin{theorem}\label{SEC3TH1} Fix two distinct labels $u,v\in[n]$. For a tree $T_n$, let $S_n=S(T_n;u,v)$ be the
number of vertices on the unique $u$--$v$ path (so the graph distance is
$d(u,v)=S_n-1$).
Assume $T_n$ is chosen uniformly from the $n^{n-2}$ labeled trees on $[n]$.
For $2\le k\le n$,
\begin{equation}\label{SEC3EQ1}
\mathbb{P}\{S_n=k\}
= \frac{k\,(n-2)!}{(n-k)!\,n^{\,k-1}}
= \frac{k}{n-1}\cdot \frac{(n)_k}{n^k},
\end{equation}
where $(n)_k=n(n-1)\cdots(n-k+1)$ is the falling factorial.
Moreover, as $n\to\infty$,
\[
\frac{\mathbb{E}[S_n]}{\sqrt{n}}\longrightarrow \sqrt{\frac{\pi}{2}}.
\]
\end{theorem}
\begin{figure}
    \centering
    \includegraphics[width=1\linewidth]{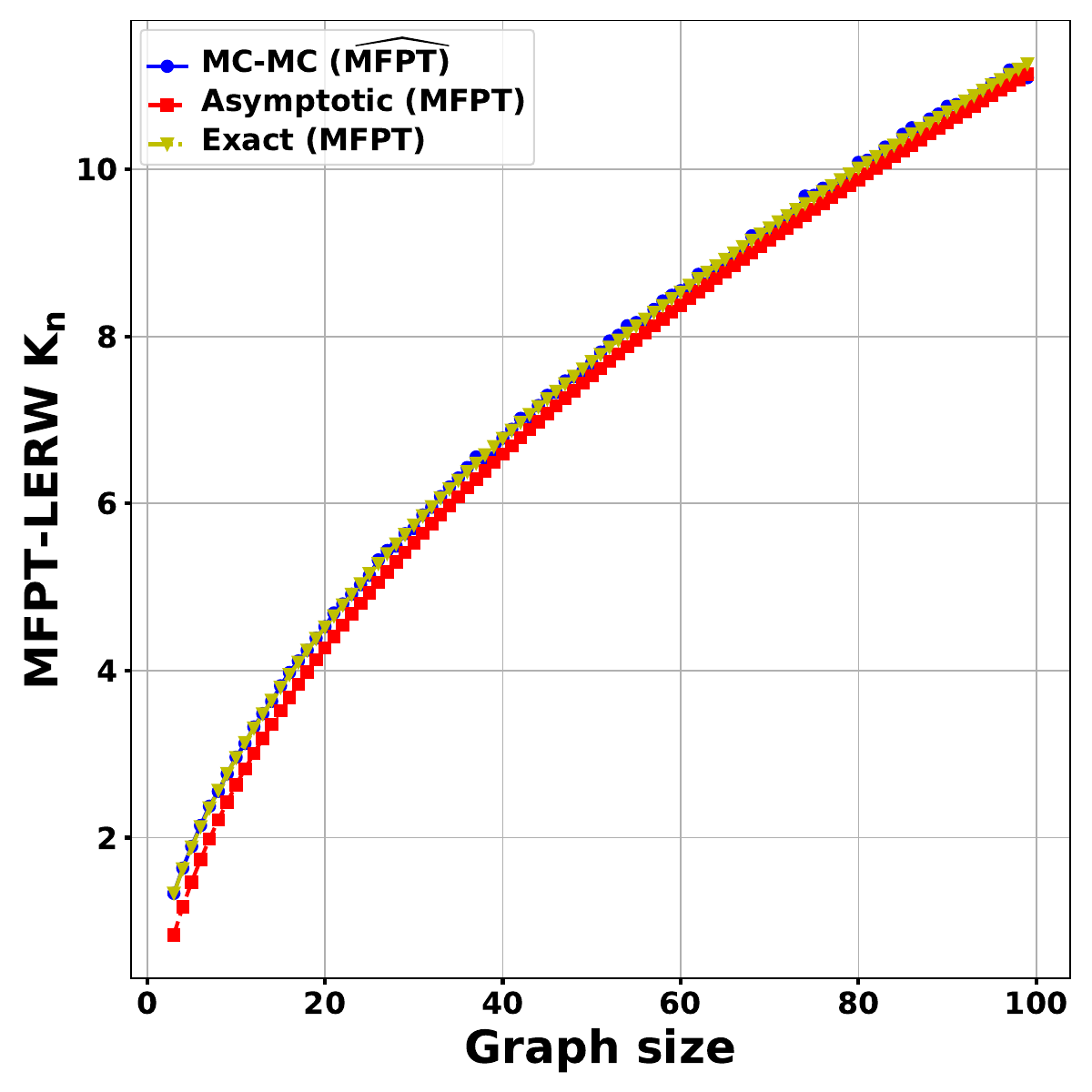}
    \caption{Mean First Passage Time of Loop Erased Random Walks between two pairs of nodes of complete graphs. }
    \label{COMP}
\end{figure}

\begin{corollary}[MFPT for LERW on $K_n$ ]
\label{cor:MFTPKn}
Let $K_n$ be the complete graph with $n$ vertices. For two vertices $i$ and $j$,  the mean first passage time (MFPT) of the loop-erased random walk (LERW) from $i$ to $j$ satisfies
\begin{align}\label{SEC3EQ2}
  \mathrm{MFPT}^{\mathrm{LERW}}(i,j)=\sum_{k=0}^{n-2}\frac{(n-2)!}{n^{k}(n-2-k)!}.  
\end{align}

\end{corollary}
\begin{proof}
Let $i\neq j$ be two vertices of $K_n$ and let $T_n$ be a uniform labelled tree on $[n]$.
The (unique) $i$--$j$ path in $T_n$ has vertex-length $S_n=S(T_n;i,j)$ and edge-length
\[
d(i,j)=S_n-1.
\]
By Theorem~\ref{SEC3TH1}, for $k=2,\dots,n$,
\[
\mathbb{P}\{S_n=k\}
=\frac{k}{n-1}\cdot\frac{(n)_k}{n^k}.
\]
Hence
\begin{align}\label{SEC3EQ3}
\mathbb{E}[d(i,j)]
=\sum_{k=2}^{n}(k-1)\,\mathbb{P}\{S_n=k\}\\
=\sum_{k=2}^{n}(k-1)\frac{k}{n-1}\frac{(n)_k}{n^k}.
\end{align}
Using $(n)_k=n(n-1)(n-2)_{k-2}$, we rewrite \eqref{SEC3EQ3} as
\[
\mathbb{E}[d(i,j)]
=\sum_{k=2}^{n}(k-1)k\,\frac{(n-2)_{k-2}}{n^{k-1}}.
\]
Set $j=k-2$ so that $j=0,\dots,n-2$ and define
\[
a_j:=\frac{(n-2)_j}{n^j}.
\]
Then $a_{j+1}=\frac{n-2-j}{n}a_j$, hence
\begin{align*}
   a_j-a_{j+1}=\frac{j+2}{n}\,a_j\, ,
\end{align*}  
from which we get
\begin{align*}
\frac{(j+1)(j+2)}{n}\,a_j=(j+1)\bigl(a_j-a_{j+1}\bigr). 
\end{align*}
Therefore
\begin{align*}
\mathbb{E}[d(i,j)]
=\sum_{j=0}^{n-2}\frac{(j+1)(j+2)}{n}\,a_j\\
=\sum_{j=0}^{n-2}(j+1)\bigl(a_j-a_{j+1}\bigr).    
\end{align*}
Summation by parts (telescoping) yields
\[
\sum_{j=0}^{n-2}(j+1)\bigl(a_j-a_{j+1}\bigr)
=\sum_{j=0}^{n-2}a_j-(n-1)a_{n-1}.
\]
Since $(n-2)_{n-1}=0$, we have $a_{n-1}=0$, and thus
\[
\mathbb{E}[d(i,j)]=\sum_{j=0}^{n-2}a_j=\sum_{j=0}^{n-2}\frac{(n-2)_j}{n^j}.
\]
Finally, by using $(n-2)_j=\dfrac{(n-2)!}{(n-2-j)!}$, we obtain
\[
 \mathrm{MFPT}^{\mathrm{LERW}}(i,j)=\mathbb{E}[d(i,j)]
=\sum_{j=0}^{n-2}\frac{(n-2)!}{n^{j}(n-2-j)!},
\]
which is exactly the claimed formula.
\end{proof}

We can refine the leading asymptotic obtained from Theorem~\ref{SEC3TH1} by expanding the exact sum~\eqref{SEC3EQ2} for large $n$. Set
\[
b_n(j)=\frac{(n-2)_j}{n^j},\qquad 0\le j\le n-2.
\]
For $j=x\sqrt n$, uniformly for $x$ in any fixed compact subset of $[0,\infty)$, Taylor expansion gives
\begin{align*}
\log b_n(j)
&=\sum_{r=0}^{j-1}\log\!\left(1-\frac{r+2}{n}\right)\\
&=-\frac{j(j+3)}{2n}
-\frac{j(2j^2+9j+13)}{12n^2}
+O\!\left(\frac{j^2}{n^2}+\frac{j^4}{n^3}\right)\\
&=-\frac{x^2}{2}
-\frac{1}{\sqrt n}\left(\frac{3x}{2}+\frac{x^3}{6}\right)
+O(n^{-1}).
\end{align*}
Consequently,
\[
b_n(j)=e^{-x^2/2}\left[1-\frac{1}{\sqrt n}
\left(\frac{3x}{2}+\frac{x^3}{6}\right)+O(n^{-1})\right].
\]
Summing by an Euler--Maclaurin \cite{Guth2008EulerMaclaurin} approximation (with the tails controlled by
Gaussian bounds), we obtain
\begin{align*}
\sum_{j=0}^{n-2}b_n(j)
&=\sqrt n\int_0^\infty e^{-x^2/2}\,dx+\frac12\\
&\quad-\int_0^\infty e^{-x^2/2}
\left(\frac{3x}{2}+\frac{x^3}{6}\right)\,dx
+O(n^{-1/2})\\
&=\sqrt{\frac{\pi n}{2}}+\frac12-\frac32-\frac13
+O(n^{-1/2}).
\end{align*}
Thus the correct constant-order refinement is
\begin{equation}\label{SEC3EQ4}
\mathrm{MFPT}^{\mathrm{LERW}}(i,j)
=\sqrt{\frac{\pi n}{2}}-\frac43+O(n^{-1/2}).
\end{equation}
In particular, the leading asymptotic is
$\mathrm{MFPT}^{\mathrm{LERW}}(i,j)\sim\sqrt{\pi n/2}$.

In Fig.~\ref{COMP}, we used the exact formula \ref{SEC3EQ2}, the asymptotic formula \ref{SEC3EQ4} and a Monte Carlo approach to compute the MFPT between two nodes of different complete graphs.
The results obtained using the different approaches show very good agreement. The corresponding absolute errors relative to the exact MFPT are presented in Appendix~\ref{app:complete-graph-errors}, Fig.~\ref{fig:comp-absolute-errors}. In particular, the error associated with the asymptotic approximation decreases as the graph size increases, tending to zero in the large-graph limit, while the Monte Carlo estimates remain close to the exact values.
In the complete graph $K_n$, the number of edges is
\begin{align}
    |E(K_n)| = \frac{n(n-1)}{2}.
\end{align}
Therefore the probability that a spanning tree contains exactly $\ell-1$ edges of some simple cycle of length $\ell$ in the complete graph $K_n$ is,
\begin{align}
    P_{st}(\ell,K_n)
    =
    \frac{2}{n(n-1)}
    \sum_{(u,v)\in E}~~
    \sum_{\delta\in\Gamma_{\ell-1}(u,v)}
    p_\delta.
\end{align}
From the theorem \ref{SEC3TH1} we obtain
\begin{align}
    P_{st}(\ell,K_n)
    &= \frac{2}{n(n-1)}
    \sum_{(u,v)\in E}
    \mathbb{P}\{S_n=\ell\} \\
    &= \mathbb{P}\{S_n=\ell\} = \frac{\ell}{n-1}\cdot \frac{(n)_\ell}{n^\ell}.
\end{align}
Thus the stationary probability that a newly added step
creates and erases a loop of size $\ell$ is 
\begin{align}
    P(\ell,K_n)
    = \frac{1}{n-1}\cdot \frac{(n)_\ell}{n^\ell}.
\end{align}
\begin{figure}
    \centering
    \includegraphics[width=1\linewidth]{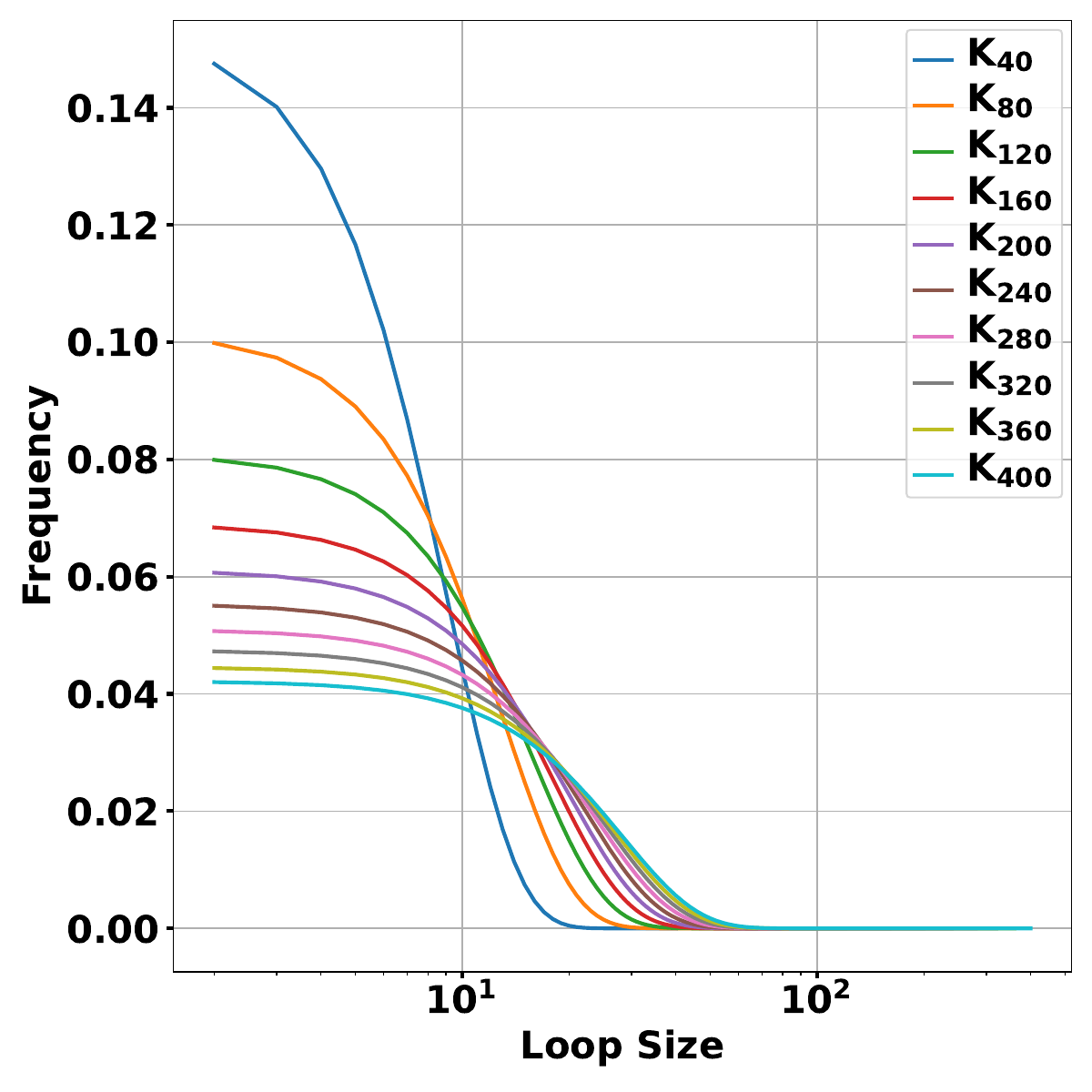}
    \caption{Distribution of erased loops obtained with the relation \ref{SEC3EQ5} on different complete graphs.}
    \label{COMP1}
\end{figure}

This is the well-known solution of the classical birthday problem ~\cite{AlappattuPitman2021ColoredLERW}.
Consequently, the proportion of loop of size $\ell$ erased $\pi_{K_n}(\ell)$is
\begin{align}\label{SEC3EQ5}
    \pi_{K_n}(\ell)
    &= 
    \frac{
        \frac{1}{n-1}\cdot \frac{(n)_\ell}{n^\ell}
    }{
        \sum_{k\geq2}
        \frac{1}{n-1}\cdot \frac{(n)_k}{n^k}
    } = 
    \frac{
        \frac{(n)_\ell}{n^\ell}
    }{
        \sum_{k\geq2}
        \frac{(n)_k}{n^k} 
    }.
\end{align}
We may therefore write
\begin{align}
    \pi_{K_n}(\ell)
    =
    \frac{a_n(\ell)}{\sum_{k=2}^{n} a_n(k)},
\end{align}
where we define
\begin{align}
    a_n(\ell)
    &= \frac{(n)_\ell}{n^\ell} = \prod_{i=0}^{\ell-1}
       \left(1 - \frac{i}{n}\right).
\end{align}
We now study the asymptotic behavior of $a_n(\ell)$, retaining the
finite-size terms that affect the approximation. Taking logarithms and using
$\log(1-z)=-z-z^2/2+O(z^3)$ gives
\begin{align}
\log a_n(\ell)
&=\sum_{j=0}^{\ell-1}\log\left(1-\frac{j}{n}\right)\notag\\
&=-\frac{\ell(\ell-1)}{2n}
-\frac{\ell(\ell-1)(2\ell-1)}{12n^2}
+O\!\left(\frac{\ell^4}{n^3}\right).
\label{SEC3EQLOG}
\end{align}
The expansion is valid with the stated error for $\ell=o(n)$; the
remainder is $o(1)$ for $\ell=o(n^{3/4})$. In the principal scaling regime
$\ell=O(\sqrt n)$, the remainder is $O(n^{-1})$. Dropping the quadratic
correction yields $a_n(\ell)\sim\exp[-\ell^2/(2n)]$ when
$\ell=O(\sqrt n)$.

The exact normalizer in~\eqref{SEC3EQ5} is
\[
Z_n=\sum_{k=2}^n a_n(k).
\]
To obtain its constant-order correction, set $k=x\sqrt n$. Then
\[
a_n(k)=e^{-x^2/2}\left[1+\frac{1}{\sqrt n}
\left(\frac{x}{2}-\frac{x^3}{6}\right)+O(n^{-1})\right]
\]
uniformly for $x$ in compact intervals. The Riemann-sum expansion,
including the endpoint term, yields
\begin{align*}
\sum_{k=0}^{n}a_n(k)
&=\sqrt{n}\int_0^\infty e^{-x^2/2}\,dx+\frac12\\
&\quad+\int_0^\infty e^{-x^2/2}
\left(\frac{x}{2}-\frac{x^3}{6}\right)\,dx+O(n^{-1/2})\\
&=\sqrt{\frac{\pi n}{2}}+\frac23+O(n^{-1/2}).
\end{align*}
Because $a_n(0)=a_n(1)=1$, the required asymptotic behavior normalizer is
\begin{equation}\label{SEC3EQZN}
Z_n=\sqrt{\frac{\pi n}{2}}-\frac43+O(n^{-1/2}).
\end{equation}
This normalizer is identical, at finite $n$, to the MFPT sum
in~\eqref{SEC3EQ2}, even though their first two asymptotic terms agree.
For $\ell=\lfloor x\sqrt n\rfloor$ with fixed $x>0$ we obtain the
 the refined finite $n$ approximation
\begin{equation}\label{SEC3EQ6}
\pi_{K_n}(\ell)\approx
\frac{\exp\!\left[-\frac{\ell(\ell-1)}{2n}
-\frac{\ell(\ell-1)(2\ell-1)}{12n^2}\right]}
{\sqrt{\pi n/2}-4/3}.
\end{equation}
For $\ell=O(\sqrt n)$, the relative error in this refined
approximation is $O(n^{-1})$; neither asymptotic formula is intended to
be uniformly accurate for all $2\le \ell\le n$.
\begin{figure}
    \centering
    \includegraphics[width=1\linewidth]{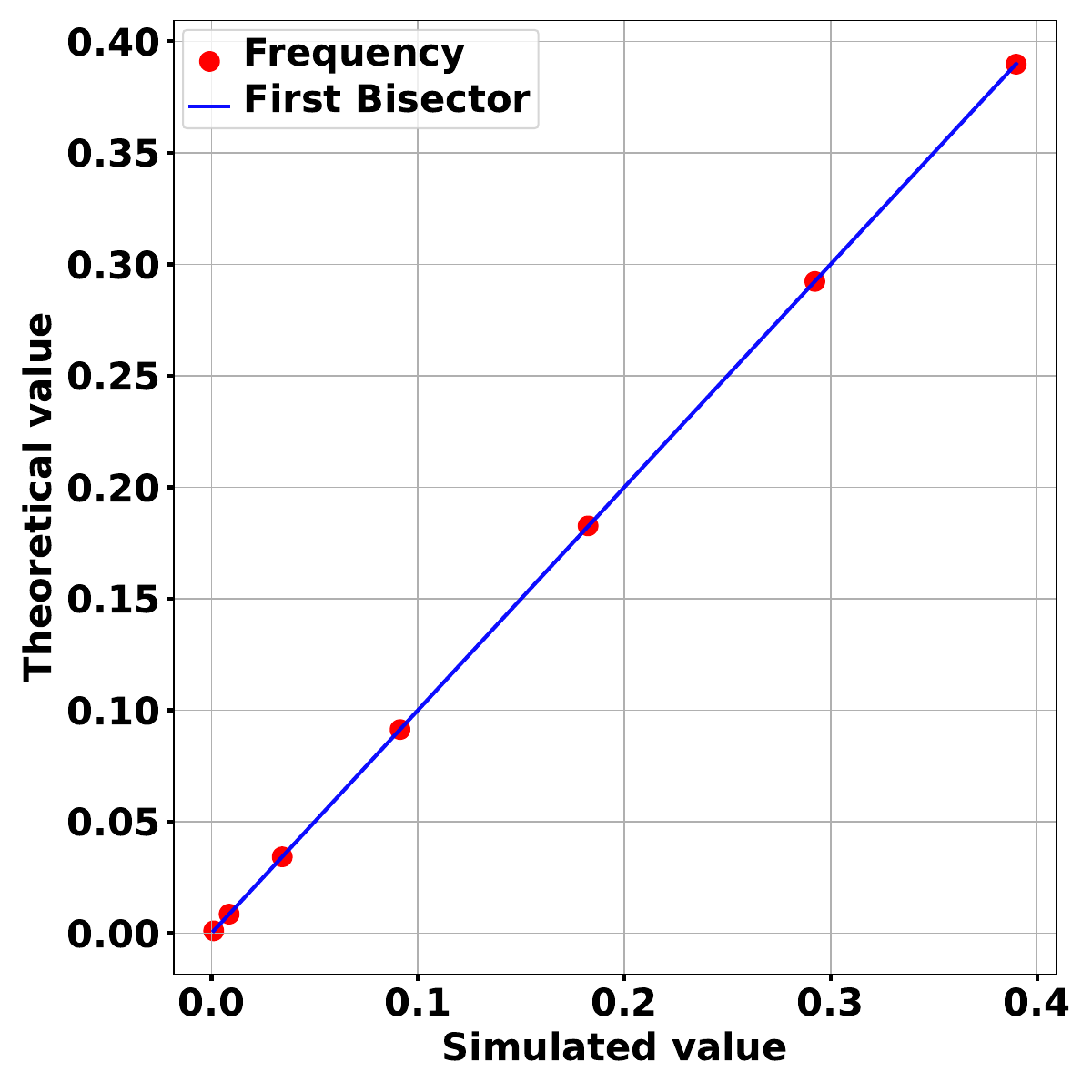}
    \caption{Comparison between Frequencies obtained by the relation~\eqref{SEC3EQ5} and Frequencies observed with true loop erase random walk on $K_{8}$. On the y-axis, we have the frequency obtained with the exact formula~\eqref{SEC3EQ5} and on the x-axis, we have the frequency obtained by counting size of erased loops.}
    \label{COMP4}
\end{figure}
\begin{figure}
    \centering
    \includegraphics[width=1\linewidth]{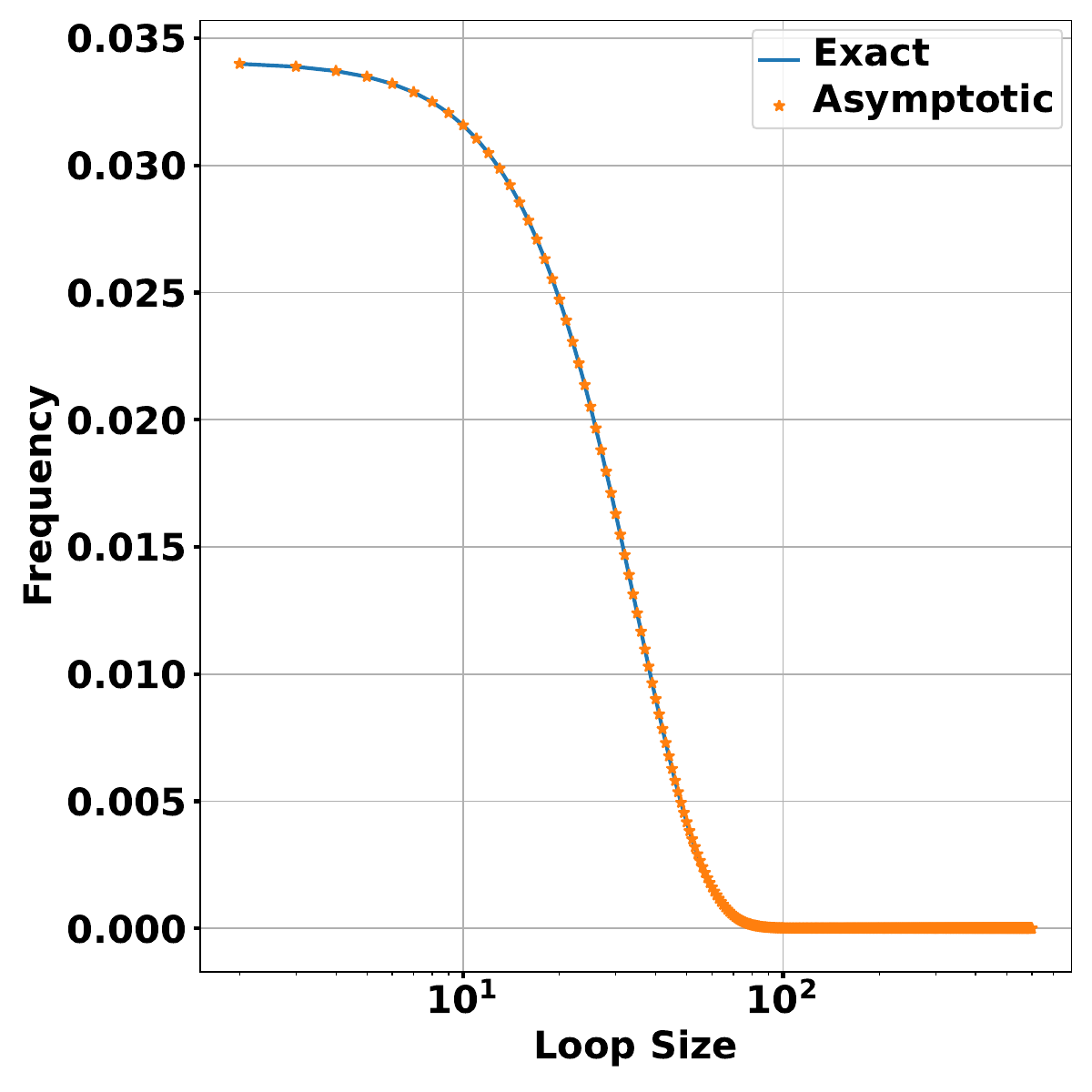}
    \caption{Comparison between distribution of erased loops obtained with the exact relation \eqref{SEC3EQ5} and the asymptotic relation \eqref{SEC3EQ6} on the complete graph $K_{500}$.}
    \label{COMP2}
\end{figure}
We used formula~\eqref{SEC3EQ5} to compute the distribution of erased loops on different complete graphs. 
The results are presented in Fig.~\ref{COMP1}. 
In Fig.~\ref{COMP4}, we compare the theoretical distribution obtained from this formula and the empirical distribution obtained by counting the sizes of erased loops produced by loop-erased random walks. 
This comparison gives $\varepsilon_{\max}=2.2\times 10^{-4}$, providing a quantitative check of the proposed formula. In addition, Fig.~\ref{COMP2} shows that the asymptotic relation~\eqref{SEC3EQ6} provides a good approximation for large complete graphs. 
As the size of the graph increases, the asymptotic expression becomes increasingly accurate and closely matches the exact distribution.

\section{Case of cycle Graphs}
\begin{definition}[Cycle graph]
For an integer $n\ge 3$, the cycle graph $C_n$ ~\cite{West1996IntroductionGraphTheory} is the simple graph with vertex set
\[
V(C_n)=\{v_1,v_2,\dots,v_n\},
\]
and edge set
\[
E(C_n)=\bigl\{\{v_1,v_2\},\{v_2,v_3\},\dots,\{v_{n-1},v_n\},\{v_n,v_1\}\bigr\}.
\]
Equivalently, $C_n$ is a single simple cycle of length $n$. In particular, every vertex has degree $2$, and
\[
|V(C_n)|=n,\qquad |E(C_n)|=n.
\]
\end{definition}
\begin{theorem}[Number of spanning trees of $C_n$]
Let $C_n$ be the cycle graph with $n\ge 3$, and denote by $\mathcal{T}(C_n)$ the set of spanning trees of $C_n$. Then
\[
|\mathcal{T}(C_n)|=n.
\]
\end{theorem}
\begin{proof}
A spanning tree of a connected graph on $n$ vertices must contain exactly $n-1$ edges. In $C_n$ there is exactly one cycle (the whole graph). Therefore, to obtain a spanning tree from $C_n$, we must break this unique cycle by deleting one edge. After removing one edge, the graph remains connected and becomes a simple path on $n$ vertices, hence a spanning tree. 
Conversely, removing any single edge of $C_n$ indeed produces a connected acyclic graph (a path), so every choice of an edge yields a spanning tree. Since $C_n$ has $n$ edges, there are $n$ possible choices, hence
\[
|\mathcal{T}(C_n)| = |E(C_n)| = n.
\]
\end{proof}
\begin{theorem}[MFPT for LERW on $C_n$ as a function of the distance]
Let $C_n$ be the cycle graph with $n\ge 3$. For two vertices $i$ and $j$, let
\[
d(i,j)=k
\]
be the length of the shortest path between $i$ and $j$ in $C_n$. Then the mean first passage time (MFPT) of the loop-erased random walk (LERW) from $i$ to $j$ satisfies
\begin{equation}\label{SEC4EQ1}
     \mathrm{MFPT}^{\mathrm{LERW}}(i,j)=\frac{2k(n-k)}{n}.
\end{equation}
\end{theorem}
\begin{proof}
On a cycle, between two distinct vertices $i$ and $j$ there are exactly two simple (self-avoiding) paths: one of length $k$ and the other of length $n-k$. The loop-erasure mechanism forces the trajectory to be compatible with one of these two simple paths once the walk is ``projected'' onto a spanning tree description.
Now recall that the spanning trees of $C_n$ are obtained by deleting exactly one edge. Fix $i$ and $j$. Consider the two arcs connecting $i$ to $j$: the short arc of length $k$, the long arc of length $n-k$. If the deleted edge lies on the long arc, then the unique simple path between $i$ and $j$ in the resulting spanning tree is the short arc (length $k$). There are exactly $(n-k)$ edges on the long arc, hence $(n-k)$ spanning trees for which the $i$--$j$ path has length $k$. Similarly, if the deleted edge lies on the short arc, then the unique path between $i$ and $j$ in the spanning tree is the long arc (length $n-k$). There are $k$ edges on the short arc, hence $k$ spanning trees for which the $i$--$j$ path has length $n-k$. Therefore, when spanning trees are sampled uniformly from $\mathcal{T}(C_n)$, the $i$--$j$ distance in the tree equals $k$ with probability $\frac{n-k}{n}$ and equals $n-k$ with probability $\frac{k}{n}$. Plugging these two contributions into the LERW MFPT expression for the cycle yields
\[
 \mathrm{MFPT}^{\mathrm{LERW}}(i,j)=\frac{2k(n-k)}{n}.
\]
\end{proof}
\begin{theorem}
Let the cycle graph $C_n$ (with $n \ge 3$). Then the average mean first passage time of LERW between two distinct vertices $i$ and $j$ of $C_n$ is
\[
\overline{ \mathrm{MFPT}}^{\mathrm{LERW}} = \frac{n+1}{3}.
\]
\end{theorem}

\begin{proof}
From the previous result, if $d(i,j)=k$ then
\[
 \mathrm{MFPT}^{\mathrm{LERW}}(i,j)=\frac{2k(n-k)}{n}.
\]
So to compute the average over all unordered pairs $\{i,j\}$, we must weight by the number of pairs at each distance $k$.

\textbf{Case 1: $n$ odd.}
Write $n=2m+1$, so the possible distances are $k=1,2,\dots,m$.
For each fixed $k\in\{1,\dots,m\}$, every vertex $i$ has exactly two vertices at distance $k$ (one in each direction along the cycle), hence the number of ordered pairs $(i,j)$ with $d(i,j)=k$ is $2n$. Therefore the number of unordered pairs $\{i,j\}$ at distance $k$ is
\[
N_k=\frac{2n}{2}=n.
\]
Thus,
\begin{align*}
\sum_{1\le i<j\le n} \mathrm{MFPT}^{\mathrm{LERW}}(i,j)
&=\sum_{k=1}^{m} N_k \,\frac{2k(n-k)}{n}\\
&=\sum_{k=1}^{m} 2k(n-k).
\end{align*}
Let us now compute
\[
\sum_{k=1}^{m} 2k(n-k)
=2\left(n\sum_{k=1}^{m}k-\sum_{k=1}^{m}k^2\right).
\]
Using
\[
\sum_{k=1}^{m}k=\frac{m(m+1)}{2}\text{ and }
\sum_{k=1}^{m}k^2=\frac{m(m+1)(2m+1)}{6},
\]
we get
\begin{align*}
   \sum_{1\le i<j\le n} \mathrm{MFPT}^{\mathrm{LERW}}(i,j)
=&2m(m+1)\left(\frac{n}{2}-\frac{2m+1}{6}\right)\\
=&\frac{(n-1)n(n+1)}{6}.
\end{align*}
Finally, the number of unordered pairs is \(
\binom{n}{2}=\frac{n(n-1)}{2},\)
so the average is
\begin{align*}
  \overline{\mathrm{MFPT}}^{\mathrm{LERW}}
&=\frac{(n-1)n(n+1)}{6}\frac{2}{n(n-1)}=\frac{n+1}{3}.  
\end{align*}
 \textbf{Case 2: $n$ even.}
Write $n=2m$. The possible distances are $k=1,2,\dots,m$.
For $k=1,\dots,m-1$, each vertex has exactly two vertices at distance $k$, so as before the number of unordered pairs at distance $k$ is $N_k=n$.
For $k=m=n/2$, each vertex has exactly one opposite vertex at distance $m$, hence the number of ordered pairs is $n$ and the number of unordered pairs is
\[
N_m=\frac{n}{2}=m.
\]
Therefore
\begin{align*}
\sum_{1\le i<j\le n} \mathrm{MFPT}^{\mathrm{LERW}}(i,j)
=&\sum_{k=1}^{m-1} \frac{2kn(n-k)}{n}
\;+\;  \frac{2m^2(n-m)}{n}\\
=&\sum_{k=1}^{m-1} 2k(2m-k) \;+\; m^2.
\end{align*}
Compute the first sum
\[
\sum_{k=1}^{m-1} 2k(2m-k)
=2\left(2m\sum_{k=1}^{m-1}k-\sum_{k=1}^{m-1}k^2\right).
\]
By recalling that
\[
\sum_{k=1}^{m-1}k=\frac{(m-1)m}{2} \text{ and }
\sum_{k=1}^{m-1}k^2=\frac{(m-1)m(2m-1)}{6},
\]
we obtain
\[
\sum_{k=1}^{m-1} 2k(2m-k)=\frac{(m-1)m(4m+1)}{3}.
\]
Hence
\begin{align*}
  \sum_{1\le i<j\le n} \mathrm{MFPT}^{\mathrm{LERW}}(i,j)
&=\frac{(m-1)m(4m+1)}{3}+m^2\\
&=\frac{m(4m^2-1)}{3}.  
\end{align*}
Now $\binom{n}{2}=\binom{2m}{2}=m(2m-1)$, so
\begin{align*}
   \overline{\mathrm{MFPT}}^{\mathrm{LERW}}
&=\frac{\frac{m(4m^2-1)}{3}}{m(2m-1)}=\frac{2m+1}{3}=\frac{n+1}{3}
\end{align*}
In both cases, we conclude that
\[
\overline{\mathrm{MFPT}}^{\mathrm{LERW}}=\frac{n+1}{3}.
\]
\end{proof}

\begin{figure}
    \centering
    \includegraphics[width=1\linewidth]{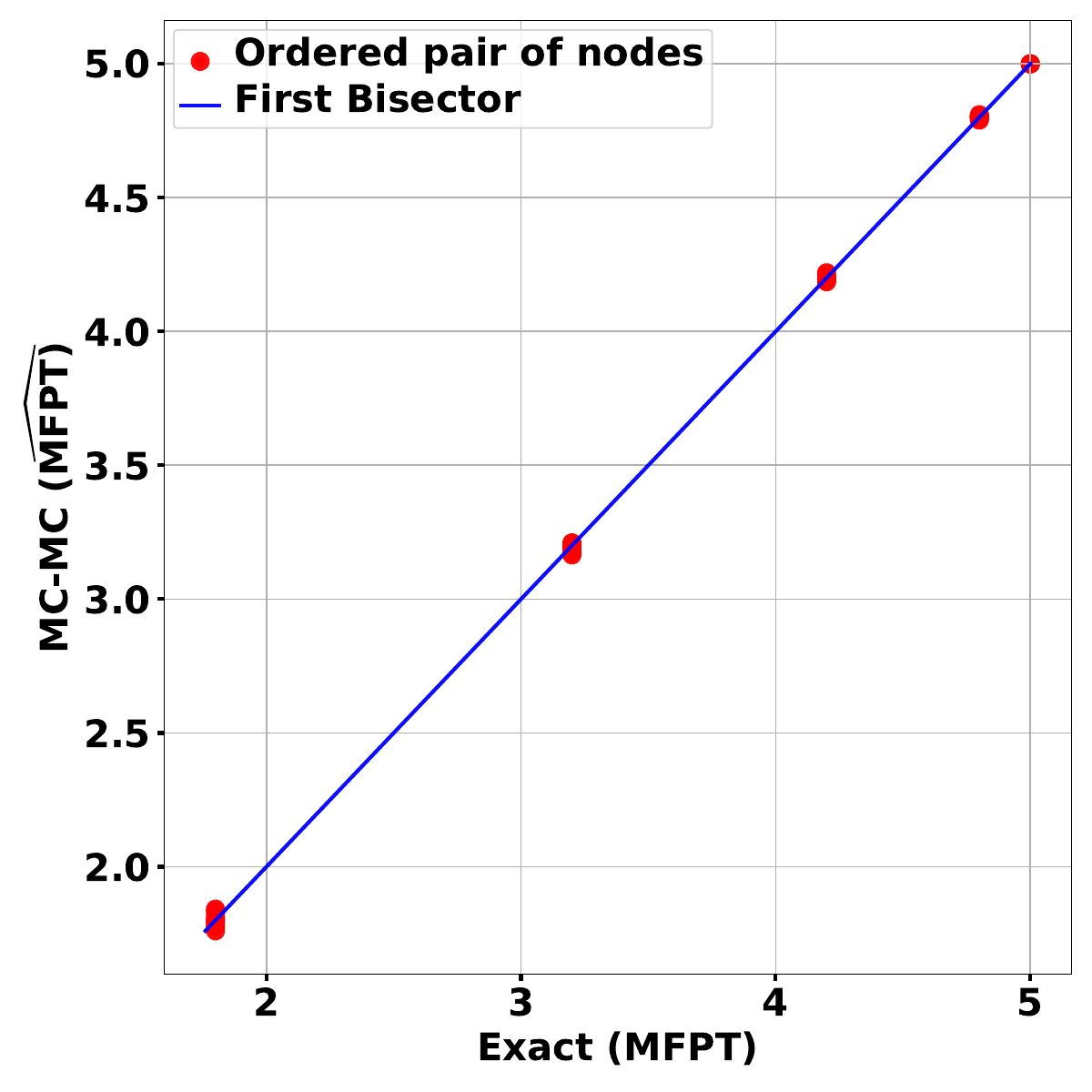}
    \caption{Mean First Passage Time of Loop Erased Random Walks between different pairs of nodes of cycle graph $C_{10}$. On the x-axis, we have the MFPT obtained with the algorithm~\ref{alg:1} and on the y-axis, we have the MFPT obtained with an MC-MC approach with empirical average.}
    \label{CYC}
\end{figure}

In Fig.~\ref{CYC}, we used the exact formula \eqref{SEC4EQ1} and a Monte Carlo approach to compute the MFPT between two nodes of cycle graph $C_{10}$. The empirical estimate $\widehat{\mathrm{MFPT}}^{\mathrm{LERW}}_{\mathrm{MC}}$ closely follows the exact MFPT, with $\varepsilon_{\max}=3.6\times 10^{-2}$.

In the cycle graph $C_n$, every spanning tree is obtained by deleting exactly one edge from the unique cycle. 
As a consequence, each spanning tree is a simple path on $n$ vertices. 
This linear structure severely restricts the possible fundamental cycles that can be formed. 
If we consider an edge that already belongs to the tree, it corresponds to a loop of size $2$, 
which represents the trivial backtracking along a tree edge. 
On the other hand, there is exactly one non-tree edge, namely the edge removed from the original cycle. 
Adding this missing edge reconnects the two endpoints of the path. 
This operation reconstructs the unique global cycle of the graph. 
The length of this cycle is precisely $n$. 
Therefore, in $C_n$, erased loops can only have size $2$ or size $n$. We obtain then 
\begin{align*}
     P_{st}(2,C_n) 
&=\frac{n-1}{n} \text{ and } P(2,C_n)
=\frac{n-1}{2n}.
\end{align*}
We also get
\begin{align*}
    P_{st}(n,C_n)&=\frac{1}{|E|}
\sum_{(u,v)\in E}
\frac{1}{|\mathcal T(G)|}
\sum_{T\in\mathcal T(G)}
\mathbf 1_{\{d_T(u,v)=n-1\}}\\
&=\frac{1}{n}
\sum_{(u,v)\in E}
\frac{1}{n}
\sum_{T\in\mathcal T(G)}
\mathbf 1_{\{d_T(u,v)=n-1\}}\\
&=\frac{1}{n}
\sum_{(u,v)\in E}
\frac{n}{n}=\frac{1}{n}.
\end{align*}
From that,
\begin{align*}
     P(n,C_n)=\frac{1}{n^2}.
\end{align*}
This gives the proportion of loop of size $2$ erased, $\pi_{C_n}(2)$
\begin{align}\label{SEC4EQ2}
   \pi_{C_n}(2)= \displaystyle \frac{\frac{n-1}{2n}}{\frac{n-1}{2n}+\frac{1}{n^2}}= \displaystyle \frac{\frac{n-1}{2n}}{\frac{n^2-n+2}{2n^2}}= \displaystyle \frac{n^2-n}{n^2-n+2}
\end{align}
and the proportion of loop of size $n$ erased, $\pi_{C_n}(n)$
\begin{align}\label{SEC4EQ3}
   \pi_{C_n}(n)= \displaystyle \frac{\frac{1}{n^2}}{\frac{n-1}{2n}+\frac{1}{n^2}}= \displaystyle \frac{\frac{1}{n^2}}{\frac{n^2-n+2}{2n^2}}= \displaystyle \frac{2}{n^2-n+2}.
\end{align} 
We can observe that 
\begin{align*}
   \lim_{n\rightarrow{\infty}} \pi_{C_n}(n)=  \lim_{n\rightarrow{\infty}} \frac{2}{n^2-n+2}
   =  0
\end{align*}
and 
\begin{align*}
   \lim_{n\rightarrow{\infty}} \pi_{C_n}(2)&=  \lim_{n\rightarrow{\infty}} \frac{n^2-n}{n^2-n+2}
   =  1.
\end{align*}

\begin{figure}
    \centering
    \includegraphics[width=1\linewidth]{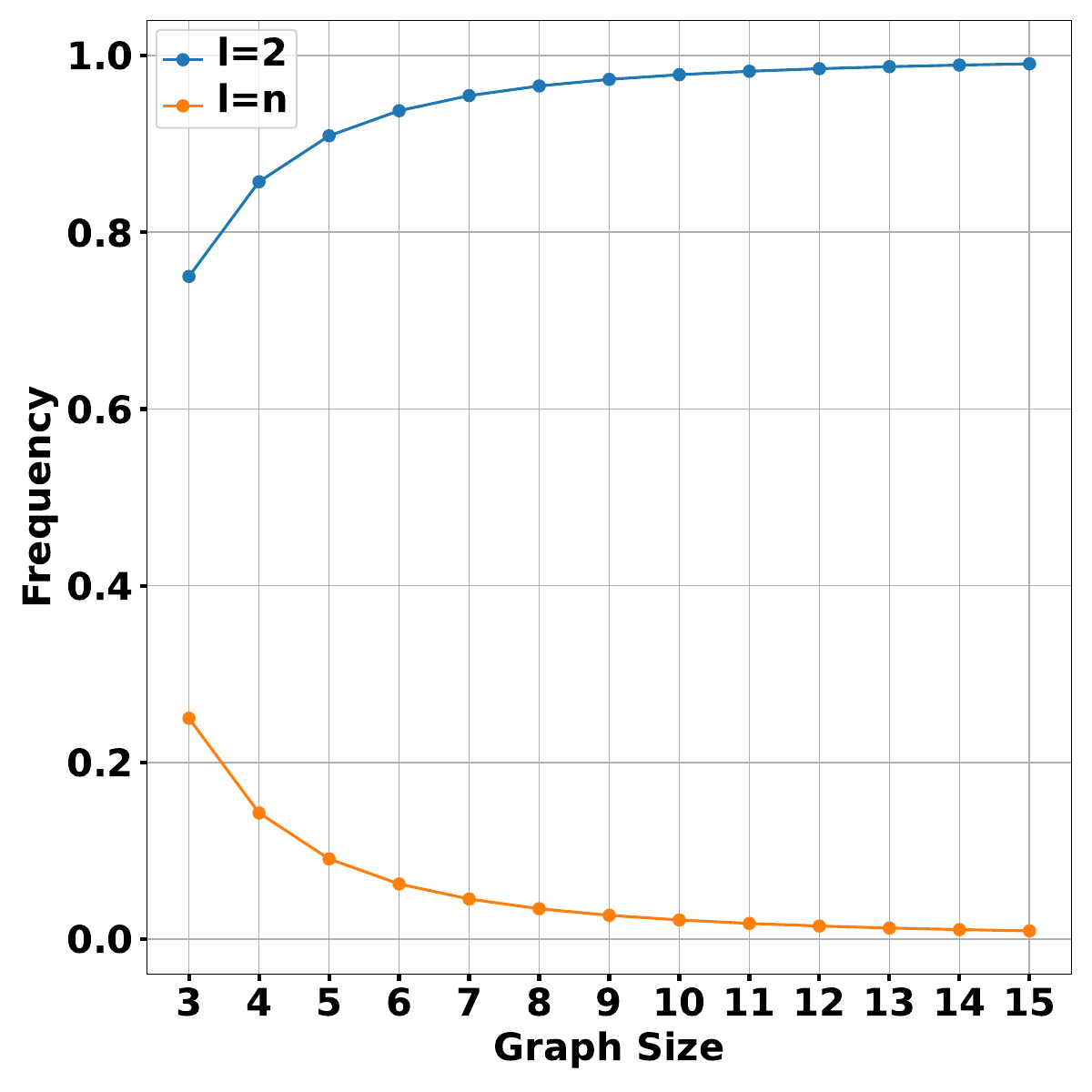}
    \caption{Distribution of erased loops obtained with the relation \eqref{SEC4EQ3} on different cycle graphs.}
    \label{CYC1}
\end{figure}
\begin{figure}
    \centering
    \includegraphics[width=1\linewidth]{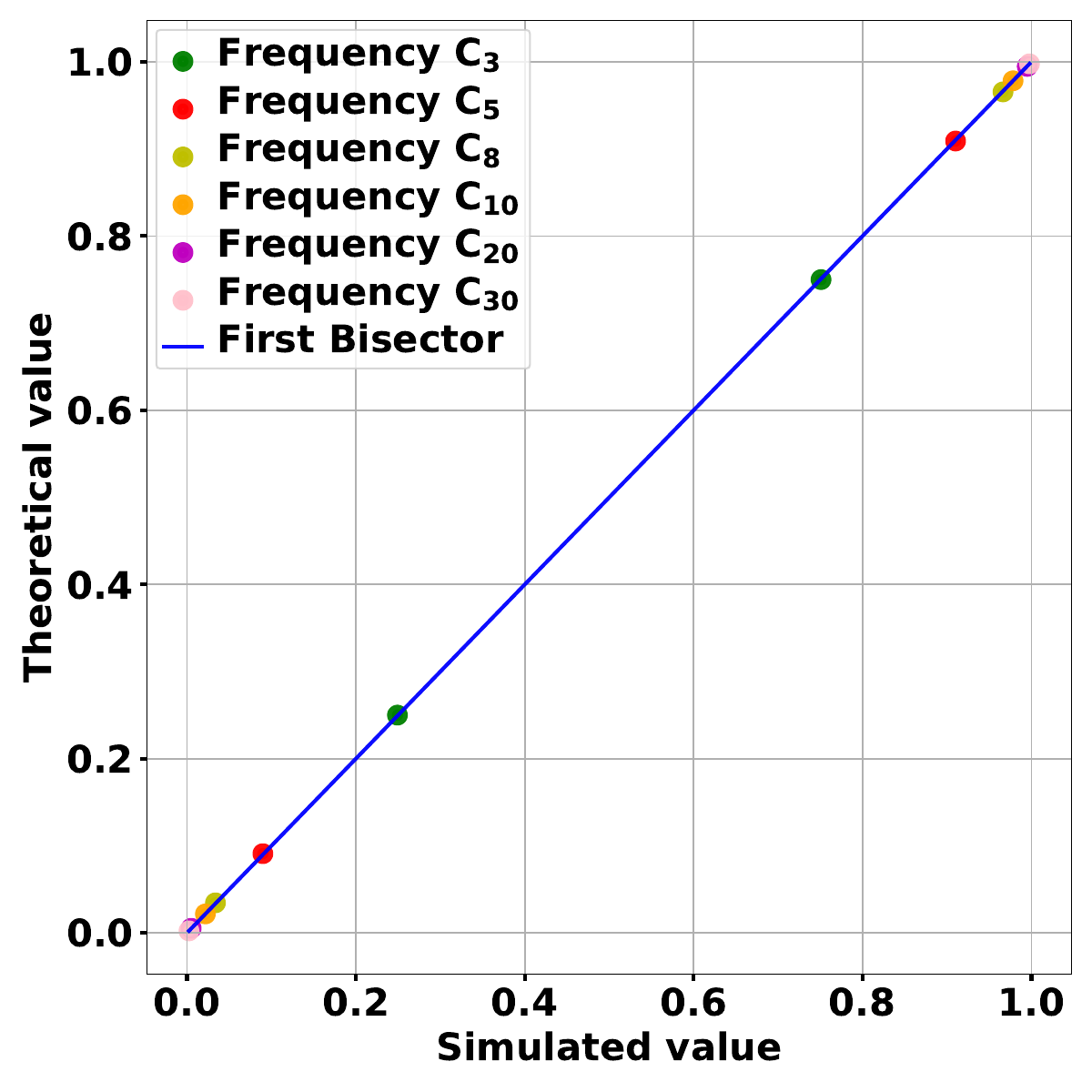}
    \caption{Comparison between Frequencies obtained by the relation \eqref{SEC4EQ3} and Frequencies observed with true loop erase random walk on different cycle graphs. On the y-axis, we have the frequency obtained with the exact formula \eqref{SEC4EQ3} and on the x-axis, we have the frequency obtained by counting size of erased loops.}
    \label{CYC2}
\end{figure}
We used formula~\eqref{SEC4EQ3} to compute the distribution of erased loops on different cycle graphs. 
The results are presented in Fig.~\ref{CYC1}. From this figure, we can clearly observe a limiting effect. 
As the size of the cycle graph increases, the probability that an erased loop has size two becomes increasingly dominant. In other words, for large cycle graphs, most erased loops correspond to simple backtracking events.
This behavior is consistent with the theoretical predictions obtained from the analytical formula. 
To validate the accuracy of the formula, we compared these theoretical results with empirical measurements obtained by directly counting the sizes of erased loops generated by loop-erased random walks. 
The comparison is shown in Fig.~\ref{CYC2}. The maximum discrepancy in Fig.~\ref{CYC2} is $\varepsilon_{\max}=1.0\times 10^{-6}$. 
Therefore, the proposed formula accurately captures the distribution of erased loop sizes on cycle graphs. 
This comparison provides a numerical check of the proposed spanning-tree model; the exact agreement with the dynamic erasure process requires the correspondence discussed in Appendix~\ref{APPA}.

\section{Motif-based network case}
We build a graph $G=(V,E)$ by gluing together small subgraphs called motifs ~\cite{Milo2002R,KarrerNewman2010RandomGraphsSubgraphs}.  
Each motif $\sigma\in\mathcal M$ lives on a subset of vertices $V(\sigma)\subseteq V$ and contributes its own internal
edge set $E(\sigma)\subseteq \binom{V(\sigma)}{2}$. The full graph is obtained by taking the union of all motif edges
\[
E=\bigcup_{\sigma\in\mathcal M}E(\sigma).
\]
Motifs may share vertices, but they cannot share edges. Furthermore, the edges belonging to them are simply the original edges of $G$.
A convenient way to encode how motifs are attached is the bipartite factor graph
\[
\mathcal B=\bigl(V\cup\mathcal M,\;E_{\mathcal B}\bigr),
\qquad
(u,\sigma)\in E_{\mathcal B}\iff u\in V(\sigma).
\]
So $\mathcal B$ has two types of nodes: original vertices and motif-nodes; an edge in $\mathcal B$ simply means
“vertex $u$ belongs to motif $\sigma$”.
In what follows we assume that this factor graph $\mathcal B$ is a tree. This is a crucial structural
assumption: it means motifs are connected together without creating additional loops between motifs.
As a consequence, any cycle in $G$ must be contained entirely inside a single motif; cycles cannot be formed by
walking through several motifs.

\begin{figure}[tbp]
    \centering
    \begin{subfigure}{0.40\textwidth}
    \centering
    \includegraphics[width=1.0\linewidth]{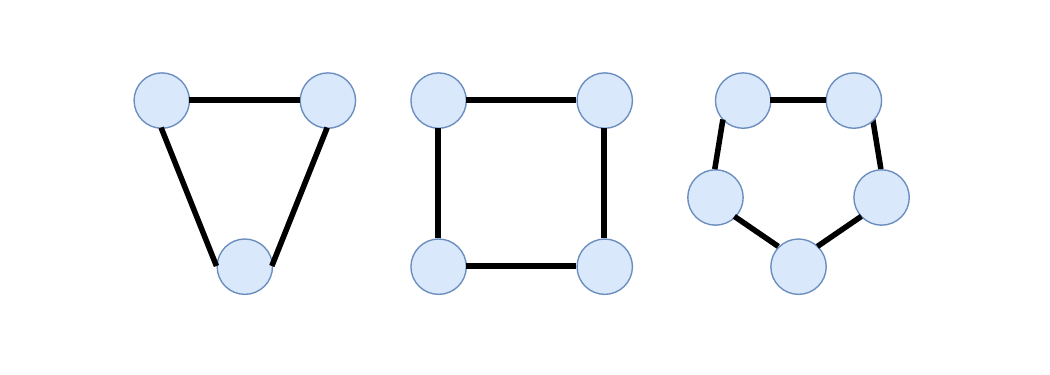}
    \caption{Cycle Motifs}
\end{subfigure}
\hfill
\begin{subfigure}{0.40\textwidth}
    \centering
    \includegraphics[width=1.0\linewidth]{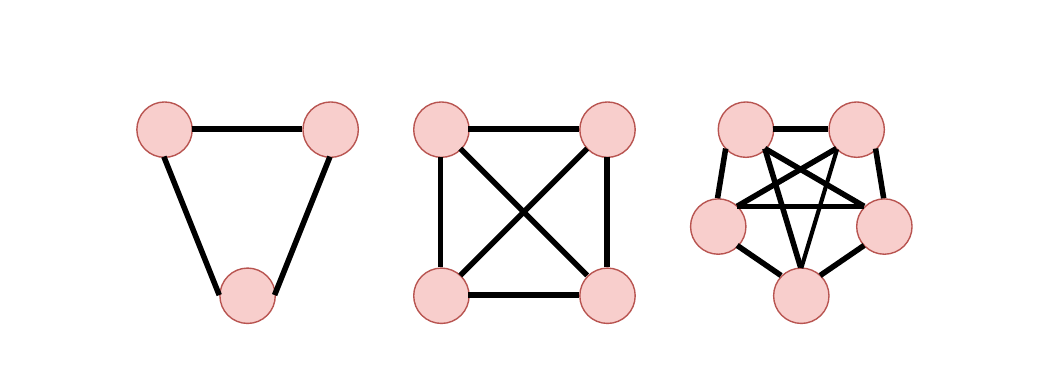}
    \caption{Clique Motifs}
\end{subfigure}
    \caption{Motifs: (a) The cycle motifs of three nodes, four nodes and five nodes. (b) The clique motifs of three nodes, four nodes and five nodes.}
    \label{MOT}
\end{figure}

Let us first consider the case where all the motifs are cycles.
\begin{definition}
A motif $\sigma$ is a cycle motif of length $\ell\ge 3$ ~\cite{Newman2019SpectraShortLoops} if its vertices can be ordered
$V(\sigma)=\{v_1,\dots,v_\ell\}$ such that
\[
E(\sigma)=\bigl\{\{v_i,v_{i+1}\}: i=1,\dots,\ell\bigr\},
\qquad v_{\ell+1}\equiv v_1.
\]    
\end{definition}
For convenience, a single edge not belonging to any longer cycle is treated as a “cycle motif of length $2$”. In panel (a) of Fig.~\ref{MOT} we show few examples of cycle motifs. 

Let
\[
S_u:=\{\sigma\in\mathcal M:\;u\in V(\sigma)\}
\]
be the motifs incident to a vertex $u$. The degree splits motif-by-motif
\[
k_u=\sum_{\sigma\in S_u}k_u^{(\sigma)},
\]
where $k_u^{(\sigma)}$ is the degree of $u$ inside $\sigma$.
If $\sigma$ is a cycle $C_\ell$ and $u\in V(\sigma)$ then $k_u^{(\sigma)}=2$, hence if $u$ belongs to $c_u$ cycles
(and to no other motifs) then $k_u=2c_u$. If motifs are edge-disjoint (no edge is repeated across motifs), then each cycle $C_{\ell_\sigma}$ contributes exactly
$\ell_\sigma$ edges, and therefore
\[
|E(G)|=\sum_{\sigma\in\mathcal M}|E(\sigma)|,
\]
where $|E(\sigma)|=\ell_\sigma$ for a genuine cycle motif,
whereas an isolated-edge motif has $|E(\sigma)|=1$.
The convention $\ell_\sigma=2$ for immediate backtracking
refers to an erased-loop size, not to two physical graph edges.

\begin{theorem}[Spanning trees of a cycle-motif tree]\label{SEC5TH1}
Let $G$ be connected and built from cycle motifs (and isolated edges as $2$-cycles), with edge-disjoint motifs and a
factor graph that is a tree. Let $\mathcal C\subseteq\mathcal M$ be the set of cycle motifs with length $\ge 3$.
Then the number of spanning trees satisfies
\[
|\mathcal{T}(G)|=\prod_{\sigma\in\mathcal C}\ell_\sigma.
\]
\end{theorem}

\begin{proof}
Only genuine cycles create choices. For each $\sigma\in\mathcal C$, choose exactly one edge of $\sigma$ to be deleted.
There are $\ell_\sigma$ options for motif $\sigma$, so the number of global choices is
$\prod_{\sigma\in\mathcal C}\ell_\sigma$. Deleting one edge from each $\sigma\in\mathcal C$ breaks the unique cycle in that motif while keeping it connected,
and isolated edges remain unchanged. Because the factor graph is a tree, no cycle can use edges from multiple motifs.
Therefore, after breaking every genuine cycle, the resulting graph is connected and acyclic, hence a spanning tree. Conversely, any spanning tree cannot contain all edges of any cycle motif $\sigma\in\mathcal C$ (otherwise it would
contain a cycle). It also cannot remove two or more edges from $\sigma$ (that would disconnect the motif). So it removes exactly on edge from each $\sigma\in\mathcal C$, yielding a unique choice. This gives a bijection between
spanning trees and choices, and the product formula follows.
\end{proof}
\begin{figure}[tbp]
    \centering
    \begin{subfigure}{0.23\textwidth}
    \centering
    \includegraphics[width=1.0\linewidth]{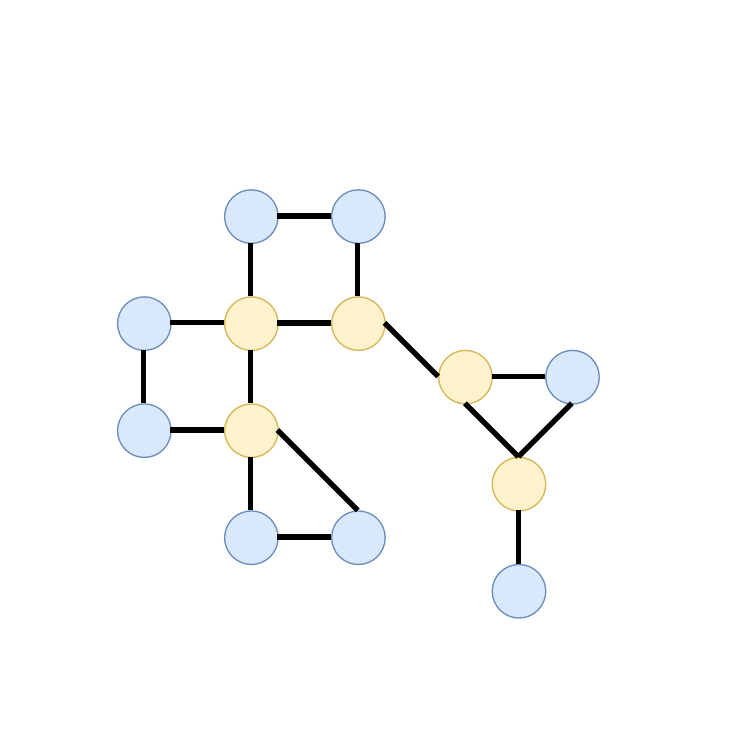}
    \caption{Original graph}
\end{subfigure}
\hfill
\begin{subfigure}{0.23\textwidth}
    \centering
    \includegraphics[width=1.0\linewidth]{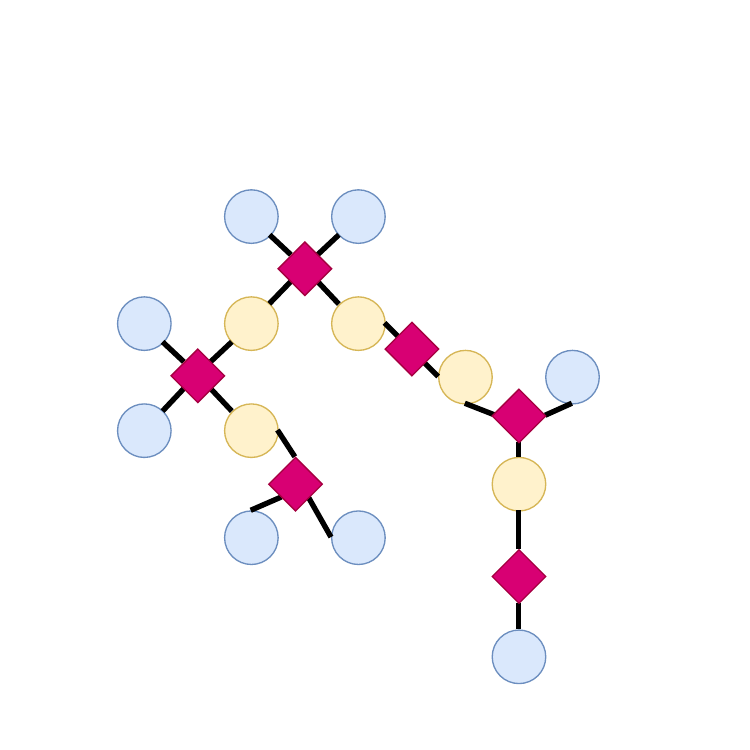}
    \caption{Factor graph }
\end{subfigure}
    \caption{An example of graph with cycle motifs (a) and its factor graph (b). The nodes in yellow are cut vertices. In the factor graph the red diamond nodes are the nodes representing the cycle motifs. The circle nodes are the nodes of the original graph. We can observe that the factor graph is bipartite.}
   \label{FAC1}
\end{figure}
This class of graphs is easy to handle because all cycles are local (inside motifs), and the tree-like factor
structure prevents “global” cycles from appearing when motifs are glued. An example of this kind of graph is presented in Fig.~\ref{FAC1} together with its factor graph.

Now, let $T$ be a uniformly random spanning tree of $G$ and $d_T(u,v)$ the distance between $u$ and $v$ inside $T$.
Using the indicator of whether an edge lies on the unique $u$--$v$ path in $T$,
\[
d_T(u,v)=\sum_{e\in E}\mathbf 1\{e\text{ is on the }u\leftrightarrow v\text{ path in }T\},
\]
and therefore
\[
\mathbb E[d_T(u,v)]
=\sum_{e\in E}\mathbb P\bigl(e\text{ is on the }u\leftrightarrow v\text{ path in }T\bigr).
\]
Fix an undirected edge $e=\{a,b\}$ and orient it as $e=(a,b)$. Let $N(u,a\mid b,v)$ be the number of spanning trees in
which the $u$--$v$ path uses $e$ in direction $a\to b$, and define $N(u,b\mid a,v)$ similarly. Define the signed current
\[
i_e:=\frac{N(u,a\mid b,v)-N(u,b\mid a,v)}{|\mathcal T|}.
\]

\begin{theorem}[Expected tree distance as an $\ell^1$ current]\label{SEC5TH2}
Under the same assumptions (cycle motifs, edge-disjoint, factor graph a tree), for any $u,v\in V$,
\[
\mathbb E[d_T(u,v)]=\sum_{e\in E}|i_e|.
\]
\end{theorem}

\begin{proof}
Define, for an oriented edge $e=(a,b)$, a signed indicator $X_e(T)\in\{-1,0,1\}$ equal to $+1$ if the $u$--$v$ path
in $T$ uses $e$ as $a\to b$, $-1$ if it uses it as $b\to a$, and $0$ otherwise. Then
$|X_e(T)|=\mathbf 1\{e \text{ is on the path}\}$, so
\[
\mathbb E[d_T(u,v)]=\sum_{e\in E}\mathbb E[|X_e(T)|].
\]
By uniformity, $\mathbb E[X_e(T)]=i_e$. Now use the tree-like motif structure. Because motifs are edge-disjoint and the factor graph is a tree, removing the
(unique) motif containing $e$ separates the global structure into two sides. If the $u$--$v$ path uses $e$, its
direction is determined by which side contains $u$ and which contains $v$. So the two orientations cannot both occur
among spanning trees: one of $N(u,a\mid b,v)$ and $N(u,b\mid a,v)$ is zero. Hence
\begin{align*}
|i_e|
=\frac{N(u,a\mid b,v)+N(u,b\mid a,v)}{|\mathcal T|}\\
=\mathbb P\bigl(e\text{ is on the }u\leftrightarrow v\text{ path}\bigr),    
\end{align*}
and summing over $e$ gives the result.
\end{proof}

Finally, let $L = D - A$ be the Laplacian matrix, and let $\varphi$ be a solution to
\[
L\varphi = b_{uv}\, ,
\]
which is defined up to an additive constant, and where $b_{uv}\in \mathbb{R}^n$ (with $n=|V|$) is the vector defined by
\(
(b_{uv})_u = 1\), \( (b_{uv})_v = -1\), and \( (b_{uv})_x = 0\)  for all other vertices $x$, i.e.\ $i_e$ is the potential drop across the edge ~\cite{BenjaminiLyonsPeresSchramm2001_UniformSpanningForests,DoyleSnell2000_RandomWalksElectricNetworks,Grimmett2012_ProbabilityOnGraphsNotes}.
Then, for an edge $e=(a,b)$, we have
\[
i_e = \varphi(b) - \varphi(a),
\]
 This means that the MFPT is given by
\begin{align}
    \mathrm{MFPT}^{\mathrm{LERW}}(u,v)=\sum_{(a,b)\in E}|\varphi(b)-\varphi(a)|.
\end{align}

This formula is fundamental because it turns a probabilistic quantity, the expected length of a loop-erased random walk, into a deterministic quantity computed from the graph potential. 
Rather than averaging over many random walks or summing over all possible self-avoiding paths, the $\mathrm{MFPT}^{\mathrm{LERW}}$ is obtained by adding the absolute potential drops along the edges. 
Thus the complexity of the problem is reduced from a path-based computation to an edge-based computation.
This gives a clear interpretation of the LERW MFPT: it measures the total amount of effective $u$--$v$ flow carried by the graph. 
Edges with larger values of $|\varphi(b)-\varphi(a)|$ are more relevant for the loop-erased trajectory, while edges with small potential drops contribute little. 
Consequently, the formula identifies which parts of the graph are most important for the effective passage from $u$ to $v$. It provides an exact and practical way to compute the MFPT, and it is particularly valuable for large graphs built from cycle motifs where enumerating all possible loop-erased paths would be computationally expensive.

We can now consider the case where all the motifs are complete graphs, i.e., cliques.
\begin{definition}
A motif $\sigma$ is a clique motif of size $r\ge 2$ if it induces a complete graph $K_r$
\[
|V(\sigma)|=r,
\qquad
E(\sigma)=\bigl\{\{x,y\}:x,y\in V(\sigma),~x\neq y\bigr\}.
\]
When $r=2$ this is just a single edge.
\end{definition}
\begin{figure}[tbp]
    \centering
    \begin{subfigure}{0.23\textwidth}
    \centering
    \includegraphics[width=1.0\linewidth]{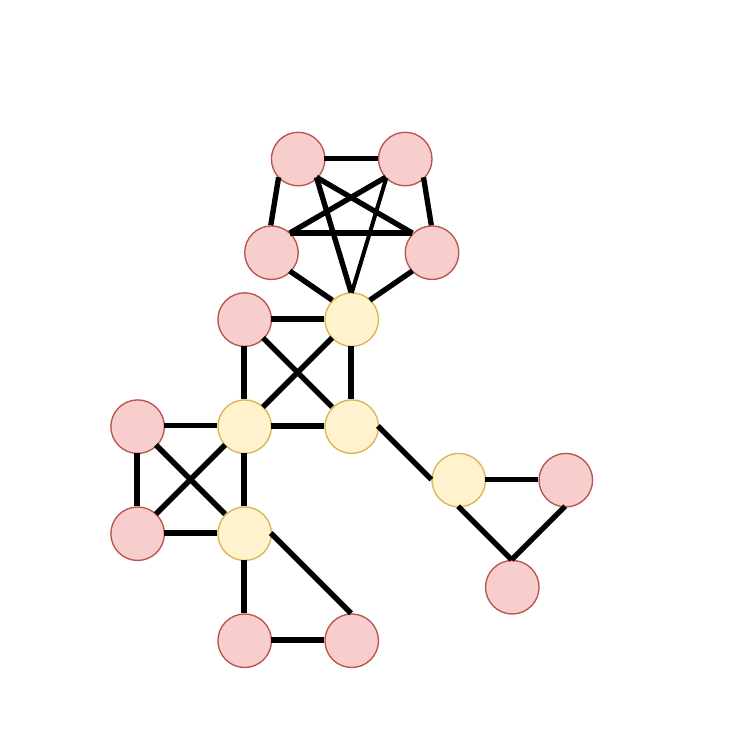}
    \caption{Original graph}
\end{subfigure}
\hfill
\begin{subfigure}{0.23\textwidth}
    \centering
    \includegraphics[width=1.0\linewidth]{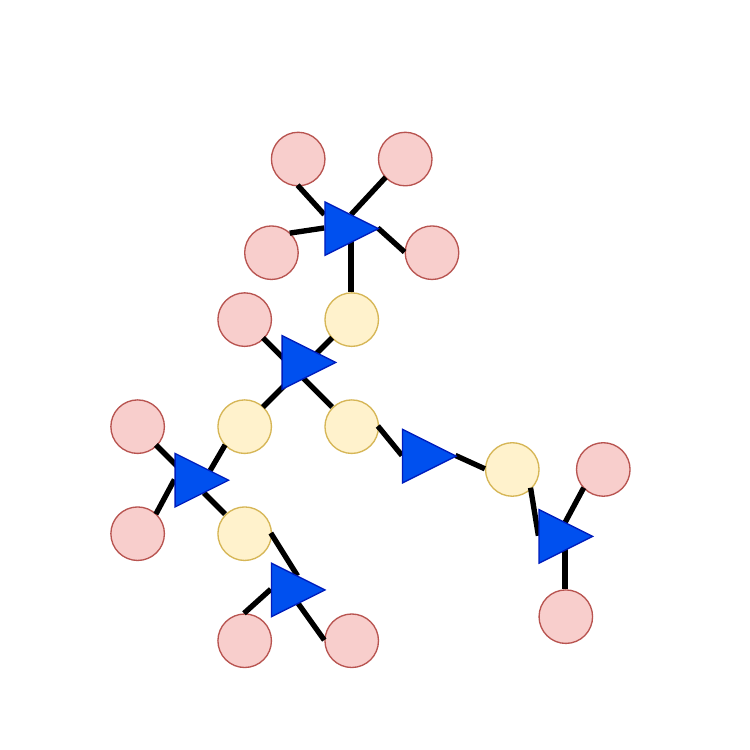}
    \caption{Factor graph}
\end{subfigure}
    \caption{An example of graph with clique motifs (a) and its factor graph (b). The nodes in yellow are cut vertices. In the factor graph the blue triangle nodes are the nodes representing the clique motifs. The circle nodes are the nodes of the original graph. We can observe that the factor graph is bipartite.}
    \label{FAC2}
\end{figure}
We consider graphs obtained as unions of clique motifs, assuming motifs are edge-disjoint and the factor graph
between vertices and motifs is a tree.  See panel (b) in Fig.~\ref{MOT} for few examples of clique motifs and Fig.~\ref{FAC2} for an example of graph obtained as union of clique motifs and its factor graph.

\begin{definition}[Tree of cliques / block graph]\label{SEC5DEF1}
A connected graph $G$ is a tree of cliques if
every motif is a clique $K_{r_\sigma}$,  motifs are edge-disjoint, and the vertex--motif factor graph is a tree.
\end{definition}

This implies that different cliques intersect in at most one vertex, and any cycle in $G$ is contained in a single
clique. For a vertex $u$, the degree decomposes as
\[
k_u=\sum_{\sigma\in S_u}k_u^{(\sigma)},
\qquad
k_u^{(\sigma)}=r_\sigma-1 \ \text{ if } \sigma\simeq K_{r_\sigma}.
\]
Also, each clique contributes $\binom{r_\sigma}{2}$ edges, so
\[
|E(G)|=\sum_{\sigma\in\mathcal M}\binom{r_\sigma}{2}.
\]

\begin{theorem}[Spanning trees of a tree of cliques]\label{SEC5TH3}
If $G$ is a tree of cliques, then
\[
|\mathcal{T}(G)|=\prod_{\sigma\in\mathcal M}|\mathcal{T}(K_{r_\sigma})|
=\prod_{\sigma\in\mathcal M}r_\sigma^{\,r_\sigma-2}.
\]
\end{theorem}

\begin{proof}
Inside each clique $\sigma\simeq K_{r_\sigma}$ there are $|\mathcal{T}(K_{r_\sigma})|=r_\sigma^{r_\sigma-2}$ spanning trees.
Choose independently one spanning tree inside each motif, and take their union. Because motifs attach through a tree of
shared cut-vertices, the union stays connected; because different motifs share at most one vertex, no cycle can be formed
by combining edges from different motifs. Hence the union is a spanning tree of $G$. Conversely, any spanning tree of $G$ restricts to a spanning tree inside each clique (there are no other edges), so it
induces exactly one local choice per motif. This is a bijection, giving the product formula.
\end{proof}

Let $V$ be the vertex set and let $\mathcal B$ be a family of motifs (blocks). Each block $b\in\mathcal B$ is a subgraph
on $V(b)\subseteq V$. We assume blocks are glued through shared cut-vertices so that the bipartite factor graph
\[
T=(V\cup\mathcal B,E),
\qquad (v,b)\in E\iff v\in V(b),
\]
is a tree. Then for any $u,v\in V$ there is a unique alternating path in $T$
\[
u=v_0-b_0-v_1-b_1-\cdots-b_{m-1}-v_m=v.
\]

We model the mean first passage time from $u$ to $v$ as the sum of the expected crossing times of the successive blocks
along that unique route
\begin{equation}
\mathrm{MFPT}^{\mathrm{LERW}}(u,v)
=\sum_{i=0}^{m-1}\tau_{b_i}(v_i\to v_{i+1}),
\label{SEC5EQ2}
\end{equation}
where $\tau_b(x\to y)$ is the MFPT inside block $b$ from entry vertex $x$ to exit vertex $y$ for the simple random
walk restricted to $b$. With the tree incidence structure, there is no ambiguity about which sequence of blocks is used,
so this additive decomposition is natural: each block contributes exactly one local passage.

Let $b=C_n$. For $x,y\in V(b)$ define the cyclic distance
\begin{equation}
k=d_{C_n}(x,y)=\min(\delta,\;n-\delta),
\label{SEC5EQ3}
\end{equation}
where $\delta$ is the distance along one chosen cyclic ordering. For a simple random walk on $C_n$,
\begin{equation}
\tau_{C_n}(x\to y)=\frac{2k(n-k)}{n}.
\label{SEC5EQ4}
\end{equation}

Let $b=K_n$. By symmetry, the MFPT from any $x$ to any $y\neq x$ depends only on $n$
\begin{equation}
\tau_{K_n}(x\to y)=\tau_{K_n}(n),
\qquad x\neq y.
\label{SEC5EQ5}
\end{equation}
A convenient closed form is
\begin{equation}
\tau_{K_n}(n)=\sum_{r=0}^{n-2}\frac{(n-2)!}{n^{r}(n-2-r)!}.
\label{SEC5EQ6}
\end{equation}

 If every block is a cycle $C_{n_i}$, then from \eqref{SEC5EQ2} and \eqref{SEC5EQ4},
\begin{equation}
\mathrm{MFPT}^{\mathrm{LERW}}(u,v)
=\sum_{i=0}^{m-1}\frac{2k_i(n_i-k_i)}{n_i},
\label{SEC5EQ7}
\end{equation}
where $n_i=|V(b_i)|$ and $k_i=d_{C_{n_i}}(v_i,v_{i+1})$.

If every block is a clique $K_{n_i}$, then from \eqref{SEC5EQ2} and \eqref{SEC5EQ6},
\begin{equation}
\mathrm{MFPT}^{\mathrm{LERW}}(u,v)
=\sum_{i=0}^{m-1}\sum_{r=0}^{n_i-2}\frac{(n_i-2)!}{n_i^{r}(n_i-2-r)!}.
\label{SEC5EQ8}
\end{equation}

\begin{figure}
    \centering
    \includegraphics[width=1\linewidth]{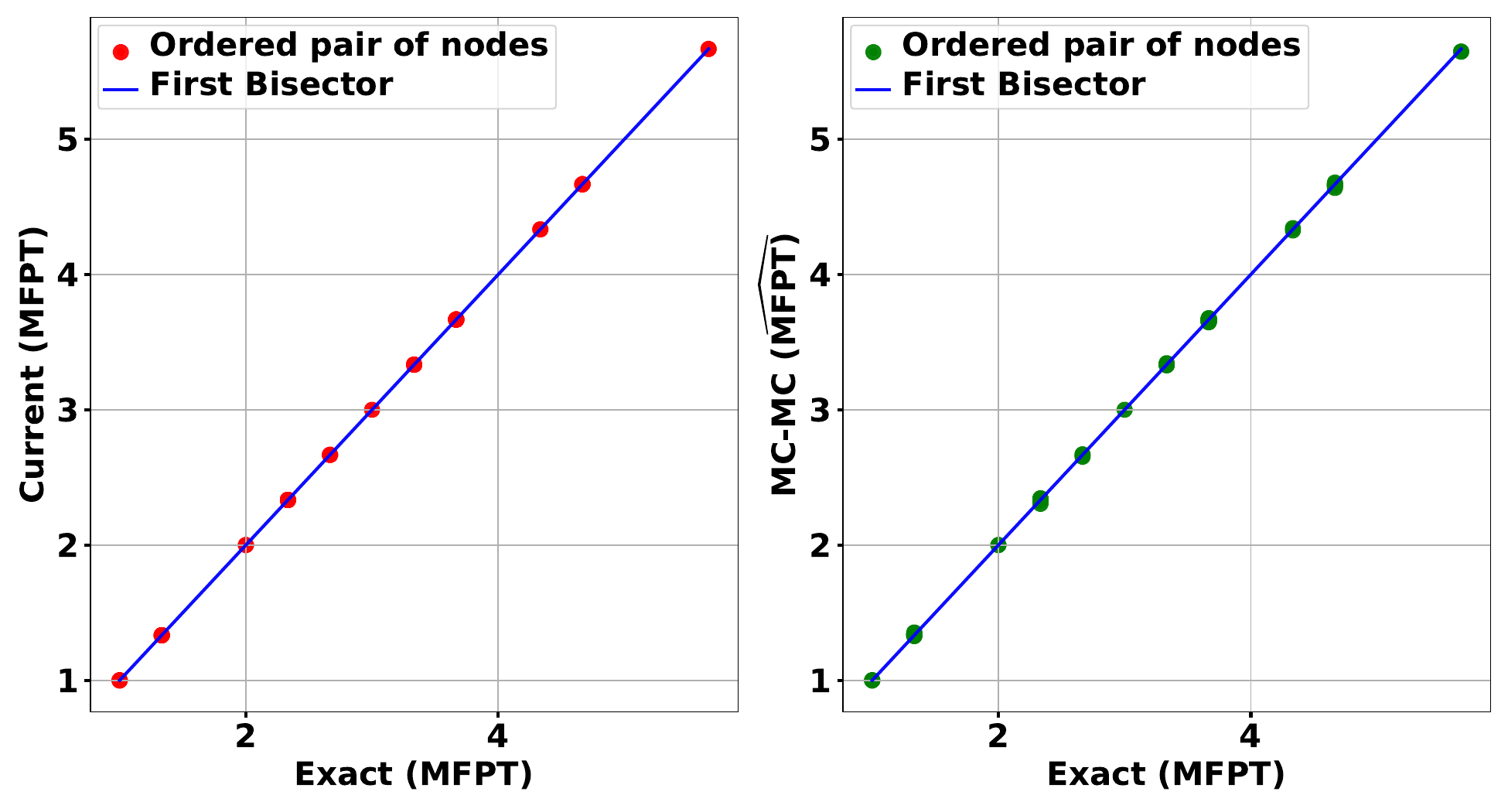}
    \caption{Mean First Passage Time of Loop Erased Random Walks between different pairs of nodes of a graph built with cycle motifs. On the left panel, on the y-axis, we have the MFPT obtained with the formula of the theorem \ref{SEC5TH2} and on the x-axis, we have the MFPT obtained with  the formula \eqref{SEC5EQ7}. On the right panel, on the y-axis,  we have the MFPT obtained with an MC-MC approach with empirical average and on the x-axis, we have again the MFPT obtained with the formula \eqref{SEC5EQ7}.}
    \label{MOTIF}
\end{figure}

The factor graph being a tree forces a unique block sequence from $u$ to $v$. The global MFPT is then just the sum of
the expected times needed to move from the entry cut-vertex to the exit cut-vertex inside each visited motif.
We consider a small graph built on cyle motifs for which the mean first-passage time between any ordered pair of nodes can be computed exactly. We evaluate the MFPT for all distinct node pairs with three different approaches. The Exact one corresponds to the factor graph approach with the formula \eqref{SEC5EQ7}. The Current method represents the MFPT computed with the formula of theorem \ref{SEC5TH2}. Finally, the MC–MC curve is obtained by Monte Carlo simulation, where the MFPT is estimated as the empirical average of many loop-erased random walk lengths. Fig.~\ref{MOTIF} displays the MFPT values for all node pairs. We observe that the points are on the first bisector, showing that the current implementation exactly reproduces the analytical MFPT and is consistent with Monte Carlo estimates. The maximum deviations in the analytical and Monte Carlo panels are $\varepsilon_{\max}^{\mathrm{left}}=0$ and $\varepsilon_{\max}^{\mathrm{right}}=3.8\times 10^{-2}$, respectively. This supports the numerical consistency of the methods on elementary cycle motifs.

\begin{figure}
    \centering 
    \begin{subfigure}{0.23\textwidth}
    \centering
    \includegraphics[width=1.0\linewidth]{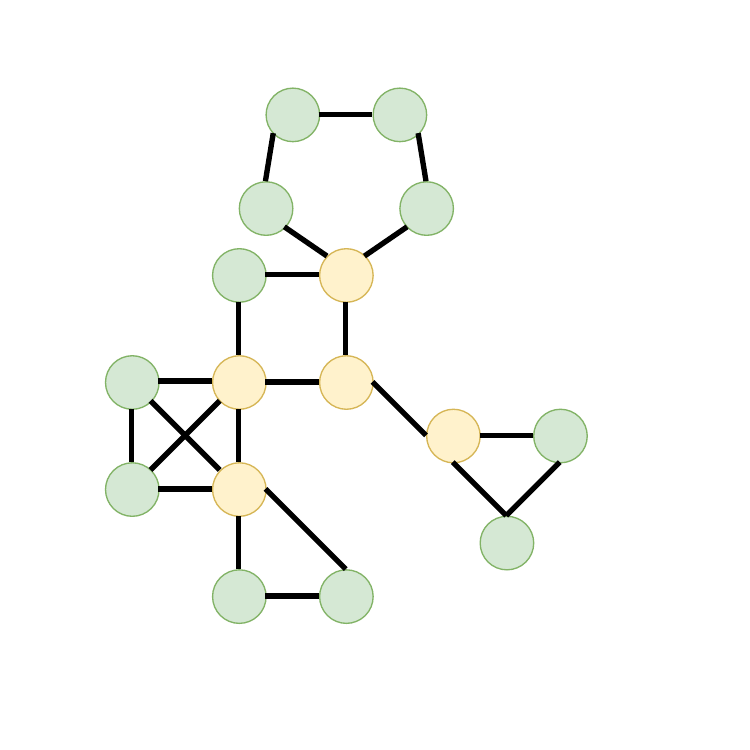}
    \caption{Original graph}
\end{subfigure}
\hfill
\begin{subfigure}{0.23\textwidth}
    \centering
    \includegraphics[width=1.0\linewidth]{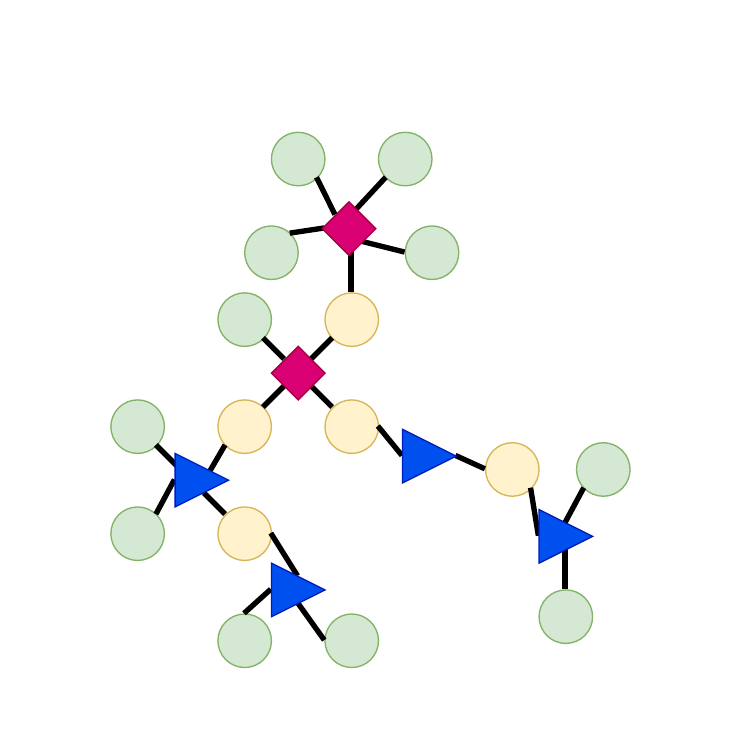}
    \caption{Factor graph}
\end{subfigure}
     \caption{An example of graph with both cycle and clique motifs (a) and its factor graph (b). The nodes in yellow are cut vertices. In the factor graph the blue triangle nodes are the nodes representing the clique motifs and  the red diamond nodes are the nodes representing the cycle motifs. The clique motifs of three nodes are the same as the cycle motifs of three node,s so we could have use the red diamond to represent them. This is the same for edge motif. The circle nodes are the nodes of the original graph. We can observe that the factor graph is bipartite.}
    \label{FAC3}
\end{figure}

\begin{figure}
    \centering
    \includegraphics[width=1\linewidth]{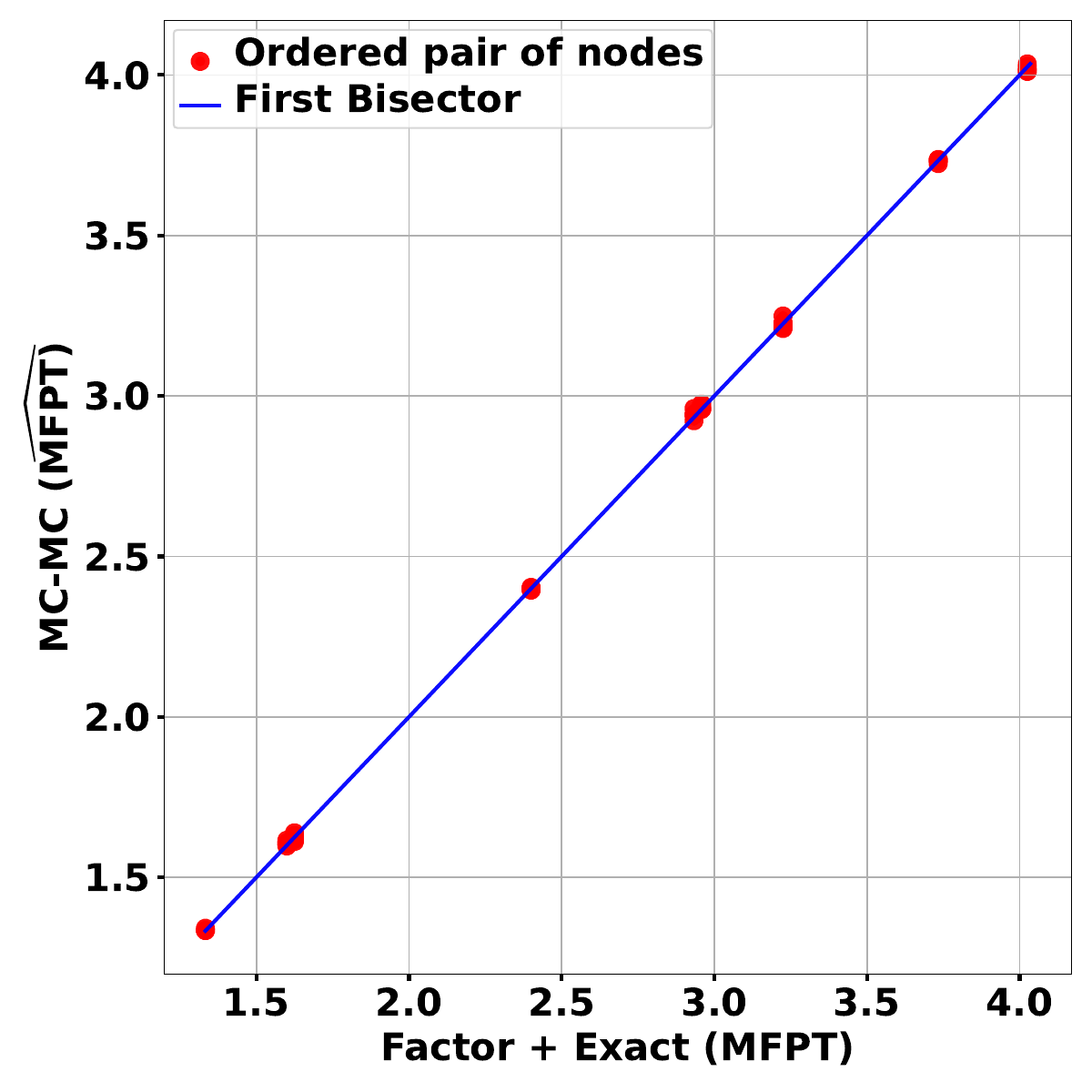}
    \caption{Mean First Passage Time of Loop Erased Random Walks between the different pairs of nodes of a graph built by gluing cycle and clique motifs. On the x-axis, we have the MFPT obtained with the global formula\eqref{SEC5EQ2} with the use of local formulas and factorization and on the y-axis, we have the MFPT obtained with an MC-MC approach with empirical average.}
    \label{MFPT3}
\end{figure}

We now apply the proposed factorization-based method to a graph built with cycle and clique motifs. See an example of this kind of graph in Fig.~\ref{FAC3}. 
For a set of distinct ordered node pairs, we compute the mean first-passage time (MFPT) of loop-erased random walks. The MC–MC curve is obtained by direct Monte Carlo simulation. The Factor curve corresponds to the MFPT computed using the factor-tree decomposition of the graph, where the total MFPT between two nodes is expressed as a sum of local contributions along the unique path between their corresponding blocks \eqref{SEC5EQ6},\eqref{SEC5EQ4}. Figure~\ref{MFPT3} shows the MFPT values for all considered node pairs. The factor-based MFPT follows the MFPT find with the MC-MC closely, with $\varepsilon_{\max}=1.6\times10^{-2}$. This result validates the accuracy of the factorization approach on graphs of intermediate size while significantly reducing the computational effort compared to full Monte Carlo estimation.
% \begin{figure}[tbp]
%     \centering
%     \includegraphics[width=1\linewidth]{FIG/MOTIF2.pdf}
%     \caption{Mean First Passage Time of Loop Erased Random Walks between different pairs of nodes of a graph built with cycle motifs. On the y-axis, we have the MFPT obtained with the formula \ref{SEC5TH2} and on the x-axis, we have the MFPT obtained with an MC-MC approach with empirical average.}
%     \label{MOTIF}
% \end{figure}

Now we will establish formulas for the distributions of erased loops for motif based networks. By using the results obtained in the previous sections, we now derive a simplified expression for the distribution of erased loops in networks constructed from motifs. In the motif-based framework introduced in this section, the graph $G$ is constructed by gluing together a family of motifs $\sigma \in \mathcal{M}$. Each motif $\sigma$ contributes a set of internal edges $E(\sigma)$ and the full edge set satisfies
\[
E(G)=\bigcup_{\sigma\in\mathcal{M}} E(\sigma).
\]

We assume that motifs are edge-disjoint and that the associated vertex--motif factor graph is a tree. Under this structural assumption, any cycle in $G$ must lie entirely inside a single motif. In particular, loops created by the loop-erased random walk cannot involve edges belonging to multiple motifs. This property allows us to decompose the global loop distribution into independent contributions coming from the different motifs.

For each motif $\sigma$, we define a local distribution $\pi_{\sigma}(\ell)$, representing the proportion of loops of size $\ell$ that are generated inside the motif $\sigma$. These local distributions can be computed explicitly for the motifs studied in the previous sections.
In the spanning-tree edge-sampling representation, an edge is chosen
uniformly from $E(G)$, so the chance that it belongs to motif $\sigma$
is $|E(\sigma)|/|E(G)|$. However, the local distribution
$\pi_\sigma(\ell)$ is already conditioned on the occurrence of a loop.
Therefore, edge counts alone are not the correct weights for mixing
these conditional distributions. Define the local unnormalized event
probabilities and loop-creation probabilities by
\[
P_\sigma(\ell)=\frac{P_{\mathrm{st},\sigma}(\ell)}{\ell},
\qquad q_\sigma=\sum_{k\ge2}P_\sigma(k),
\qquad \pi_\sigma(\ell)=\frac{P_\sigma(\ell)}{q_\sigma}.
\]
Because the spanning-tree measure factorizes across the blocks of a
tree-like factor graph, the global unnormalized distribution is
\[
P_G(\ell)=\sum_{\sigma\in\mathcal M}
\frac{|E(\sigma)|}{|E(G)|}P_\sigma(\ell).
\]
Normalizing only after summing the contributions yields
\begin{equation}\label{SEC5EQ9}
\pi_G(\ell)=
\frac{\displaystyle\sum_{\sigma\in\mathcal M}
|E(\sigma)|q_\sigma\pi_\sigma(\ell)}
{\displaystyle\sum_{\sigma\in\mathcal M}
|E(\sigma)|q_\sigma}.
\end{equation}
The simpler edge-weighted average is valid only if all $q_\sigma$
coincide. Equation~\eqref{SEC5EQ9} applies to the explicitly defined
spanning-tree loop-size model; comparison to dynamically erased loops
requires the correspondence assumed in Theorem~\ref{SEC2TH1}.

Suppose that the graph $G$ is constructed exclusively from cycle motifs. Let
\[
\mathcal{M} = \{C_{n_1},C_{n_2},\dots,C_{n_m}\}
\]
be the set of cycle motifs, where $C_{n_i}$ denotes a cycle of length $n_i$.
Since each cycle contributes exactly $n_i$ edges, we have
\[
|E(C_{n_i})| = n_i.
\]

For a cycle $C_{n_i}$, the local loop-creation probability is
\[
q_{C_{n_i}}=\frac{n_i-1}{2n_i}+\frac1{n_i^2}.
\]
Using the corrected conditional-mixture weights gives
\begin{equation}
\pi(\ell)=
\sum_{i=1}^{m}
\frac{n_i q_{C_{n_i}}}{\sum_{j=1}^{m} n_j q_{C_{n_j}}}
\,\pi_{C_{n_i}}(\ell).
\end{equation}

In particular, since cycles only generate loops of size $2$ or $n_i$, the global probabilities are
\begin{align}
\pi(2) &=
\sum_{i=1}^{m}
\frac{n_i q_{C_{n_i}}}{\sum_{j=1}^{m} n_j q_{C_{n_j}}}
\,\frac{n_i^2-n_i}{n_i^2-n_i+2}, \\
\pi(\ell)\big|_{\ell\ge 3} &=\sum_{i=1}^{m}
\frac{n_i q_{C_{n_i}}}{\sum_{j=1}^{m} n_j q_{C_{n_j}}}
\,\frac{2}{n_i^2-n_i+2}\,\mathbf{1}_{\{\ell=n_i\}}.
\end{align}
Thus, in a network composed only by cycles, erased loops correspond either to local backtracking events ($\ell=2$) or to the traversal of the entire cycle motif.

Assume now that the network is constructed exclusively from clique motifs
\[
\mathcal{M}=\{K_{r_1},K_{r_2},\dots,K_{r_m}\},
\]
where $K_{r_i}$ denotes a clique with $r_i$ vertices.

Each clique contributes
\[
|E(K_{r_i})|=\binom{r_i}{2}=\frac{r_i(r_i-1)}{2}
\]
edges to the network.

From Section~3, the loop-size distribution inside a clique $K_n$ is given by
\[
\pi_{K_n}(\ell)=
\frac{a_n(\ell)}{\sum_{k\ge2} a_n(k)},
\]
where
\[
a_n(\ell)=\frac{(n)_\ell}{n^\ell}
=
\prod_{i=0}^{\ell-1}\left(1-\frac{i}{n}\right).
\]

Here the local event probability is
\[
q_{K_r}=\frac{1}{r-1}\sum_{k=2}^{r}\frac{(r)_k}{r^k}.
\]
Using the corrected conditional-mixture weights, the global distribution becomes
\begin{equation}
\pi(\ell)
=
\sum_{i=1}^{m}
\frac{\binom{r_i}{2}q_{K_{r_i}}}{\sum_{j=1}^{m}\binom{r_j}{2}q_{K_{r_j}}}
\,\pi_{K_{r_i}}(\ell).
\end{equation}
\begin{figure}
    \centering
    \includegraphics[width=1\linewidth]{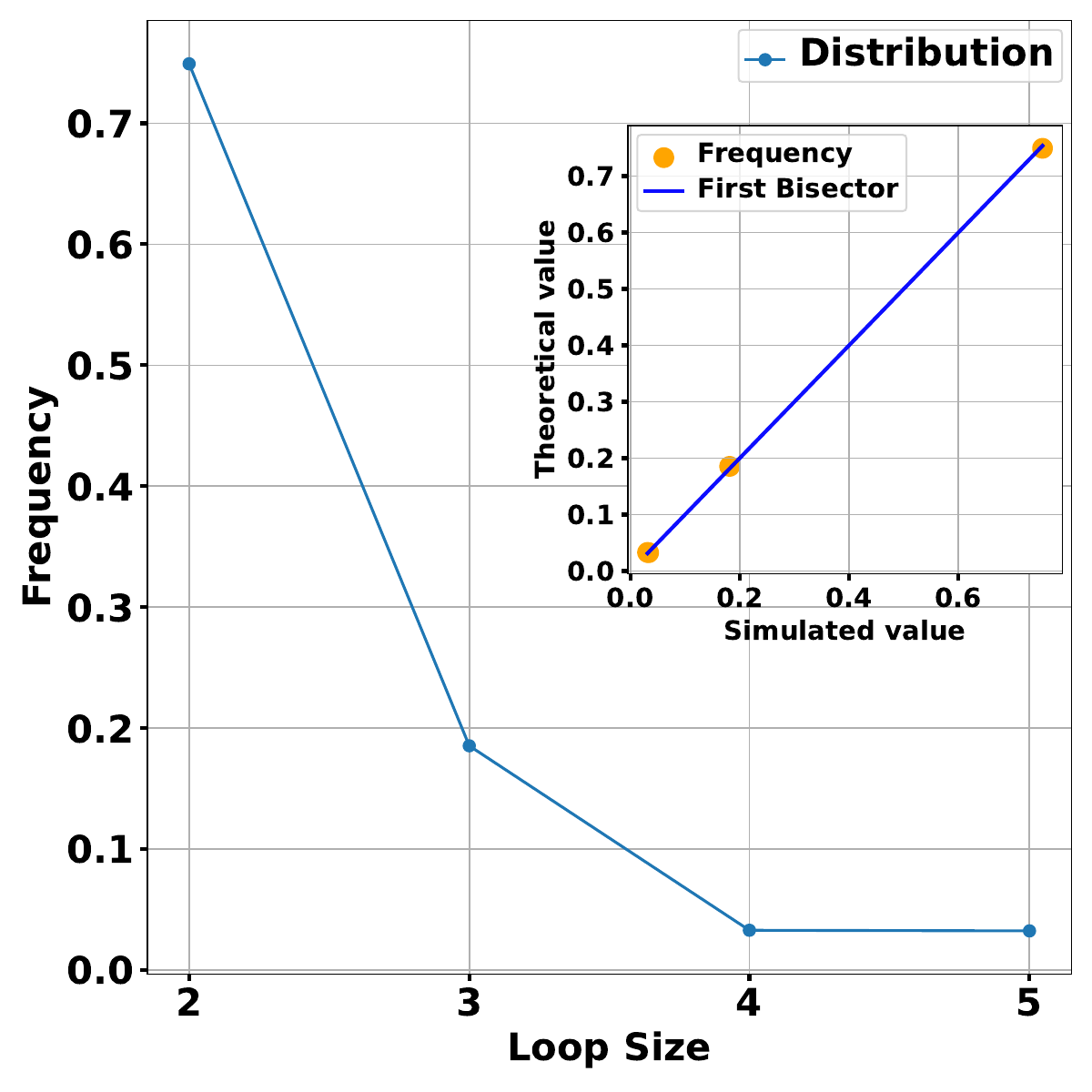}
    \caption{Main panel shows the distribution of erased loops obtained with the relation \eqref{SEC5EQ9}. Inset shows the comparison between Frequencies obtained by the relation \eqref{SEC5EQ9} and Frequencies observed with true loop erase random walk. On the y-axis, we have the frequency obtained with the exact formula \eqref{SEC5EQ9} and on the x-axis, we have the frequency obtained by counting size of erased loops.}
    \label{LOOPMOTIF1}
\end{figure}
Substituting the explicit expression of $\pi_{K_{r_i}}(\ell)$ gives
\begin{equation}
\pi(\ell)
=
\sum_{i=1}^{m}
\frac{\binom{r_i}{2}q_{K_{r_i}}}{\sum_{j=1}^{m}\binom{r_j}{2}q_{K_{r_j}}}
\,
\frac{a_{r_i}(\ell)}{\sum_{k\ge2} a_{r_i}(k)},
\end{equation}

with
\[
a_{r_i}(\ell)=
\prod_{j=0}^{\ell-1}\left(1-\frac{j}{r_i}\right).
\]

Hence, when the network is composed entirely of clique motifs, the loop-size distribution is a weighted mixture of the birthday-type distributions generated inside each clique.

We used formula~\eqref{SEC5EQ9} to compute the distribution of erased loops on a graph buit by gluing cliques and cycles. The results are presented in the main panel of Fig.~\ref{LOOPMOTIF1}. To validate the accuracy of the formula, we compared these theoretical results with empirical measurements obtained by directly counting the sizes of erased loops generated by loop-erased random walks. The comparison is shown in the inset of Fig.~\ref{LOOPMOTIF1}. The largest difference in the inset of Fig.~\ref{LOOPMOTIF1} is $\varepsilon_{\max}=5.8\times10^{-2}$. Therefore, the proposed formula accurately captures the distribution of erased loop sizes on this kind of graphs. This comparison provides a numerical check of the proposed spanning-tree model; the exact agreement with the dynamic erasure process requires the correspondence discussed in Appendix~\ref{APPA}.

\section{Networks with tree-like factor graph }

\begin{figure}[tbp]
\centering
\begin{subfigure}{0.23\textwidth}
    \centering
    \includegraphics[width=1.0\linewidth]{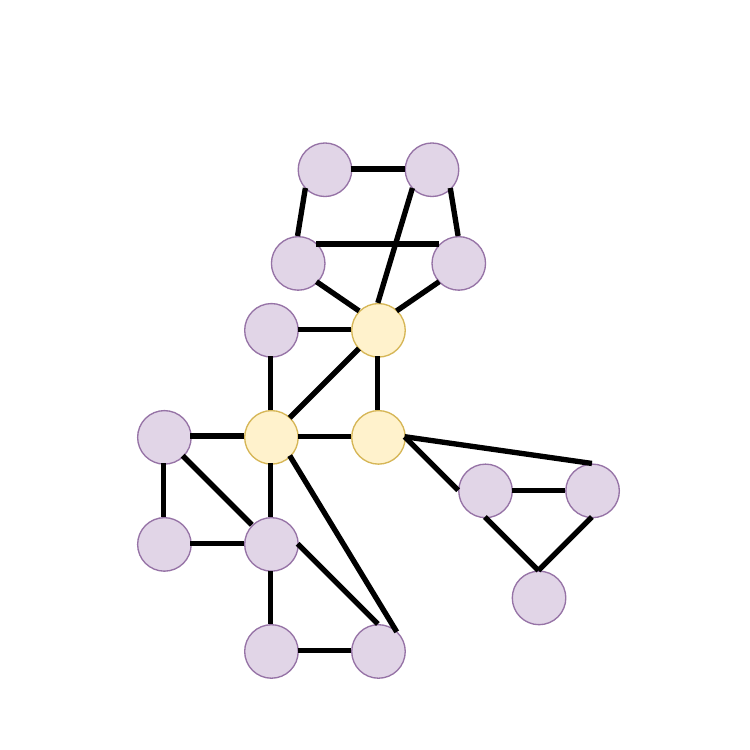}
    \caption{Original graph}
\end{subfigure}
\hfill
\begin{subfigure}{0.23\textwidth}
    \centering
    \includegraphics[width=1.0\linewidth]{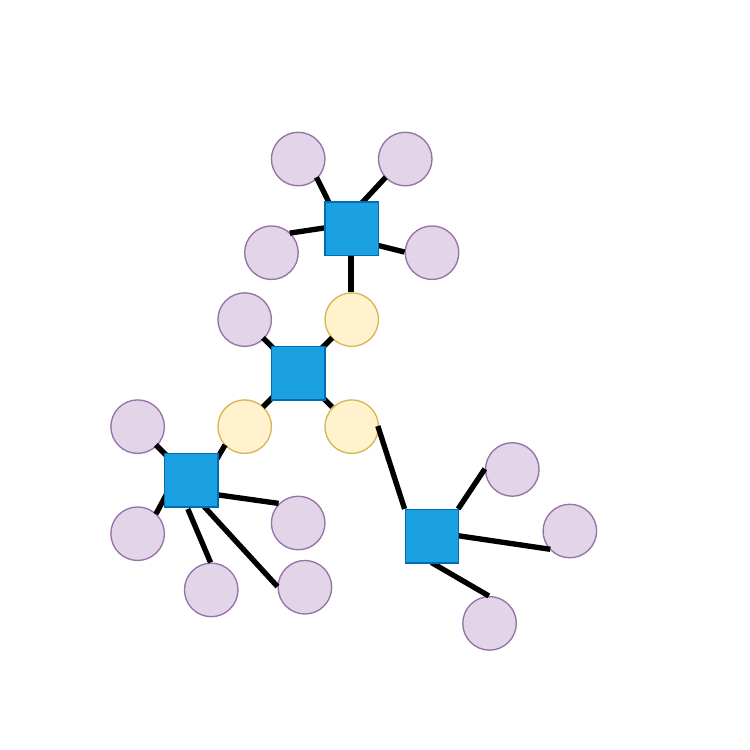}
    \caption{Factor graph}
\end{subfigure}
 \caption{An example of graph where motifs are not clique or cycle (a) and its factor graph (b). The nodes in yellow are cut vertices. In the factor graph the blue square  nodes are the nodes representing the motifs. The circle nodes are the nodes of the original graph. We can observe that the factor graph is bipartite.}
\label{FAC4}
\end{figure}

This section extends the motif-based framework to general networks by using their factor graph decomposition. See Fig.~\ref{FAC4} where the motifs are not necessarily cycles or cliques.
The approach captures the mean first passage structure of a more complex graph through its
biconnected components ~\cite{West1996IntroductionGraphTheory}.
\begin{definition}
 Let $G=(V,E)$ be a finite, connected, undirected graph.
A vertex $a\in V$ is called a cut vertex (or articulation point)  ~\cite{Tarjanrober1972,West1996IntroductionGraphTheory}
if removing it disconnects the graph, i.e.
\[
\kappa(G-a)>\kappa(G),
\]
where $\kappa(\cdot)$ denotes the number of connected components.   
\end{definition}
We can see cut vertices in the yellow color in Fig.~\ref{FAC1}, Fig.~\ref{FAC2}, Fig.~\ref{FAC3} and Fig.~\ref{FAC4}.
\begin{definition}
 A block (or biconnected component) ~\cite{Tarjanrober1972,West1996IntroductionGraphTheory} of $G$ is a maximal subgraph $b$
that is $2$-vertex-connected (no cut vertex).
Equivalently, $b$ is maximal among subgraphs in which every pair of vertices
has two internally vertex-disjoint paths.   
\end{definition}

Let $\mathcal{B}$ be the family of all blocks of $G$ and $V(b)\subseteq V$ the vertex set of block $b$.
The way blocks are attached is encoded by the bipartite factor graph
\[
T=\bigl(V\cup\mathcal{B},\,E_T\bigr), \qquad (v,b)\in E_T\iff v\in V(b).
\]
Nodes of $T$ are of two types: vertices of $G$ and block-nodes.
An edge $(v,b)$ indicates that vertex $v$ belongs to block $b$.
When $G$ is connected, $T$ is a tree—connected and acyclic known as the factor graph. Let $\mathrm{Art}(G)$ denotes the set of all cut vertices of $G$.
In this representation
\begin{align*}
  x\in \mathrm{Art}(G)
&\iff
\bigl|\{\,b\in\mathcal{B}:\ x\in V(b)\,\}\bigr|\ge2\\
&\iff
\deg_T(x)\ge2.
\end{align*}
Such vertices are called the attachment points of the decomposition.

In practice to construct $T$,
one first computes the set of biconnected components $\mathcal{B}$ of $G$
(e.g. via a DFS-based algorithm such as Tarjan's algorithm ~\cite{Tarjanrober1972}). Then one builds the bipartite
graph $T$ by creating one node for every $v\in V$ and one node for every block $b\in\mathcal{B}$,
and adding an edge $(v,b)$ whenever $v\in V(b)$.
This produces exactly the factor tree representation described above. Assume now that $G$ is decomposed into blocks $\mathcal{B}$ and that the factor graph
$T=(V\cup\mathcal{B},E_T)$ is a tree. Fix two distinct vertices $u,v\in V$.
Because $T$ is a tree, there is a unique alternating path
\[
u=v_0 - b_0 - v_1 - b_1 - \cdots - b_{m-1} - v_m=v
\quad \text{in } T,
\]
which identifies the unique sequence of blocks that must be crossed to go from $u$ to $v$,
and the intermediate vertices $v_1,\dots,v_{m-1}$ are precisely the attachment (cut) vertices ~\cite{ferreira2012optimal}. The procedure is described in Algorithm~\ref{ALGO4}. In addition, for the simple case in which blocks are motifs (e.g., cycles or cliques), we introduce Algorithm~\ref{ALGO2} and its reverse counterpart, Algorithm~\ref{ALGO3}.

For each block $b\in\mathcal{B}$, define a local mean first passage time
\[
\tau_b(x\to y),
\qquad x,y\in V(b),
\]
as the MFPT to go from $x$ to $y$ for the random-walk model restricted to the subgraph
induced by $V(b)$ (with the same transition structure/weights as in $G$).
This local quantity can be computed exactly for special motif families (cycles, cliques),
or numerically/algorithmically for general blocks (e.g. via spanning-tree/LERW kernels).

In the tree-like setting, the global MFPT from $u$ to $v$ can be modeled as the sum of
the successive expected crossing times of the blocks along the unique route in $T$
\begin{equation}\label{SEC6EQ1}
\mathrm{MFPT}^{\mathrm{LERW}}(u,v)
\;=\;
\sum_{i=0}^{m-1} \tau_{b_i}(v_i\to v_{i+1}).
\end{equation}
The intuition is that the factor tree structure removes any ambiguity about the order in which
blocks are visited: there is exactly one block sequence from $u$ to $v$ in $T$.
Hence the passage from $u$ to $v$ decomposes into $m$ local passages, each occurring
inside a single block between two well-defined boundary vertices (entry $v_i$ and exit $v_{i+1}$).

This decomposition provides a natural generalization of the motif-based framework:
each block plays the role of a “super-motif,” and the factor tree specifies how these local structures
connect hierarchically.  
The additivity in \eqref{SEC6EQ1} follows from the acyclic nature of $T$,
which guarantees that each block contributes exactly one traversal term and no revisiting occurs.  
This principle underlies all higher-level MFPT computations in modular or tree-like networks. This factorized approach reduces the effective dimensionality of the problem by working at the level of blocks rather than individual elements, which significantly lowers computational complexity for large networks. 

\begin{figure}
    \centering
    \includegraphics[width=1\linewidth]{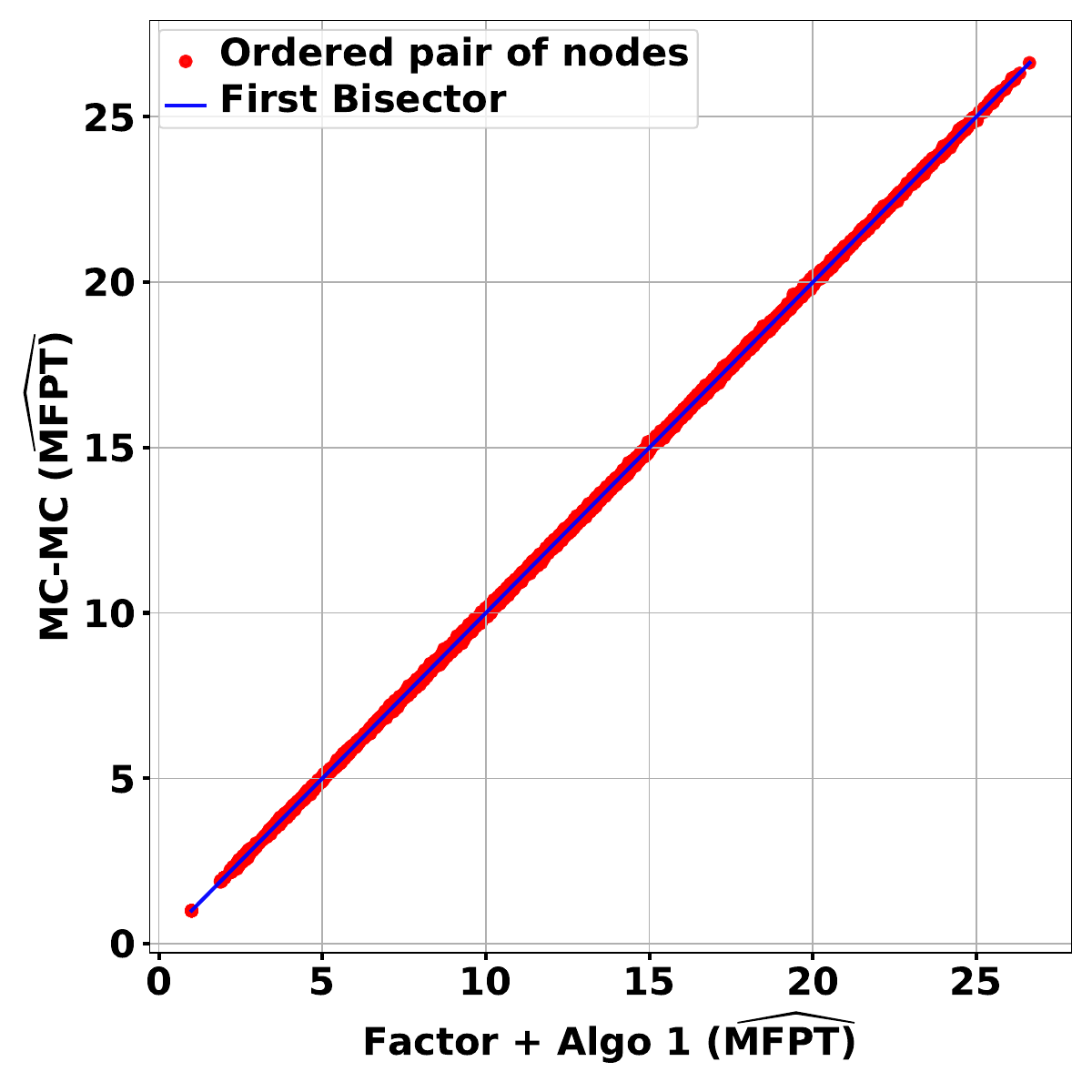}
    \caption{Mean First Passage Time of Loop Erased Random Walks between different pairs of nodes of a graph built by gluing Watts–Strogatz communities: $10$ communities of size $12$ (total $120$ nodes), each with mean degree $k=4$ and rewiring probability $p=0.1$, connected through gateway nodes. On the x-axis, we have the MFPT obtained with the algorithm \ref{alg:1} with factorization approach and on the y-axis, we have the MFPT obtained with an MC-MC approach with empirical average.}
    \label{MFPT4}
\end{figure}
\begin{figure}
    \centering
    \includegraphics[width=1\linewidth]{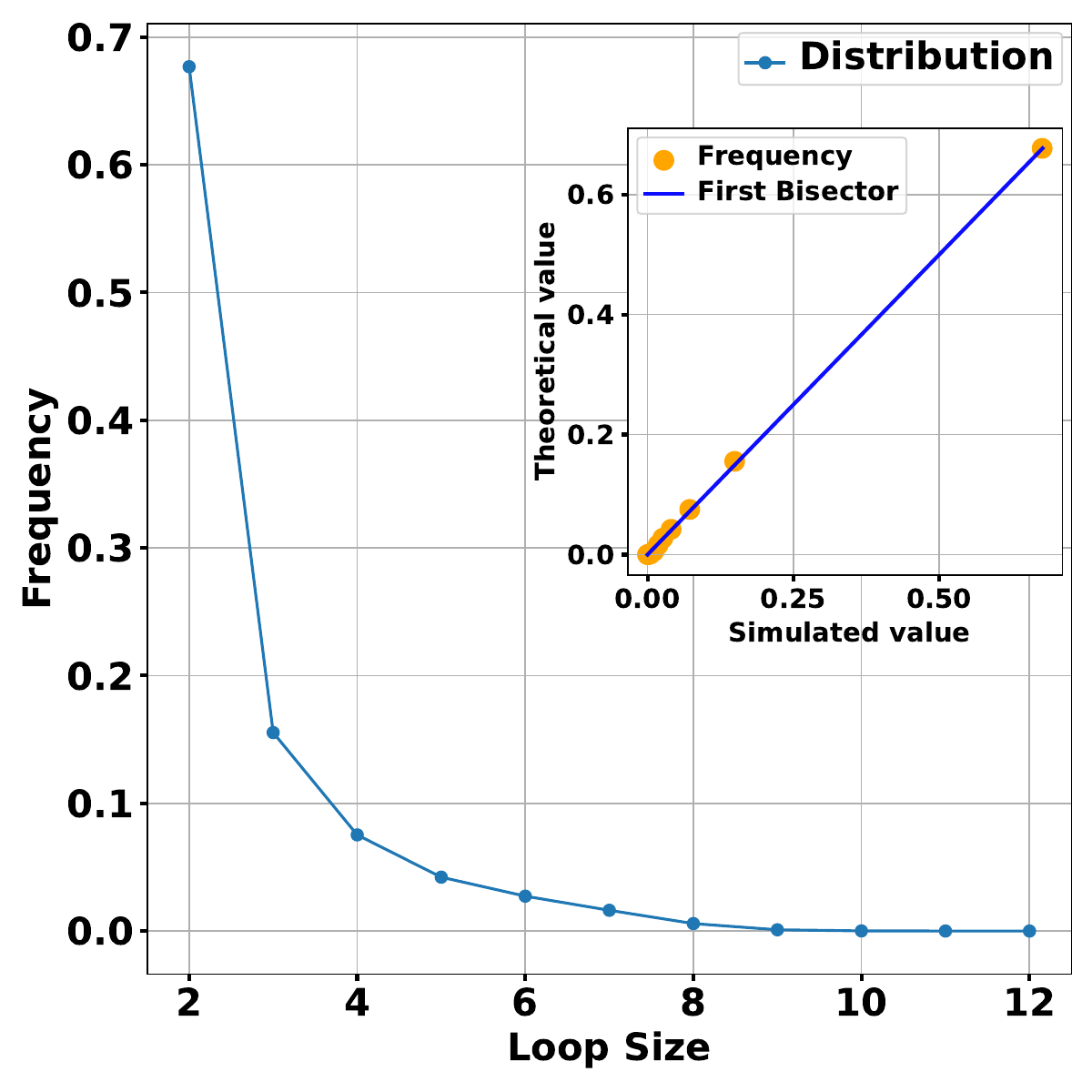}
    \caption{Main panel shows the distribution of erased loops obtained with the relation \eqref{SEC6EQ2} on a graph by gluing Watts–Strogatz communities: $10$ communities of size $12$ (total $120$ nodes), each with mean degree $k=4$ and rewiring probability $p=0.1$, connected through gateway nodes. Inset shows the comparison between Frequencies obtained by the relation \eqref{SEC6EQ2} and Frequencies observed with true loop erase random walk on this graph. On the y-axis, we have the frequency obtained with the exact formula \eqref{SEC5EQ9} and on the x-axis, we have the frequency obtained by counting size of erased loops.}
    \label{LOOPFACTOR1}
\end{figure}
We generate a graph by gluing Watts–Strogatz communities: $10$ communities of size $12$ (total $120$ nodes), each with mean degree $k=4$ and rewiring probability $p=0.1$, connected through gateway nodes. For every ordered node pair $(u,v)$, we estimate the MFPT in two ways. First, the MC-MC values are obtained by Monte Carlo simulation: for each pair, we sample many random walks from $u$ until hitting $v$, apply loop-erasure, and approximate the MFPT by the empirical mean of LERW lengths (with a step cap). Second, the Factor values uses the proposed factorization approach: the graph is decomposed into a factor tree (block–articulation structure), and the MFPT between $u$ and $v$ is computed as a sum of local contributions along the unique path in this tree. Inside each block, these local expectations are computed exactly via the transfer matrix approach presented in the first section. Fig.~\ref{MFPT4} plots the MFPT values for the first $50$ distinct node pairs (indexed on the x-axis). The factorization result closely follows the one obtained with the MC-MC ($\varepsilon_{\max}=2.5\times 10^{-1}$).

For the same reason as the previous section the global distribution of erased loops can be written as a weighted mixture of the local bloc distributions,
\begin{equation}\label{SEC6EQ2}
\pi_G(\ell)
=\frac{\displaystyle\sum_{b\in\mathcal B}
|E(b)|q_b\pi_b(\ell)}
{\displaystyle\sum_{b\in\mathcal B}|E(b)|q_b},
\qquad q_b=\sum_{k\ge2}P_b(k).
\end{equation}

We used formula~\eqref{SEC6EQ2} to compute the distribution of erased loops on a graph by gluing Watts–Strogatz communities: $10$ communities of size $12$ (total $120$ nodes), each with mean degree $k=4$ and rewiring probability $p=0.1$, connected through gateway nodes. The results are presented in the main panel of Fig.~\ref{LOOPFACTOR1}. To validate the accuracy of the formula, we compared these theoretical results with empirical measurements obtained by directly counting the sizes of erased loops generated by loop-erased random walks. The comparison is shown in the inset of Fig.~\ref{LOOPFACTOR1}. For the inset of Fig.~\ref{LOOPFACTOR1}, we obtain $\varepsilon_{\max}=9.6\times 10^{-3}$. Therefore, the proposed formula accurately captures the distribution of erased loop sizes on this kind of graphs. This comparison provides a numerical check of the proposed spanning-tree model; the exact agreement with the dynamic erasure process requires the correspondence discussed in Appendix~\ref{APPA}.

 \section{Conclusion}

In this work, we studied loop-erased random walks on finite graphs through their connection with uniform spanning trees. 
The goal was to describe two complementary aspects of the loop-erasure process. 
The first is the mean first passage time (MFPT) of the loop-erased random walk, which measures the length of the effective self-avoiding path that remains after all loops have been removed. 
The second is the distribution of erased loop sizes, which describes the redundant excursions produced by the random walk before they are erased. 
Together, these two quantities give a more complete description of random exploration on graphs: one describes the part of the trajectory that is kept, while the other describes the part that is removed.

The first contribution of this work is a general method for computing the LERW MFPT on arbitrary finite graphs. 
Using the transfer-current kernel, we expressed the probability of each simple path in terms of spanning-tree probabilities. 
This provides an exact framework in which the MFPT is computed from theoretical path weights, rather than only from empirical averages. 
The numerical comparisons with Monte Carlo simulations show very good agreement, both for unweighted and weighted graphs. 
This confirms that the transfer-current approach gives a reliable way to quantify effective random-walk exploration.
The second contribution is the derivation of a general formula for the distribution of erased loop sizes. 
This result shows that the probability of erasing a loop of a given size is also determined by spanning-tree probabilities. 
Thus, the same mathematical framework describes both sides of the loop-erasure process: the effective path that remains and the loops that are removed. 
This is important because these two observables contain different information. 
The MFPT measures the length of the effective route, while the loop-size distribution measures the cyclic redundancy generated by the graph structure.

We applied the general theory to several classes of graphs. 
For complete graphs, we obtained explicit formulas for both the LERW MFPT and the erased-loop distribution. 
The results show that the typical loop-erased path length grows like the square root of the graph size. 
The erased-loop distribution has the same scaling and is closely related to the classical birthday problem. 
Thus, the complete graph provides a useful reference case for dense networks, where many alternative routes exist between two vertices.
For cycle graphs, the behavior is very different. 
Because the graph contains only one global cycle, the possible erased loops are strongly restricted. 
In large cycle graphs, almost all erased loops have size two, corresponding to simple backtracking events. 
This contrast with the complete graph case shows that LERW observables are sensitive to the topology of the graph. 
They can therefore distinguish between dense graphs with many possible loops and sparse cyclic graphs where the loop structure is much more constrained.

We also considered motif-based networks. 
In these graphs, the local structure of the motifs determines the global loop-erasure behavior. 
The MFPT can be decomposed along the sequence of motifs connecting two vertices, and the erased-loop distribution can be written as a mixture of the local distributions inside the motifs. 
This gives a clear interpretation: local building blocks control the effective random-walk behavior of the whole network. 
The numerical simulations again confirm the theoretical predictions.
A particularly important case is given by graphs built from cycle motifs. 
When the motifs are connected in a tree-like way, the path from one vertex to another crosses a unique sequence of motifs. 
In this situation, the global LERW MFPT can be computed by summing local contributions inside the crossed motifs. 
Equivalently, it can be obtained by summing the potential drops. 
This is a key simplification because the global problem is reduced to smaller computations on elementary cycle blocks. 
Therefore, graphs with cycle motifs form a special class where the theory is not only exact, but also computationally efficient and easy to interpret.

Overall, these results show that loop-erased random walks provide a meaningful way to study exploration on finite graphs. 
They give an intermediate point of view between shortest paths and ordinary random walks. 
Shortest paths ignore the random exploration process, while ordinary random walks keep all repeated visits. 
Loop-erased random walks keep only the effective route, and the erased-loop distribution adds complementary information by describing the loops that were removed. 
Together, these quantities provide a richer description of graph exploration and of the role played by cycles in the network.
From the perspective of exploration in unknown environments, a loop-erased
trajectory may also serve as a recorded return route: after locating a target,
an agent capable of remembering visited positions can traverse the retained
self-avoiding route in reverse rather than initiate a new random search.
This needs not be a shortest route and is not, by itself, an implemented
robotic navigation algorithm. Testing such a strategy with finite memory,
localization errors, and changes in the environment would be an interesting
future extension.
Several perspectives follow naturally from this work. 
A first direction is to improve the computational part of the method. 
Although the transfer-current formula gives exact weights, the enumeration of simple paths may become expensive on large dense graphs. 
Future work could use importance sampling, path truncation, or more efficient graph-search algorithms to reduce the computational cost.
Finally, the motif-based procedure is an efficient method for studying large networks built from repeated local structures   and could lead to possible applications in network design.
By choosing suitable motifs and connections between them, one could control the typical effective path length or the size of loops produced during exploration. 
This could be useful for designing networks where one wants to balance efficient navigation with controlled redundancy.
In conclusion, this work shows that the combination of loop-erased random walks, spanning trees, and transfer-current methods provides a powerful framework for studying finite graphs. 
It gives exact formulas, clear interpretations, and strong agreement with simulations. 
The special case of cycle-motif graphs is especially important, because it shows that the global LERW MFPT can be obtained by summing local potential contributions. 
This opens several directions for future work, both theoretical and applied.

 \section*{Acknowledgments}
Computational resources have been provided by the Consortium des Équipements de Calcul Intensif (CÉCI), funded by the Fonds de la Recherche Scientifique de Belgique (F.R.S.-FNRS) under Grant No. 2.5020.11 and by the Walloon Region. K.F. research is funded by ARC project "EMOTIONS".
\bibliography{article}
 
 \begin{appendices}
 \section{Absolute errors for complete graphs}
\label{app:complete-graph-errors}
\begin{figure}[htbp]
    \centering
    \includegraphics[width=\linewidth]{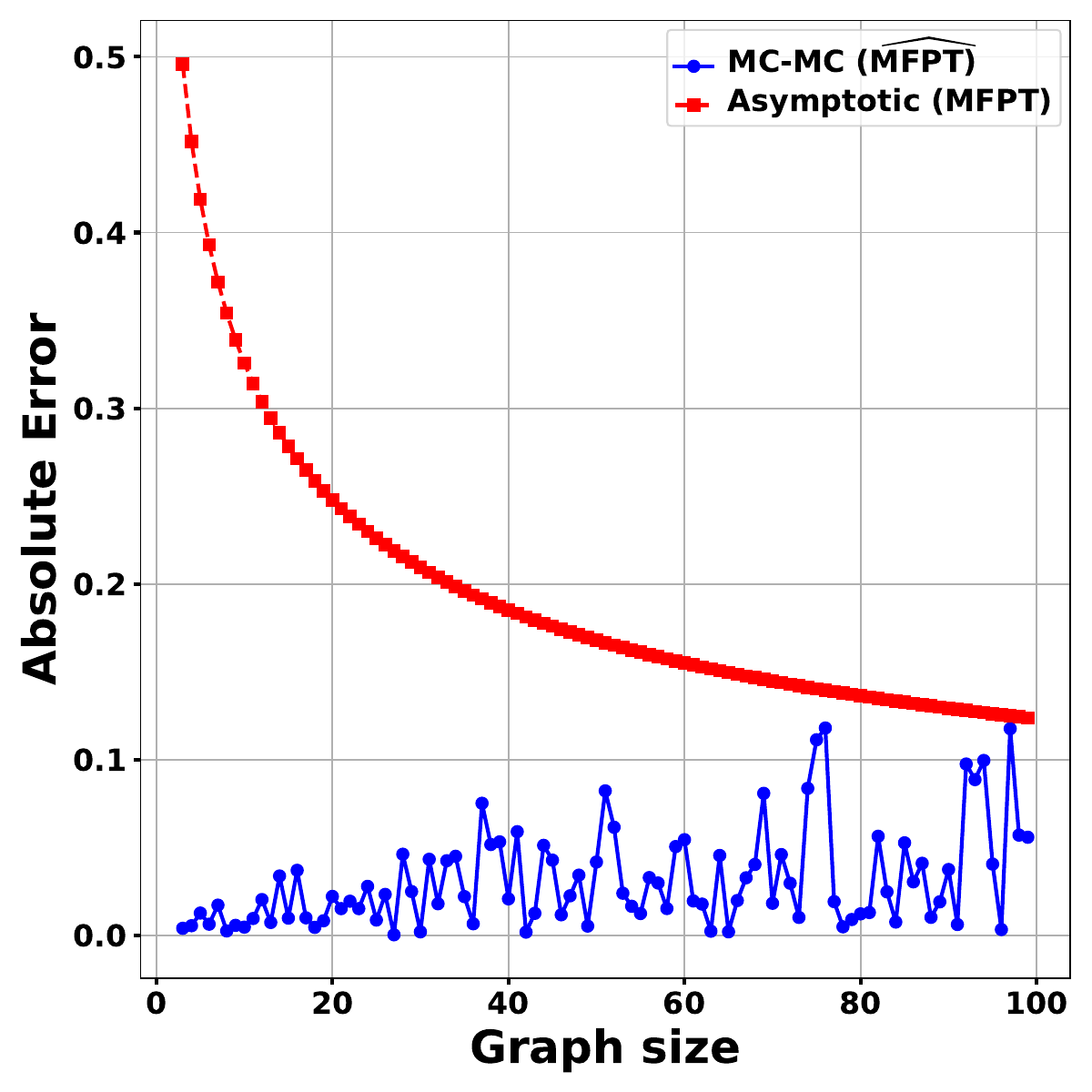}
    \caption{Absolute errors of the Monte Carlo estimate and the asymptotic approximation of Eq.~\eqref{SEC3EQ4}, both measured relative to the exact MFPT of Eq.~\eqref{SEC3EQ2}, as functions of the complete-graph size $n$.}
    \label{fig:comp-absolute-errors}
\end{figure}

Figure~\ref{fig:comp-absolute-errors} compares the absolute deviations of the Monte Carlo estimate and the asymptotic approximation from the exact MFPT on complete graphs.

\section{Algorithms}
In this appendix we schematically present the algorithms introduced in the main text.

\newcounter{appendixalgo}[section]
\renewcommand{\theappendixalgo}{\thesection.\arabic{appendixalgo}}
\begin{tcolorbox}[breakable,
  colback=beigeback,
  colframe=blueheader,
  title=\textbf{Algorithm C.1: Construction of the factor graph $T$ from motifs (cycle, clique) based graph $G$},
  fonttitle=\bfseries,
  coltitle=black,
  boxrule=0.8pt,
  arc=2mm,
  left=2mm,right=2mm,top=1mm,bottom=1mm
]\refstepcounter{appendixalgo}\label{ALGO2}

\textbf{Input.} An undirected graph $G=(V,E)$ whose biconnected components are cliques and/or simple cycles.

\vspace{0.15cm}

\textbf{Output.} A bipartite factor graph $T$ (blocks $\leftrightarrow$ vertices).

\vspace{0.15cm}

\textbf{1. Initialization.} Create an empty graph $T$.  
For each vertex $v\in V$, add a node $v$ to $T$ with attribute \texttt{kind = vertex}.

\vspace{0.15cm}

\textbf{2. Biconnected decomposition.} Compute the set of biconnected components
\[
\mathcal{B} = \{B_1,\dots,B_m\}
\]
of $G$.

\vspace{0.15cm}

\textbf{3. Block type identification.} For each block $B_i\in\mathcal{B}$:
\begin{itemize}[leftmargin=*,nosep]
  \item If $|B_i|=2$, treat $B_i$ as a clique $K_2$.
  \item Otherwise, test whether the induced subgraph $G[B_i]$ is a clique.
  \item Otherwise, test whether $G[B_i]$ is a simple cycle $C_k$ and reconstruct a cyclic order.
  \item If none applies, raise an error (unsupported block type).
\end{itemize}

\vspace{0.15cm}

\textbf{4. Block node construction.} Create a block node $b_i$ in $T$ with attributes $\texttt{kind = block}$, $\texttt{block\_type}\in\{\text{clique},\text{cycle}\}$, $\texttt{size}=|B_i|$. If $b_i$ is a cycle, store its cyclic ordering \texttt{cycle\_order}.

\vspace{0.15cm}

\textbf{5. Block--vertex connections.} For each vertex $v\in B_i$, add an edge $(b_i,v)$ to $T$.

\vspace{0.15cm}

\textbf{6. Structural check.} If $G$ is connected, verify that $T$ is a tree; otherwise raise an error.

\end{tcolorbox}

\begin{tcolorbox}[breakable,
  colback=beigeback,
  colframe=blueheader,
  title=\textbf{Algorithm C.2: Construction of the base graph $G$ from the factor graph $T$},
  fonttitle=\bfseries,
  coltitle=black,
  boxrule=0.8pt,
  arc=2mm,
  left=2mm,right=2mm,top=1mm,bottom=1mm
]\refstepcounter{appendixalgo}\label{ALGO3}

\textbf{Input.} A bipartite factor graph $T$ (assumed to be a tree), with nodes of type \texttt{vertex} and \texttt{block}.

\vspace{0.15cm}

\textbf{Output.} A graph $G$.

\vspace{0.15cm}

\textbf{1. Initialization.} Create an empty graph $G$.  
Add to $G$ all nodes of $T$ whose attribute is \texttt{kind = vertex}.

\vspace{0.15cm}

\textbf{2. Block traversal.} For each block node $b$ in $T$:
\begin{itemize}[leftmargin=*,nosep]
  \item Read its type \texttt{block\_type}.
  \item Extract the set of neighboring vertices $V_b$.
\end{itemize}

\vspace{0.15cm}

\textbf{3. Edge insertion.}
\begin{itemize}[leftmargin=*,nosep]
  \item If $b$ is a clique, add all possible edges between vertices in $V_b$.
  \item If $b$ is a cycle and a cyclic order is available, add consecutive cycle edges.
  \item Otherwise, add cycle edges in an arbitrary order.
\end{itemize}

\vspace{0.15cm}

\textbf{4. Output.} Return the reconstructed graph $G$.

\end{tcolorbox}

\begin{tcolorbox}[breakable,
  colback=beigeback,
  colframe=blueheader,
  title=\textbf{Algorithm C.3: General Factor Graph construction and related utilities},
  fonttitle=\bfseries,
  coltitle=black,
  boxrule=0.8pt,
  arc=2mm,
  left=2mm,right=2mm,top=1mm,bottom=1mm
]\refstepcounter{appendixalgo}\label{ALGO4}

\textbf{Input.} An undirected graph $G=(V,E)$ and a block prefix.

\vspace{0.15cm}

\textbf{Output.} (1) A factor graph $T$, (2) block-induced subgraphs, (3) articulation points inferred from $T$.

\vspace{0.15cm}

\textbf{1. Initialization of $T$.} Create an empty graph $T$.  
For each vertex $v\in V$, add $v$ to $T$ with attribute \texttt{kind = vertex}.

\vspace{0.15cm}

\textbf{2. Biconnected decomposition.} Compute the set of biconnected components
$\mathcal{B}=\{B_1,\dots,B_m\}$ of $G$.

\vspace{0.15cm}

\textbf{3. Block node construction.} For each block $B_i\in\mathcal{B}$:
\begin{itemize}[leftmargin=*,nosep]
  \item Create a block node $b_i$ in $T$ with attributes
  \[
  \texttt{kind = block},\quad \texttt{size}=|B_i|,\quad \texttt{nodes}=B_i.
  \]
  \item Connect $b_i$ to all vertices $v\in B_i$.
\end{itemize}

\textbf{4. Block subgraph extraction.} Given a block node $b$, return the induced subgraph $G[\texttt{nodes}(b)]$.

\vspace{0.15cm}

\textbf{5. Attachment point identification.} Define the set of articulation points as
$\mathcal{A}$, the set of all $v\in V$ adjacent to at least two block nodes in $T$.  
Return $\mathcal{A}$.
\end{tcolorbox}

%\end{appendices}
\section{ERROR BOUNDS AND CONVERGENCE OF THE MFPT ESTIMATOR}
\label{APP_MFPT_ERROR}

In this appendix, we establish a deterministic error bound
for the estimator introduced in Algorithm~\ref{alg:1}
and prove its almost-sure convergence to the exact
mean loop-erased path length.

Let $G=(V,E)$ be a finite, connected, simple, undirected
graph, and let $u,v\in V$ be two distinct vertices.
We denote by $\Gamma(u,v)$ the set of all simple paths
connecting $u$ to $v$.

For each $\gamma\in\Gamma(u,v)$, let
\begin{equation}
p_\gamma
=
\det\left(Y_{E(\gamma),E(\gamma)}\right)
\end{equation}
be its probability under the loop-erased random-walk
measure, as established in Eq.~\eqref{SEC1EQ3}.

The exact mean loop-erased path length is
\begin{equation}
\mu_{u,v}
=
\mathrm{MFPT}^{\mathrm{LERW}}(u,v)
=
\sum_{\gamma\in\Gamma(u,v)}
p_\gamma |\gamma|_E.
\label{eq:app_exact_mfpt}
\end{equation}

Since the set of simple paths forms a partition of
the possible loop-erased trajectories, their
probabilities satisfy
\begin{equation}
\sum_{\gamma\in\Gamma(u,v)}p_\gamma=1.
\label{eq:app_total_mass}
\end{equation}

Suppose that $N$ independent random-walk
realizations are generated and chronologically
loop-erased. Let
$\mathcal{G}_N\subseteq\Gamma(u,v)$
denote the set of distinct paths discovered
after these $N$ realizations.

We define the captured probability mass as
\begin{equation}
M_N
=
\sum_{\gamma\in\mathcal{G}_N}p_\gamma,
\qquad 0<M_N\leq1.
\label{eq:app_captured_mass}
\end{equation}

The estimator returned by Algorithm~\ref{alg:1}
can then be written as
\begin{equation}
\widehat{\mu}_N
=
\frac{1}{M_N}
\sum_{\gamma\in\mathcal{G}_N}
p_\gamma|\gamma|_E.
\label{eq:app_estimator}
\end{equation}

Unlike a standard Monte Carlo estimator,
$\widehat{\mu}_N$ is not obtained by averaging
the lengths of the simulated trajectories.
Instead, the simulations identify candidate
paths, whose probability weights are subsequently
computed exactly using the transfer-current kernel.

\begin{proposition}[Deterministic error bound]
\label{prop:app_error_bound}

Let $d_G(u,v)$ denote the shortest-path distance
between $u$ and $v$ in $G$.

For every nonempty set of discovered paths
$\mathcal{G}_N$, the estimator
$\widehat{\mu}_N$ satisfies
\begin{equation}
\left|
\mu_{u,v}-\widehat{\mu}_N
\right|
\leq
(1-M_N)
\left(
|V|-1-d_G(u,v)
\right).
\label{eq:app_error_bound}
\end{equation}
\end{proposition}

\begin{proof}
We decompose the exact expectation into
the contributions of discovered and
undiscovered paths:
\begin{align}
\mu_{u,v}
&=
\sum_{\gamma\in\mathcal{G}_N}
p_\gamma|\gamma|_E
\nonumber\\
&\quad+
\sum_{\gamma\in\Gamma(u,v)
\setminus\mathcal{G}_N}
p_\gamma|\gamma|_E.
\label{eq:app_decomposition}
\end{align}

By definition of $\widehat{\mu}_N$,
the first contribution equals
\[
\sum_{\gamma\in\mathcal{G}_N}
p_\gamma|\gamma|_E
=
M_N\widehat{\mu}_N.
\]

If $M_N<1$, define the conditional mean
length of the undiscovered paths by
\begin{equation}
\mu_N^{\mathrm{miss}}
=
\frac{1}{1-M_N}
\sum_{\gamma\in\Gamma(u,v)
\setminus\mathcal{G}_N}
p_\gamma|\gamma|_E.
\label{eq:app_missing_mean}
\end{equation}

Consequently, Eq.~\eqref{eq:app_decomposition}
becomes
\begin{equation}
\mu_{u,v}
=
M_N\widehat{\mu}_N
+
(1-M_N)\mu_N^{\mathrm{miss}}.
\label{eq:app_convex_combination}
\end{equation}

Subtracting $\widehat{\mu}_N$ from both
sides gives
\begin{equation}
\mu_{u,v}-\widehat{\mu}_N
=
(1-M_N)
\left(
\mu_N^{\mathrm{miss}}-\widehat{\mu}_N
\right).
\label{eq:app_error_identity}
\end{equation}

Since every loop-erased trajectory is
a simple path connecting $u$ to $v$,
its length satisfies
\begin{equation}
d_G(u,v)
\leq
|\gamma|_E
\leq
|V|-1.
\label{eq:app_length_bounds}
\end{equation}

The lower bound follows from the definition
of the shortest-path distance, whereas the
upper bound follows from the fact that a
simple path visits each vertex at most once.

Both $\widehat{\mu}_N$ and
$\mu_N^{\mathrm{miss}}$ are convex
combinations of simple-path lengths.
Therefore,
\begin{align}
d_G(u,v)
&\leq\widehat{\mu}_N\leq |V|-1,\\
d_G(u,v)
&\leq\mu_N^{\mathrm{miss}}\leq |V|-1.
\end{align}

It follows that
\begin{equation}
\left|
\mu_N^{\mathrm{miss}}-\widehat{\mu}_N
\right|
\leq
|V|-1-d_G(u,v).
\end{equation}

Combining this inequality with
Eq.~\eqref{eq:app_error_identity}
yields
\[
\left|
\mu_{u,v}-\widehat{\mu}_N
\right|
\leq
(1-M_N)
\left(
|V|-1-d_G(u,v)
\right),
\]
which proves the result.

If $M_N=1$, all the probability mass
has been captured and
$\widehat{\mu}_N=\mu_{u,v}$,
so the inequality also holds.
\end{proof}

\begin{corollary}[Bounds on the exact MFPT]
\label{cor:app_mfpt_interval}

Under the assumptions of
Proposition~\ref{prop:app_error_bound},
the exact mean loop-erased path length satisfies
\begin{align}
&M_N\widehat{\mu}_N
+(1-M_N)d_G(u,v)
\leq\mu_{u,v}
\nonumber\\
&\qquad\leq
M_N\widehat{\mu}_N
+(1-M_N)(|V|-1).
\label{eq:app_mfpt_interval}
\end{align}
\end{corollary}

\begin{proof}
Equation~\eqref{eq:app_convex_combination}
expresses the exact expectation as
\[
\mu_{u,v}
=
M_N\widehat{\mu}_N
+
(1-M_N)\mu_N^{\mathrm{miss}}.
\]

Using
\[
d_G(u,v)
\leq
\mu_N^{\mathrm{miss}}
\leq
|V|-1
\]
immediately gives the claimed inequalities.
The case $M_N=1$ follows by equality.
\end{proof}

These bounds are deterministic and do not rely
on a confidence level or an asymptotic normal
approximation. They provide an interval
containing the exact MFPT once the captured
probability mass $M_N$ is known.

Moreover, the interval width is
\begin{equation}
(1-M_N)
\left(|V|-1-d_G(u,v)\right),
\end{equation}
which decreases as the captured probability
mass approaches one.

\begin{proposition}[Almost-sure convergence]
\label{prop:app_consistency}

Assume that the $N$ loop-erased trajectories
used to construct $\mathcal{G}_N$ are generated
independently according to the loop-erased
random-walk distribution.

Then
\begin{equation}
M_N\xrightarrow[N\to\infty]{\mathrm{a.s.}}1,
\end{equation}
and
\begin{equation}
\widehat{\mu}_N
\xrightarrow[N\to\infty]{\mathrm{a.s.}}
\mu_{u,v}.
\label{eq:app_consistency}
\end{equation}
\end{proposition}

\begin{proof}
Since $G$ is finite, the set
$\Gamma(u,v)$ of simple paths
between $u$ and $v$ is finite.

Furthermore, every simple path
$\gamma\in\Gamma(u,v)$ has strictly
positive probability $p_\gamma>0$
under the loop-erased random-walk
measure.

Indeed, the original random walk
can follow the prescribed simple path
without revisiting any vertex, an event
having positive probability on a finite
connected graph.

For a given path $\gamma$, the probability
that it is not observed in any of the first
$N$ independent realizations is
\begin{equation}
\mathbb{P}
\left(
\gamma\notin\mathcal{G}_N
\right)
=
(1-p_\gamma)^N.
\end{equation}

Since $p_\gamma>0$,
\[
(1-p_\gamma)^N\longrightarrow0.
\]

Because $\Gamma(u,v)$ is finite,
the probability that at least one
simple path remains undiscovered
satisfies the union bound
\begin{equation}
\mathbb{P}
\left(
\mathcal{G}_N\neq\Gamma(u,v)
\right)
\leq
\sum_{\gamma\in\Gamma(u,v)}
(1-p_\gamma)^N
\longrightarrow0.
\end{equation}

In addition, the sets
$\mathcal{G}_N$ are increasing with $N$.
Each path is eventually discovered almost
surely, and therefore all paths are
eventually discovered almost surely.

Consequently,
\[
M_N
=
\sum_{\gamma\in\mathcal{G}_N}p_\gamma
\xrightarrow{\mathrm{a.s.}}1.
\]

Applying the deterministic bound
\eqref{eq:app_error_bound}, we obtain
\[
0\leq
\left|
\mu_{u,v}-\widehat{\mu}_N
\right|
\leq
(1-M_N)
\left(|V|-1-d_G(u,v)\right).
\]

Since the right-hand side converges
to zero almost surely, the estimator
converges almost surely to the exact
mean loop-erased path length.
\end{proof}

\section{PROOF OF THEOREM \ref{SEC2TH1}}
\label{APPA}

Following Dhar and Dhar \cite{DharDhar1997ErasedLoops}, we relate the formation
of erased loops to cycles created by adding an
edge to a random spanning tree.
Let $T$ be a uniformly chosen spanning tree of
$G$, and let $e$ be an edge selected uniformly
from $E$, independently of $T$.
For $e=\{u,v\}$, define
\begin{equation}
L_T(e)=d_T(u,v)+1,
\end{equation}
where $d_T(u,v)$ is the length of the unique
path connecting $u$ and $v$ in $T$.
If $e\notin T$, adding $e$ to $T$ creates a
unique cycle of size $L_T(e)$.
If $e\in T$, we adopt the convention that
$L_T(e)=2$, corresponding to immediate
backtracking along the same edge.
We denote by $P_{\mathrm{st}}(\ell)$ the
probability that the preceding spanning-tree
experiment produces a cycle of size $\ell$.
We now explain the relation between this
probability and the probability $P(\ell)$
that a step of the random walk creates and
erases a loop of size $\ell$.
Consider first a fixed simple cycle
$C$ of size $\ell>2$, whose edges are
\begin{equation}
E(C)=\{e_1,e_2,\ldots,e_\ell\}.
\end{equation}
For $j=1,\ldots,\ell$, let
$P_{\mathrm{st}}(C,e_j)$ be the probability
that adding the marked edge $e_j$ to a
spanning tree creates precisely the cycle $C$.
Since $T$ is a spanning tree, this event is
equivalent to
\begin{equation}
E(C)\setminus\{e_j\}\subseteq E(T),
\qquad e_j\notin E(T).
\end{equation}
The probability $P_{\mathrm{st}}(C,e_j)$
does not depend on $j$. Indeed, every
spanning unicycle whose unique cycle is $C$
can be transformed into a spanning tree by
removing any one of the $\ell$ edges of $C$.
This gives a bijection between the possible
marked-edge choices.
Therefore,
\begin{equation}
P_{\mathrm{st}}(C,e_j)
=
\frac{P_{\mathrm{st}}(C)}{\ell},
\label{eq:app_marked_cycle}
\end{equation}
where
\begin{equation}
P_{\mathrm{st}}(C)
=
\sum_{j=1}^{\ell}
P_{\mathrm{st}}(C,e_j).
\end{equation}
At time $N$, the loop-erased random walk can
be represented by the directed path from its
starting vertex $o$ to the current position
in the associated last-exit tree.
Adding an edge to this tree may create a
cycle. However, not every such cycle is
necessarily erased from the loop-erased path.
For a cycle $C$ closed by the marked edge
$e_j$, an additional condition is required:
the directed path from $o$ must enter the
cycle at the appropriate vertex, so that
the new step closes a loop on the current
loop-erased path.
Let $q_j(C)$ denote the corresponding
conditional probability in the Dhar--Dhar \cite{DharDhar1997ErasedLoops}
last-exit-tree construction.
The key observation is that,
for a fixed cycle $C$, the relevant entry
events associated with the possible marked
closing edges partition the possibilities.
Thus,
\begin{equation}
\sum_{j=1}^{\ell}q_j(C)=1.
\label{eq:app_partition}
\end{equation}
Under the assumed stationary
correspondence, the probability that the
marked edge $e_j$ creates and erases the
cycle $C$ is
\begin{equation}
P(C,e_j)
=
P_{\mathrm{st}}(C,e_j)\,q_j(C).
\end{equation}
Summing over the possible closing edges and
using Eq.~\eqref{eq:app_marked_cycle}, we obtain
\begin{align}
P(C)
&=
\sum_{j=1}^{\ell}P(C,e_j)\\
&=
\sum_{j=1}^{\ell}
P_{\mathrm{st}}(C,e_j)q_j(C)\\
&=
\frac{P_{\mathrm{st}}(C)}{\ell}
\sum_{j=1}^{\ell}q_j(C)\\
&=
\frac{P_{\mathrm{st}}(C)}{\ell}.
\end{align}
Finally, summing over all simple cycles of
size $\ell>2$, we obtain
\begin{equation}
P(\ell)
=
\frac{P_{\mathrm{st}}(\ell)}{\ell}.
\label{eq:app_dhar}
\end{equation}
For $\ell=2$, a loop corresponds to
immediate backtracking along an edge.
This case is treated
using the two orientations of the same edge,
giving the same relation, subject to the
assumed stationary correspondence
\begin{equation}
P(2)=\frac{P_{\mathrm{st}}(2)}{2}.
\end{equation}

\section{PROOF OF THEOREM \ref{SEC3TH1}}
To prove theorem~\ref{SEC3TH1} we need the following Theorem~\cite{Moon1970CountingLabelledTrees}.

\begin{theorem}\label{APPBTH1}
Let $F_n$ be a fixed forest on the labeled vertex set $[n]$ with $\ell=\ell(F_n)\ge 1$
connected components. Let the component sizes be $j_1,\dots,j_\ell$ (so
$\sum_{i=1}^{\ell} j_i=n$), and define
\[
p=p(F_n):=\prod_{i=1}^{\ell} j_i.
\]
Let $T(F_n)$ denote the number of labelled trees on $[n]$ that contain $F_n$ as a
spanning subgraph (i.e.\ that contain every edge of $F_n$). Then
\begin{equation}\label{APPBEQ1}
T(F_n)= p\, n^{\ell-2}.
\end{equation}
\end{theorem}

\begin{proof}\label{APPB}
\textbf{Part A: exact distribution \eqref{SEC3EQ1}.}
Fix $u$ and $v$. We count how many labelled trees have $S_n=k$. Count the possible $u$--$v$ paths with $k$ vertices. Such a path consists of $u$ and $v$ plus $k-2$ distinct intermediate vertices chosen from the remaining $n-2$ labels. Choose the intermediate vertex set in $\binom{n-2}{k-2}$ ways, and order them along
the path in $(k-2)!$ ways. Hence the number of distinct labelled paths
\[
u=w_0 - w_1 - \cdots - w_{k-2} - w_{k-1}=v
\]
with $k$ vertices equals
\[
\binom{n-2}{k-2}(k-2)! = \frac{(n-2)!}{(n-k)!}.
\]
For a fixed such path, count trees containing it.
Fix one concrete labelled $k$-vertex path $P$ from $u$ to $v$.
Consider the forest $F_n$ consisting of this path component (size $k$)
together with the remaining $n-k$ isolated vertices (each of size $1$).
This forest has
\[
\ell = 1+(n-k)=n-k+1
\]
components, and
\[
p = k\cdot 1^{n-k} = k.
\]
By The theorem \ref{APPBTH1}, the number of trees containing this forest (hence containing the path)
is
\[
T(F_n) = p\,n^{\ell-2} = k\,n^{(n-k+1)-2}=k\,n^{n-k-1}.
\]
 Multiply and normalize.
Therefore the total number of labelled trees with $S_n=k$ is
\[
\frac{(n-2)!}{(n-k)!}\cdot k\,n^{n-k-1}.
\]
Since there are $n^{n-2}$ labelled trees in total, we obtain
\[
\mathbb{P}\{S_n=k\}
=\frac{\frac{(n-2)!}{(n-k)!}\,k\,n^{n-k-1}}{n^{n-2}}
=\frac{k\,(n-2)!}{(n-k)!\,n^{k-1}},
\]
which is the first form of \eqref{SEC3EQ1}. The second form follows from
\begin{equation*}
    (n)_k=\dfrac{n!}{(n-k)!}=n(n-1)\dfrac{(n-2)!}{(n-k)!}.
\end{equation*}
\noindent\textbf{Part B: Mean scaling}
We will use the product form,
\[
\frac{(n)_k}{n^k}=\prod_{j=0}^{k-1}\left(1-\frac{j}{n}\right).
\]
For $0<t<1$, we have the inequality 
\begin{equation*}
 e^{-\frac{t}{1-t}}   <1-t< e^{-t}.
\end{equation*}
This means 
\begin{align*}
 &e^{-\frac{j/n}{1-j/n}}   <1-j/n< e^{-j/n}\\\Longrightarrow& e^{-\frac{j}{n-j}}   <1-j/n< e^{-j/n}\\
\Longrightarrow & \prod_{j=0}^{k-1} e^{-\frac{j}{n-j}}   <\prod_{j=0}^{k-1} 1-j/n< \prod_{j=0}^{k-1} e^{-j/n}\\
\Longrightarrow & e^{-\sum_{j=0}^{k-1}\frac{j}{n-j}}   <\prod_{j=0}^{k-1} 1-j/n<  e^{-\sum_{j=0}^{k-1}j/n}\\
\Longrightarrow & e^{-\sum_{j=0}^{k-1}\frac{j}{n-k}}   <\prod_{j=0}^{k-1} 1-j/n<  e^{-\sum_{j=0}^{k-1}j/n}\\
\Longrightarrow & e^{-\frac{k.(k-1)}{2(n-k)}}   <\prod_{j=0}^{k-1} 1-j/n<  e^{-\frac{k.(k-1)}{2n}}\\
\Longrightarrow & \frac{k}{n-1}e^{-\frac{k.(k-1)}{2(n-k)}}   <\mathbb{P}\{S_n=k\}<\frac{k}{n-1}  e^{-\frac{k.(k-1)}{2n}}\\
\Longrightarrow & \frac{k}{n-1}e^{-\frac{k^{2}}{2(n-k)}}   <\mathbb{P}\{S_n=k\}<\frac{k}{n-1}  e^{-\frac{k.(k-1)}{2n}}.\\
\end{align*}
We have ($\frac{k}{n-1}\sim_{\infty} \frac{k}{n}$). Thus,
we obtain the local approximation 
\[
\mathbb{P}\{S_n=k\}\sim_{\infty} \frac{k}{n}\exp\!\left(-\frac{k^2}{2n}\right).
\]
Now let $x>0$ and set $k= x\sqrt{n}$. The tail probability becomes
a Riemann sum
\begin{align*}
\mathbb{P}\left\{\frac{S_n}{\sqrt{n}}\ge x\right\}
&=\sum_{k\ge x\sqrt{n}} \mathbb{P}\{S_n=k\}
\\\;\longrightarrow\;
\int_x^\infty t\,e^{-t^2/2}\,dt
&=e^{-x^2/2}.    
\end{align*}
The limiting distribution is the Rayleigh law with density
$f(t)=t e^{-t^2/2}$ on $t>0$. Its first  moment is
\[
\mathbb{E}[R]=\int_0^\infty t\cdot t e^{-t^2/2}dt=\sqrt{\frac{\pi}{2}}.
\]
Because the tails are exponentially small, the family $S_n/\sqrt{n}$ is uniformly
integrable, and convergence in distribution implies convergence of the first two
moments
\[
\frac{\mathbb{E}[S_n]}{\sqrt{n}}\to \mathbb{E}[R]=\sqrt{\frac{\pi}{2}}.
\]
\end{proof}
\end{appendices}

\end{document}